\pdfoutput=1

\documentclass{siamltex}

\usepackage{thumbpdf}
\usepackage[%
  colorlinks%
  ,bookmarks={false}%
  ,pdftitle={}%
  ,pdfsubject={}%
  ,pdfkeywords={}%
  ,pdfauthor={Vladyslav Gapyak <vladyslav.gapyak@h-da.de>}%
  ]{hyperref}

\usepackage[numbers]{natbib}
\usepackage{a4wide}
\usepackage[a4paper, left=2cm, right=2cm, top=3cm, bottom=3cm]{geometry}
\usepackage{amsmath}
\usepackage{amssymb}
\usepackage{mathpazo}
\usepackage{graphicx}
\usepackage{enumerate}
\usepackage[utf8,latin1]{inputenc}
\usepackage[T1]{fontenc} 

\usepackage{mathtools}
\usepackage{xcolor}
\usepackage{caption}
\usepackage{subcaption}
\usepackage{tabularx}
\usepackage{multirow}
\usepackage{booktabs}
\usepackage{array}
\usepackage[percent]{overpic}
\usepackage{bm}
\usepackage{algorithm}
\usepackage{algpseudocode}

\newtheorem{rmk}{Remark}

\DeclareMathOperator{\divergence}{\text{div}}
\DeclareMathOperator{\prox}{\text{prox}}
\DeclareMathOperator{\sign}{\text{sign}}
\DeclareMathOperator{\Lip}{\text{Lip}}

\newcommand{\rom}[1]{%
	\textup{\uppercase\expandafter{\romannumeral#1}}%
}
\DeclareMathOperator{\trace}{trace}

\title{Variational Model-Based Reconstruction Techniques for Multi-Patch Data in Magnetic Particle Imaging.}

\author{Vladyslav Gapyak\footnotemark[1] \and Thomas M\"{a}rz\footnotemark[1]\
\and Andreas Weinmann\footnotemark[1]\ }

\begin{document}
\maketitle
\renewcommand{\thefootnote}{\fnsymbol{footnote}}
\footnotetext[1]{Department of Mathematics and Natural Sciences, Hochschule Darmstadt, Sch\"{o}fferstra{\ss}e 3, 64295, Darmstadt, Germany. Emails: vladyslav.gapyak@h-da.de, thomas.maerz@h-da.de, andreas.weinmann@h-da.de.}
\renewcommand{\thefootnote}{\arabic{footnote}}

\begin{abstract}
Magnetic Particle Imaging is an emerging imaging modality through which it is possible to detect tracers containing superparamagnetic nanoparticles. The exposure of the particles to dynamic magnetic fields generates a non-linear response that is used to locate the particles and produce an image of their distribution. The bounding box that can be covered by a single scan curve depends on the strength of the gradients of the magnetic fields applied, which is limited due to the risk of causing peripheral nerve stimulation (PNS) in the patients. To address this issue, multiple scans are performed. The scan data must be merged together to produce reconstructions of larger regions of interest. In this paper we propose a mathematical framework which can deal with rather general multi-patching scenarios including rigid transformations of the field of view (FoV), the specimen and of the scanner. We show the flexibility of this framework in a variety of different scanning scenarios. Moreover, we describe an iterative reconstruction algorithm that yields a reconstruction of the target distribution by minimizing a convex functional which includes positivity constraints and sparsity enforcing priors. We show its convergence to a minimizer and perform numerical experiments on simulated data.
\end{abstract}


\section{Introduction}\label{sec:introduction}

Magnetic Particle Imaging (MPI) is a promising new imaging modality introduced in 2005 by Gleich and Weizenecker in their seminal work~\cite{GleichWeizenecker2005}. The aim of MPI is to reconstruct the spatial distribution of superparamagnetic nanoparticles using the voltage induced by the non-linear response of the particles when exposed to a magnetic field. The paramagnetic particles are injected into the specimen to be imaged. The specimen is inserted into the MPI scanner, where both a static and a magnetic field are applied. In the field free point setup (contrasting with the field free line setup), these magnetic fields are such that there exists a point where they cancel out (a field free point, FFP) and the variation of the dynamic field moves the FFP in the scanning area. Near the FFP the response of the superparamagnetic particles is non-linear and close to the FFP the particles produce a voltage that is received by coils in the scanner. This is the signal used to reconstruct the original distribution of the paramagnetic particles which allows to visualize the interior of the specimen. Many important medical applications have been already thought of. For example, MPI can be used for blood flow imaging and cancer detection~\cite{BuzugKnopp2012} or stem cell imaging~\cite{zheng2013quantitative}. The recently developed human-sized MPI scanner could lead to MPI-assisted stroke detection and monitoring~\cite{Graeser2019humanbrain}. Moreover, MPI offers advantages over other established imaging modalities: for example, MPI offers a high spatial resolution and does not use radiation (this is opposed to methods such as PET~\cite{PET} and SPECT~\cite{SPECT}).

In the two and three dimensional setups it is common to employ system-matrix-based approaches, which acquire a system matrix in the following manner: for each voxel position a delta probe is scanned and yields a column of the system matrix~\cite{knopp2011prediction,WeizeneckerBorgertGleich2007,Rahmer_etal2012,Lampe_etal2012}. Once the system matrix $S$ is obtained, if the scan of the specimen is collected in a data vector $y$, one has to solve the linear system $S\rho =y$ for $\rho$ (the particle distribution). The spatial distribution of the particles can be for example obtained via regularized inversion of the linear system above~\cite{Knopp_etal2010ec}, with classical iterative reconstruction approaches~\cite{storath2016edge}, with novel iterative reconstruction approaches based on Plug-and Play~\cite{askin2022pnp}, or deep-learning-driven methods for reconstruction~\cite{gungor2023deqmpi}, by using the Deep Image Prior~\cite{DIP,dittmer2020deep}, and for the improvement of the spatial resolution of the reconstructions~\cite{Shang_2022}. The acquisition of a system matrix however is an extremely time-consuming procedure, therefore research has been pushed into the development of both acceleration methods for the calibration procedure (system matrix acquisition) and model-based techniques. Possible accelerated system matrix calibration include: compressed sensing calibration, which exploit the high compressibility of the system matrix~\cite{KnoppWeber2013}, its spatial symmetry~\cite{Weber_2015_symmetries} and the combination of both properties~\cite{Weber2015_symm_compressed}; super-resolution of system matrices acquired on down-sampled grids using deep learning~\cite{gungor2021superresolution,Schrank2022superresolution,gungor2022tranSMS,Yin2023dipSM}; extrapolation of system matrices outside the FoV~\cite{schefflerboundary2022,scheffler2023extrapolation}.
The development of model-based techniques is another important direction of research, because it allows to avoid or significantly reduce the lengthy calibration procedure. One possibility is to produce a model-based system matrix by simulating the scans of the delta distribution in each voxel using for example the Langevin model~\cite{KnoppSattel_etal2010,KnoppBiederer_etal2010} and employ any of the regularized inversion methods developed for the measured-system matrices. In the framework of model-based techniques, work has been done also to avoid the system matrix formulation of the MPI imaging task~\cite{Rahemeretal2009,Schomberg2010,GoodwillConolly2010,GoodwillConolly2011,Gruettner_etal2013,bringout2020new}. It is precisely in the framework of model-based methods not relying on the system-matrix-based approach that this paper contributes.

The application of static and dynamic magnetic fields lays at the core of the MPI technology and consequently, studies on the requirements regarding the basic safety of patients have been conducted in comparison with the currently available International Electrotechnical Commission standards for safety in MRI~\cite{schamelsafety2015,Saritas2013-jn}. According to these studies, for the patient to be safe during a scan, peripheral nerve stimulation (PNS) and tissue heating shall be avoided. In particular, the selection field employed to saturate the nanoparticles but those in the FFP is strong, for example $7Tm^{-1}$~\cite{Saritas2013-jn}, and the drive field use to steer the FFP in the FoV in the scanners currently in use employs amplitudes of $0.1-40mT$ with frequency between $5$ and $25kHz$~\cite{Saritas2013-jn}. It is known that for sinusoidal electric currents, the value of $100kHz$ is the crossover frequency from nerve stimulation to thermal heating and above this frequency, tissue heating becomes a safety concern~\cite{Saritas2013-jn}. The frequency and the amplitude of the drive field determine the size of the FoV: indeed, the partial field of view pFoV$[m]$ (FoV covered during half of the period) is related to the peak amplitude $B_{\mathrm{peak}}[T]$ of the sinusoidal drive field and $\mu_0 G[T/m]$ the selection field strength, where $\mu_0$ is the vacuum permeability, by~\cite{Saritas2013-jn} $pFoV = \frac{2B_{\mathrm{peak}}}{\mu_0G}$. Consequently, the safety requirements pose a bound on the strength of the gradients of the magnetic fields that can be employed, limiting thereby the spatial coverage of a scan -- also called the Field of View (FoV) -- to about 10-30 mm~\cite{szwargulski2018efficient,Saritas2013-jn}.
Therefore, to scan bigger regions, so-called multi-patch scan sequences are utilized~\cite{szwargulski2018efficient,gdaniec2020suppression}. To perform measurements on multiple patches so far there are two ways: the first possibility is to use an additional focus field that can shift the FoV~\cite{Weizenecker2010FastMD}. The second way is to physically move the object through the FoV~\cite{Szwargulskimovingtable2018}. In particular, the multi-patching is usually performed as follows: the FoV or the object are positioned in their initial location and a scan of the region covered by the FoV is performed: this constitutes a patch; then, either the object or the FoV is moved and positioned in the next region of interest and another scan is performed. This procedure continues in an alternating fashion with step-wise reposition until the collection of scans/patches covers the specimen. In this paper, we will refer to this procedure as \emph{standard multi-patching} in juxtaposition to \emph{generalized multi-patching}, the label with which we denote the framework presented in Section~\ref{sec:mathematical:modeling} of this paper.

While multi-patching has been done in the measurement-based approach, the objective of this paper is to provide a mathematical framework that allows to incorporate multi-patching in the model-based approach introduced in~\cite{math10183278}. To enhance further the quality of the reconstructions we introduce positivity constraints and sparsity promoting priors in the energy functional to be minimized. Further, we prove the convergence of the employed iterative algorithm to a minimizer of said functional.

\paragraph{\textbf{Related Work}}

In the measurement-based approach it is essential to collect the system matrix which contains the response of the scanning system to delta distributions in the voxel positions. Each patch however, requires the acquisition of its own system matrix~\cite{Knopp_2015}, which increases the time needed for the calibration. To overcome this problem methods that reuse a single system matrix can lead to faster calibrations and reconstruction, but produce displacement artifacts~\cite{szwargulski2018efficient} and effort has been put to alleviate such artifacts~\cite{boberggeneralized2020,schefflerboundary2022}.

To circumvent the acquisition of the system matrix in the measurement-based approach, model-based approaches have been considered as well~\cite{Rahemeretal2009,KnoppSattel_etal2010,Schomberg2010,GoodwillConolly2010,GoodwillConolly2011,Gruettner_etal2013,bringout2020new}. Concerning the model-based approach, the first reconstruction formulae where given for the one dimensional setup using a relationship between the model and the Chebyshev polynomials~\cite{Rahemeretal2009}. The preliminary work~\cite{marz2016model} of the authors provided reconstruction formulae also for the 2D and 3D cases and a two-stage algorithm, that employed a local least square estimate of the MPI core operator in the first stage and a deconvolution in the second stage. This algorithm was made more flexible by estimating the MPI core operator in the first stage with a variational approach~\cite{marz2022amee}. The flexibility and advantages of the two-stage approach introduced with the variational first stage have been showed in the authors' work~\cite{math10183278}, in which new quality-enhancing techniques have been proposed. In particular, because the first stage of the variational approach is independent of the scanning trajectory chosen, one could merge scans coming from different trajectories (in general rigid body transformations of the standard trajectory) and reconstruct on a thus better spatially distributed sampling. For example, performing scans after rotating the specimen can help enhance the quality of the reconstruction, as show in our previous conference proceeding~\cite{math10183278}.

\paragraph{\textbf{Contributions}}

With this paper we contribute a mathematical framework for model-based image reconstruction from multi-patch MPI data.
Based on this mathematical framework, we extend an algorithm we introduced in~\cite{math10183278} so that it can handle (general) multi-patch data
and analyze the mathematical properties of the proposed algorithm.
In more detail, we make the following contributions:
\begin{itemize}
	\item[(i)] we provide a general mathematical framework for multi-patching that can deal with any rigid transformation of the specimen, the FoV and the scanner. In particular, we show that upon transformation of the data acquired, multi-patching can be seen mathematically as a scan along a transformed Lissajous trajectory, with fixed scanner and particle distribution, and that the data thus collected can be fed in the reconstruction algorithm we propose;
	\item[(ii)] We show that the multi-patching approach can be incorporated in our previous two-stage algorithm \cite{marz2016model,marz2022amee,math10183278}. Furthermore, we add positivity constraints and a sparsity promoting prior to the functional in the second stage, thus introducing a non-differentiable component to the minimization problem, which is treated with established methods. We prove convergence of the iterations to a minimizer of the employed functional. In particular, we prove that the differentiable part of the functional is Lipschitz continuous and provide a Lipschitz constant;
	\item[(iii)] we perform numerical experiments on simulated data in order to show the potential of the proposed framework and algorithm.
\end{itemize}

A preliminary idea to handle multi-patch data was presented as extended conference abstract in~\cite{gapyak2022icnaam}.

\paragraph{\textbf{Outline}}

The paper is structured as follows: in Section~\ref{sec:mathematical:modeling} we describe the mathematical framework for the scanning process, introduce the generalized Lissajous trajectories and show that this framework is indeed flexible and includes the current multi-patching approaches employed as a subcases. In Section~\ref{sec:rec:algorithm} we recall the theory behind the model-based approach and the relationships between the data and the MPI Core Operator, which provide the equations needed for the two-stage algorithm. Moreover, we describe in detail the algorithms employed and the discretizations performed for the numerical treatment of the experiments. We prove that the functional employed in the second stage satisfies sufficient conditions that guarantee the convergence of the iterative algorithm to a minimizer of the functional. In Section~\ref{sec:experiments} we describe the experimental setup, the parameter choices and their optimization. Finally we describe the experiments performed and summarize the results obtained. In Section~\ref{sec:conclusions} we draw our conclusions and discuss further possible research directions.

\section{A Mathematical Framework for Generalized Multi-Patching}\label{sec:mathematical:modeling}

\begin{figure}[t]
	\centering
	\begin{overpic}[width=0.5\textwidth]{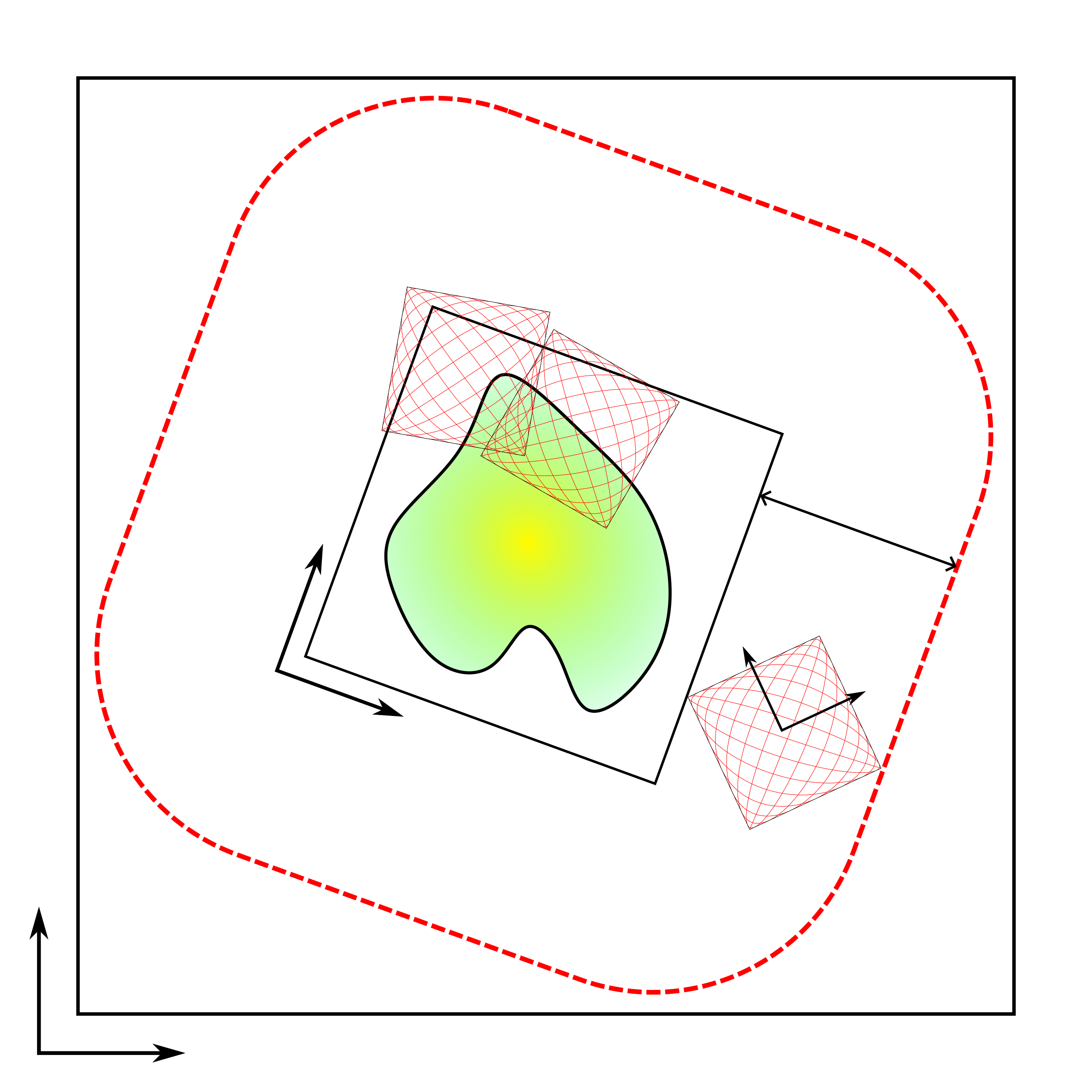}
		\put (48,49) {$\rho$}
		\put (48,28) {$\Omega$}
		\put (48,13) {$\Omega_{2d}$}
		\put (48,95) {$\Omega_S$}
		\put (78,53) {$2d$}
		\put (1,1)   {$\bm{x}_S$}
		\put (22,36) {$\bm{x}_{\Omega}$}
		\put (70,30) {$\bm{x}_{F}$}
	\end{overpic}
	\caption{Graphical representation of our generalized multi-patch framework in 2D. The physical volume inside the scanner is denoted with $\Omega_S$. The distribution $\rho$ shares the frame of reference $\bm{x}_\Omega$ of the region $\Omega$ to be reconstructed. The region $\Omega_{2d}$ is the region where the scan points can fall into and $2d$ is the diameter of the FoV. The motion of the FoV inside $\Omega_{2d}$ can be represented by the relative motion of the frame of reference $\bm{x}_{F}$ of the FoV and the frame of reference of the scanner $\bm{x}_S$. }
	\label{fig:model}
\end{figure}

In this section we describe a general mathematical framework for multi-patching in MPI, allowing for rigid transformations of the scanning trajectory as a way to overcome the limitations arising from the bound on the size of the FoV due to peripheral nerve stimulation~\cite{Saritas2013-jn}. First, we provide a mathematical description of the scanning area (or volume) inspired by the multi-patch approach using an additional focus field, which allows to move the FoV (region where the scan is performed) keeping the target particle distribution $\rho$ and the scanner fixed. Please refer to Fig.~\ref{fig:model} for a visualization of the quantities introduced in a schematic 2D scenario.

Let us denote with $\Omega_{S}\subset\mathbb{R}^n$, for $n\in\lbrace 1,2,3\rbrace$ the physical space inside the scanner, whose points have coordinates written in the frame of reference (f.o.r. from now on) $\bm{x}_S = (x_S^{(1)},\dots ,x_S^{(n)})^T$ of the scanner, and with $\Omega\subset\Omega_{S}$ the region containing the particles distribution $\rho\colon\Omega\mapsto\mathbf{R}^{+}$ to be scanned and reconstructed. Because we are dealing with a scanner which entirely contains the specimen, we will consider compactly supported distributions $\rho$, i.e., such that $\overline{\supp (\rho )}\subset \mathring{\Omega}$, where $\overline{\Omega }$ (resp. $\mathring{\Omega}$) is the closure (resp. the interior) of $\Omega$ with respect to the standard topology on $\mathbb{R}^n$ and $\supp (\rho )=\lbrace x\in\Omega\colon \rho (x)\neq 0\rbrace$.  We assume $\Omega$ to be a box, i.e.,  $\Omega = \prod_{i=1}^{n} [a_i , b_i ]$ where $a_i,b_i\in\mathbb{R}$ such that $a_i < b_i$ for all $i\in\lbrace 1,\dots ,n\rbrace$ with its own f.o.r. $\bm{x}_\Omega$. Usually the scans are performed along Lissajous trajectories~\cite{Knoppetal2009,knopp2020openmpidata} $r\colon [0,T]\to\mathbb{R}^n$, whose analytical form is:
\begin{equation}\label{eq:lissajous}
	r(t) = \bigl ( A_1\sin (2\pi m_1 t + \phi_1 ), \dots , A_n\sin (2\pi m_n t + \phi_n)\bigr )^T ,
\end{equation}
where for all $i\in\lbrace 1, \dots ,n\rbrace$, $m_i$ are the frequencies, $\phi_i$ are the phase shifts and $A_i$ the amplitudes. The centered field of view (FoV) is the subset $\Phi (0)=\prod_{i=1}^{n}[-A_i , A_i ]$ spanned by a scan along the curve $r$ in Eq.\eqref{eq:lissajous}. As discussed in the Section~\ref{sec:introduction}, the Lissajous amplitudes that are possible in a real scanner are not big enough to guarantee spatial coverage of the region $\Omega$, and multi-patching is necessary. In particular, in the standard multi-patching~\cite{szwargulski2018efficient}, one moves the FoV and performs different scans in different regions of $\Omega$. This procedure can be easily described in the following way: given an offset vector $b\in\mathbb{R}^n$, the FoV is the subset $\Phi (b )=b + \Phi (0)$ and multi-patching is the collection of scans in the regions $\Phi (b)$ for a collection of offset vectors $b$. In particular, the FoV has its own f.o.r. $\bm{x}_F$ and with respect to $\bm{x}_F$, points inside the FoV are always in the set $\Phi (0)$ and the scan trajectory has the analytic form expressed in Eq.~\eqref{eq:lissajous}. Inspired by other common scanning modalities, in particular by CT~\cite{Kalender1990SpiralVC} and MRI~\cite{Kruger2002ContinuouslyMT}, we introduce a general mathematical framework that allows to scan a specimen with an MPI scanner that can move the FoV and rotate it at the same time, both while scanning. The relative motion of the scanner, the region $\Omega$ and the FoV can be therefore expressed by relative motions between their f.o.r. $\bm{x}_S$,$\bm{x}_\Omega$ and $\bm{x}_F$ (see Fig.~\ref{fig:model}).

Consider now the sets of $\varepsilon$-distant points to $\Omega$ for some $\varepsilon >0$:
\begin{equation}
	\Omega_\varepsilon \coloneq \lbrace x\in \Omega_{S} \colon \text{dist}(x,\Omega)\leq \varepsilon\rbrace ,
\end{equation}
where $\text{dist}(x,\Omega) = \inf \lbrace \lVert x-y\rVert_2 \colon y\in\Omega\rbrace$.
Let $d \coloneq \sqrt{\sum_{i=1}^{n} A_i^2}$ be half of the diameter of the centered field of view $\Phi (0)$, then we assume the set inclusions $\Omega\subset\Omega_{2d}\subset\Omega_S$ and the following properties to hold:
\begin{itemize}
	\item $\Omega_S \setminus \Omega_{2d}$ cannot be scanned and will not contain any data;
	\item $\Omega_{2d}\setminus\Omega$ may have data, but it is not used for the reconstruction;
	\item $\Omega$ contains data and such data is used for the reconstruction.
\end{itemize}

Now, consider a time-dependent off-set vector $b\colon [0,T]\to\Omega_d$ and a function $Q\colon [0,T] \to \mathrm{SO}(n)$, i.e., $Q(t)$ is an orthogonal matrix with $\text{det}\bigl ( Q(t)\bigr ) = 1$ for every time $t\in [0,T]$.
These two functions can be used to define a trajectory in the space of rigid transformations:
\begin{align}
	[0,T] &\to \Omega_d \times \mathrm{SO}(n) \nonumber\\
	t &\mapsto   \bigl ( b(t),Q(t)\bigr )  .
\end{align}
In the 2D case, which is considered in this paper, $Q (t)$ can be defined as the rotation matrix:
\begin{equation}\label{eq:rot:matrix}
	Q_{\alpha}(t) = \begin{pmatrix}
		\cos \alpha (t) & -\sin\alpha (t)\\
		\sin\alpha (t)  & \cos\alpha (t)
	\end{pmatrix} ,
\end{equation}
for some angle function $\alpha : [0,T]\to \mathbb{R}$. In the 3D case, the rotation matrix $Q (t)$ can be constructed considering the vector function $\alpha (t) = (\phi (t)$, $\theta (t)$, $\psi (t))^T$ whose elements are the Euler angles: precession, nutation and rotation~\cite{goldsteinmechanics}. Because the transformation $Q(t)$ is determined by either a scalar or a 3-dimensional angle function $\alpha (t)$, we will write $Q_\alpha (t)$ in place of $Q(t)$ in what follows.
We introduce the generalized Lissajous trajectory of offset $b(t)$ and rotation angle $\alpha (t)$ as follows:
\begin{equation}\label{eq:gen:liss}
	\Lambda (t) = b(t)+Q_{\alpha} (t)\,  r(t)
\end{equation}
where $r(t)$ is a standard Lissajous trajectory expressed in $\bm{x}_F$ as in Eq.~\eqref{eq:lissajous}. The trajectory $\Lambda (t)$ is consequently the scanning trajectory in the f.o.r. $\bm{x}_S$ of the scanner. We point out that the term ``generalized Lissajous trajectory" is a term of convenience, because, even though in what follows the base scanning trajectory $r(t)$ is a Lissajous curve as in Eq.~\eqref{eq:lissajous}, a generalized trajectory $\Lambda (t)$ as in Eq.~\eqref{eq:gen:liss} can be formed with any choice of trajectory $r(t)$. This framework provides a way to describe different multi-patching scans and in particular helps describing the standard multi-patching approach as a generalized Lissajous trajectory in Eq.\eqref{eq:gen:liss}.
Further assumptions regarding $\rho$ are discussed in Section~\ref{sec:rec:algorithm}.

\paragraph{\textbf{Standard Multi-Patching in the Generalized Lissajous Setting}}\label{sec:gen:multipatching}

We now describe how multi-patching \cite{szwargulski2018efficient,gdaniec2020suppression} and the enhancement of the quality through merging of rotations introduced by the authors~\cite{math10183278} are particular subcases of the general framework introduced above. In particular, suppose that $T>0$ is the acquisition time of one scan and that $\Xi\in\mathbb{N}$ scans of a distribution $\rho$ are performed, each scan centered in some points $b_{\xi}\in\Omega_d$ with rotation angle $\alpha_{\xi}\in [0,2\pi )$ for $\xi \in\lbrace 0,\dots ,\Xi -1\rbrace$. Moreover, let $\tau >0$ be the time necessary to move the center of the FoV from $b_\xi$ to $b_{\xi +1}$. Then we consider $\tilde{T}\coloneq \Xi\cdot T+(\Xi -1)\tau$ the total time needed to perform the generalized multi-patch scan and moving in an alternating fashion $\Xi$ times. The intervals in which the scans are performed are denoted with $I(\xi ,T ,\tau ) = [\xi (T+\tau ),\xi (T+\tau )+T]$, with the set of times spent scanning will be denoted as $\Sigma \coloneq \bigcup_{\xi = 0}^{\Xi - 1} I (\xi , T ,\tau )$ while the time spent moving in which no scan is performed corresponds to the set $N \coloneq [0,\tilde{T}] \setminus \Sigma$. Then, it is sufficient to consider an angle function $\alpha\in C^1 ([0,\tilde{T} ], \mathbb{R})$ and an offset trajectory $b\in C^1 ([0,\tilde{T} ], \Omega_d )$ such that they are constant on the scanning intervals:
\begin{equation}\label{eq:piecewise:b:t}
	b \restriction_{I (\xi ,T, \tau )} \equiv b_\xi \quad \text{and}\quad \alpha\restriction_{I (\xi ,T, \tau )} \equiv \alpha_\xi \quad\text{for all }\xi\in\lbrace 0,\dots ,\Xi -1\rbrace.
\end{equation}
This is coherent with known interpretations of standard multi-patching~\cite{Szwargulskimovingtable2018}. Even though the framework provided considers a static distribution $\rho$ and a fixed scanner, it can deal with other setups as we shall see. Before proceeding however, we need to recall the relationship between the distribution $\rho$ and the data obtained by the scanning procedure according the model-based approach.

\paragraph{\textbf{Signal Encoding}}

We now briefly recall the principles behind the model-based approach for the field-free point (FFP) setup. The interested reader can find further details in the authors' previous paper~\cite{marz2016model} and in the book~\cite{BuzugKnopp2012}.
In Magnetic Particle Imaging the non-linearity of the magnetization response of superparamagnetic particles is exploited. Indeed, the applied magnetic field $H(x,t)\in\mathbb{R}^3$ is usually of the form $H(x,t)=Gx + H_D (t)$, i.e., the superposition of a static magnetic field $Gx$ with non-singular gradient matrix $G\in\mathbb{R}^{3\times 3}$ and the drive field $H_D (t)$, such that there is a point where the field $H$ vanishes and it is hence called the field-free point (FFP). If the location of the FFP in time is denoted with $r(t)$, it is uniquely defined by $H \bigl ( r(t)\bigr ) =0$, as $r(t)=-G^{-1}H_D (t)$ and conversely $H_D (t)= -Gr(t)$. In particular, this leads to a more convenient expression of the applied magnetic field: $H(x,t)=G(x-r(t))$. The response of the superparamagnetic particles is modeled according to the Langevin theory of superparamagnetism~\cite{chikazumi1978physics,jiles1998introduction}, that is, if $m$ is the magnetic moment of a single particle, $H_{\text{sat}}$ is a parameter depending on different physical quantities (e.g., particles' diameter, temperature) and $\rho (x)$ is the particle distribution, the magnetization response $M(x,t)\in\mathbb{R}^3$ is modeled by the following equation:
\begin{equation}\label{eq:magnetization:response}
	M(x,t) = \rho (x) m \mathcal{L}\biggl ( \frac{\lvert H(x,t)\rvert}{H_{\text{sat}}}\biggr ) \frac{H(x,t)}{\lvert H(x,t)\rvert}, \quad\text{where}\quad\mathcal{L}(\xi )=\coth (\xi )-\frac{1}{\xi}\quad\text{and}\quad H_{\text{sat}}=\frac{k_b T}{\mu_0 M_{\text{sat}}\frac{\pi}{6}d^3},
\end{equation}
where $k_b$ is the Boltzmann constant, $\mu_0$ the magnetic permeability in vacuum, $T$ the temperature of the particles, $M_{\text{sat}}$ the saturation magnetization of the particles and $d$ their diameter.
The principle behind MPI scanners~\cite{gleich2005original} is that the particles close to the FFP induce a voltage in the receiving coils of the scanner, and constitute the data collected from which the particle distribution is to be reconstructed. According to Faraday's law, this voltage is the negative time derivative of the magnetic flux, which is the spatial integral of the superposition of the applied magnetic field $H(x,t)$ and the magnetization response $M(x,t)$ of the particles in Eq.~\eqref{eq:magnetization:response}. However, only the magnetization response $M(x,t)$ depends on the target particle distribution $\rho$, and consequently contains information about it, meaning that one can subtract the signal of an empty scan and obtain a signal $s(t)\in\mathbb{R}^3$ modeled as follows:
\begin{equation}\label{eq:signal:model}
	s(t) = \mu_0 m \frac{d}{dt}\int_{\mathbb{R}^3}\rho (x) \mathcal{L}\biggl (\frac{\lvert G(r(t)-x)\rvert}{H_{\text{sat}}}\biggr ) \frac{G(r(t)-x)}{\lvert G(r(t)-x)\rvert}\, dx.
\end{equation}
Upon changes of variables and a nondimensionalization (see~\cite{marz2016model} for the details), the relationship in Eq.~\eqref{eq:signal:model} can be rewritten as follows:
\begin{equation}\label{eq:signal:dimensionless}
	s(t) = \frac{d}{dt}\int_{\mathbb{R}^3} \rho (x) \mathcal{L}\biggl (\frac{\lvert r(t)-x\rvert}{h}\biggr ) \frac{r(t)-x}{\lvert r(t)-x\rvert}\, dx ,\quad\text{where}\quad h=\frac{H_{\text{sat}}}{gL},
\end{equation}
where $h>0$ is a new, dimensionless resolution parameter which is obtained as a rescaling of $H_{\mathrm{sat}}$ by the strength $g$ of the gradient $G$ and $L$, the length of the FoV. Passing the derivative under the integral sign in Eq.~\eqref{eq:signal:dimensionless}, we obtain the following relationship between the signal $s(t)$, the scanning trajectory $r(t)$, its velocity $v(t)=\frac{d}{dt}r(t)$ and the MPI Core Operator $A_h [\rho ](r)\in\mathbb{R}^{3\times 3}$:
\begin{equation}\label{eq:MPI:core}
	s(t) = A_h [\rho ]\bigl (r(t)\bigr ) v(t)\quad \text{where}\quad A_h [\rho ] (r) = \int_{\mathbb{R}^3} \rho (x) \nabla_r \biggl (\mathcal{L}\biggl (\frac{\lvert r-x\rvert}{h}\biggr ) \frac{r-x}{\lvert r-x\rvert}\biggr )\, dx .
\end{equation}
The MPI Core Operator $A_h [\rho ]$ defined in Eq.~\eqref{eq:MPI:core} encodes the information needed to restore the distribution $\rho$. It is however not necessary to work with matrices, as it has been proved~\cite{marz2016model} that the trace of the MPI Core Operator contains all the information about the distribution. In particular, the following relationship holds true:
\begin{equation}\label{eq:trace:convolution}
	\kappa_h * \rho = \trace (A_h [\rho ])\quad\text{where}\quad \kappa_h (y) = \divergence \biggl ( \mathcal{L}\biggl (\frac{\lvert y\rvert}{h}\biggr ) \frac{y}{\lvert y\rvert}  \biggr )\quad\text{is the scalar-valued kernel defined for }y\in\mathbb{R}^3 .
\end{equation}
\begin{rmk}\label{rmk:h}
	We clarify now the meaning of the parameter $h$. Given the definitions of $H_{\mathrm{sat}}$ in Eq.~\eqref{eq:magnetization:response} and of $h$ in Eq.~\eqref{eq:signal:dimensionless}, and assuming that the temperature of the particles is $310 K$ -- temperature of the human body -- and considering a length of the FoV $L=2cm$, the value of $h$ for nanoparticles of diameter $d$ of size $20nm\leq d\leq 30nm$ is in the range $0.019\geq h\geq 0.006$~\cite{marz2016model,BuzugKnopp2012}. The value of $0.01$ in between is a reasonable choice of parameter that catches the magnitude of $h$ for a particle under reasonable conditions. We observe here that, because the Langevin function $\mathcal{L}$ is a sigmoid function, the lower the value of $h$, the better $\mathcal{L}$ approximates the sign function. In 1D this means that the convolution kernel $\kappa_h$ approximates the Dirac delta, but this property is not transferred to the general $n$-dimensional case ($n>1$), as the the convolution kernels $\kappa_h$ converges (in the sense of distributions) to the kernel $\kappa (y)=\frac{n-1}{\lvert y\rvert}$, with $y\in\mathbb{R}^{n}$~\cite[Corollary 1]{marz2016model}. Moreover, if $\alpha_h [\rho ](r)\coloneq \trace A_h [\rho ](r)$ is the trace of the MPI Core Operator and $\alpha [\rho ](r) \coloneq \lim_{h\to 0}\alpha_h [\rho ] (r)$ is the ideal limit, the relationship between the trace and its limit is given by the formula $\alpha_h [\rho ] = -\frac{\Delta\kappa_h}{8\pi}* \alpha [\rho ]$~\cite[Theorem 3.6]{marz2016model}. In particular, because $\kappa_h$ is analytic~\cite[Theorem 3.4]{marz2016model}, the trace $\alpha_h [\rho ]$ is a smoothed version of the idealized trace $\alpha [\rho ]$. This shows that a smaller $h$ is to be preferred, i.e., for the same scanning setup and environmental conditions, bigger particles are to be preferred. However, the bigger particles in the real world scenario increase the relaxation effects~\cite{Shasha2019relaxation} and need non-adiabatic assumptions of the model~\cite{Croft2012relaxation}. In this paper however, we consider adiabatic assumptions and no relaxation effect.
\end{rmk}
\vspace{\baselineskip}

As we will explain in Section~\ref{sec:rec:algorithm} and in view of Remark~\ref{remark:no:trajectory}, the data needed for a reconstruction is a set of triplets of the form $(s(t),r(t),v(t))$ or in the multi-patching setting $(s(t),\Lambda (t), \Lambda' (t))$ in the coordinate system of the reconstruction area $\Omega$, but the samples are collected by the scanner and consequently expressed in the frame of reference of the scanner. In the framework presented here however, these two frames of reference coincide and no transformation of the data is needed. This is not a limitation as other setups can be brought back to the setup that uses generalized Lissajous curves. We will discuss next how to perform the needed transformations of different setups into the generalized Lissajous framework.

\paragraph{\textbf{Other Setups}}

The problem of performing multiple scans is equivalent to collecting data $s(t)$ along the generalized Lissajous $\Lambda (t)$ defined in Eq.~\eqref{eq:gen:liss} with offset trajectory and angle function defined in Eq.~\eqref{eq:piecewise:b:t}, related to each other by the formula:
\begin{equation}\label{eq:gen:liss:data}
	s(t) = A_h [\rho ]\bigl (\Lambda (t) \bigr ) \Lambda ' (t) .
\end{equation}
For the reconstructions one considers as input data the set $\mathcal{S}$ in phase-space:
\begin{equation}\label{eq:data:points}
	\mathcal{S} = \bigl\lbrace \bigl ( s(t),\Lambda (t) , \Lambda ' (t) \bigr )\colon t\in \Sigma \bigr \rbrace .
\end{equation}
In general however, one has to pay attention to the frame of reference of each component of the elements in $\mathcal{S}$. Indeed, the data $s(t)$ is collected by the receiving coils of the scanner and are therefore in the frame of reference of the scanner itself but they must be in the frame of reference of the reconstruction region $\Omega$ for the reconstruction to take place. As mentioned, multi-patching can be achieved by either moving the FoV via a focus field or by physically moving the specimen inside the scanner. The relationships in Eq.~\eqref{eq:MPI:core} hold true only if the specimen is time independent, i.e., it does not vary in time inside the reconstruction area $\Omega$. This means that the frame of reference of $\Omega$ and $\rho$ coincide or equivalently, that there is no relative motion between $\Omega$ and $\rho$. In our framework we suppose that the frame of reference of the scanner and that of the region $\Omega$ (and hence of $\rho$) coincide, while the FoV is free to move and rotate. This is a mathematical simplification that does not exclude other possibilities, like moving the distribution through the FoV. Indeed, our first result in this paper is to show that our focus-field-inspired framework is equivalent to any other way of acquiring a multi-patching scan for time-independent distributions.
\vspace{\baselineskip}
\begin{lemma}\label{lemma:transf}
	Consider a time-independent particle distribution $\rho$ with frame of reference $\bm{x}_{\Omega}$ integral to that of $\Omega$. The following multi-patching approaches are equivalent up to transformation of the data acquired:
	\begin{itemize}
		\item[i)] both the scanner and $\Omega$ are fixed and the data is collected with rigid transformations of the FoV;
		\item[ii)] the scanner and the FoV are fixed, while $\Omega$ is moved using rigid transformations;
		\item[iii)] the region $\Omega$ is fixed while the scanner and the FoV move together with the same rigid transformation.
	\end{itemize}
\end{lemma}
\vspace{\baselineskip}
\begin{proof}
	In the following proof, let $\bm{x}_S$ denote the frame of reference (f.o.r.) of the scanner (and hence of the FoV) and $\bm{x}_{\Omega}$ the frame of reference of $\Omega$ (see Fig.~\ref{fig:model}). The data is collected by the scanner, therefore information about the scanning trajectory $\hat{r}(t)$, the velocities $\hat{v}(t)$ and the data collected $\hat{s}(t)$ will be denoted by a hat on top of the relative quantity to indicate that they are in the f.o.r. of the scanner $\bm{x}_S$. The data needed for the reconstruction must be in the f.o.r. $\bm{x}_{\Omega}$ of $\Omega$. All rigid transformations are compositions of a rotation represented by a rotation matrix $Q_{\alpha}(t)$ and a shift represented through an offset trajectory $b(t)$. Because all rigid transformations are invertible, to prove the equivalence of \emph{i)},\emph{ii)} and \emph{iii)} it is enough to provide the necessary rigid transformations of the data $(\hat{s}(t),\hat{r}(t),\hat{v}(t))$ collected by the scanner into the f.o.r. $\bm{x}_{\Omega}$.
	
	\emph{i)$\iff$ii)} Consider a fixed scanner and a moving region $\Omega$ which from the point of view of $\bm{x}_S$ moves rigidly: $\hat{\Omega} (t) = \hat{b}(t)+Q_{\alpha}(t)\Omega$. This means that we have the following transformation of coordinates:
	\begin{equation}
		\bm{x}_\Omega = \hat{b}(t)+Q_{\alpha}(t)\bm{x}_S .
	\end{equation}
	So, the Lissajous trajectory $\hat{r}(t)$ becomes a generalized Lissajous trajectory $\Lambda (t) = \hat{b}(t)+Q_{\alpha}(t)\hat{r}(t)$ in $\bm{x}_{\Omega}$ and the data $\hat{s}(t)$ must be transformed accordingly  $\hat{s}(t)\mapsto s(t) = \hat{b}(t)+Q_{\alpha}(t)\hat{s}(t)$.
	The velocities in $\bm{x}_{\Omega}$ are:
	\begin{equation}
		\Lambda '(t) = \hat{b}' (t) + \frac{d}{dt}Q_{\alpha}(t) \hat{r}(t)+Q_{\alpha}(t)\frac{d}{dt}\hat{r}(t).
	\end{equation}
	
	\emph{i)$\iff$iii)} Suppose now that $\Omega$ is fixed and the scanner (and hence the FoV) is moving w.r.t $\bm{x}_{\Omega}$ with rotation matrix $Q_{\alpha} (t)$ and offset vector $b(t)$, then the coordinate transformations have the following form:
	\begin{equation}\label{eq:direct:transf}
		\bm{x}_S = b(t)+Q_{\alpha}(t)\bm{x}_{\Omega},
	\end{equation}
	whose inverse transformation is:
	\begin{equation}\label{eq:inverse:transf}
		\bm{x}_{\Omega} = Q_{\alpha}^T (t)[\bm{x}_S -b(t)].
	\end{equation}
	The transformation in Eq.~\eqref{eq:inverse:transf} provides the formula to transform the data $s(t)=Q_{\alpha}^T (t)[\hat{s}(t)-b(t)]$ and the transformation of $\Lambda$ and $\Lambda '$ into $\hat{r}$ and $\hat{v}$ can be written upon differentiation in the following form:
	\begin{equation}\label{eq:prop1:matrix}
		\begin{pmatrix}
			\hat{r}(t) \\
			\hat{v}(t)
		\end{pmatrix}
		=
		\begin{pmatrix}
			Q_{\alpha}(t)  & \bm{0} \\
			\frac{d}{dt}Q_{\alpha}(t) & Q_{\alpha}(t)
		\end{pmatrix}
		\begin{pmatrix}
			\Lambda (t)  \\
			\Lambda ' (t)
		\end{pmatrix}
		+
		\begin{pmatrix}
			b (t)  \\
			b ' (t)
		\end{pmatrix}.
	\end{equation}
	Inversion of the system in Eq.~\eqref{eq:prop1:matrix} provides the transformation formula for the trajectory data written with respect to $\bm{x}_{S}$ into the data needed for the reconstruction:
	\begin{equation}\label{eq:lemma:iii}
		\begin{pmatrix}
			\Lambda (t)  \\
			\Lambda ' (t)
		\end{pmatrix}
		=
		\begin{pmatrix}
			Q_{\alpha}(t)  & \bm{0} \\
			\frac{d}{dt}Q_{\alpha}(t) & Q_{\alpha}(t)
		\end{pmatrix}^{-1}
		\biggl [\begin{pmatrix}
			\hat{r}(t) \\
			\hat{v}(t)
		\end{pmatrix}
		-
		\begin{pmatrix}
			b (t)  \\
			b ' (t)
		\end{pmatrix}
		\biggr ] .
	\end{equation}
	
\end{proof}

The framework however can be used also for any general multi-patching scheme that can use both a focus field moving the FoV, rigid motion of $\Omega$ and movements of the scanner, i.e., for multi-patching approaches in which all three frames of reference $\bm{x}_S$, $\bm{x}_\Omega$ and $\bm{x}_F$ are free to move:
\vspace{\baselineskip}
\begin{theorem}\label{theorem:transformations}
	Consider a time-independent particle distribution $\rho$ with frame of reference integral to that of $\Omega$, then any multi-patching scan acquired using rigid transformation of the FoV, the scanner and the region $\Omega$ can be rewritten in the form of the generalized multi-patching.
\end{theorem}
\vspace{\baselineskip}

\begin{proof}
	Consider the three frames of reference: $\bm{x}_S$ of the scanner, $\bm{x}_F$ of the FoV, $\bm{x}_{\Omega}$ of $\Omega$ and $\bm{x}_R$ an absolute f.o.r. of the room in which the scan is performed. Each frame of reference moves with respect to the room with a rigid motion, i.e., 
	\begin{equation}\label{eq:movements}
		\bm{x}_k = b_k (t)+Q_k (t) \bm{x}_R \quad\text{for all }k\in\lbrace \Omega ,F, S\rbrace ,
	\end{equation}
	for some offset vectors $b_F (t) , b_\Omega (t) , b_S (t)\in\mathbb{R}^n$ and rotation matrices $Q_F (t) , Q_\Omega (t) , Q_S (t)\in\text{SO}(n)$ as in Eq.~\eqref{eq:rot:matrix}. We will omit from now on the dependence on time for ease of notation. The inverse transformation of each f.o.r. into the coordinates of the room $\bm{x}_R$ read as follows:
	\begin{equation}\label{eq:into:room}
		\bm{x}_R = Q_k^T [\bm{x}_k -b_k ]\quad\text{for all }k\in\lbrace \Omega ,F, S\rbrace .
	\end{equation}
	The Lissajous trajectory $r=r(t)$ in Eq.~\eqref{eq:lissajous} is written in $\bm{x}_F$ and produces a generalized Lissajous trajectory $\Lambda_S = \Lambda_S (t)$ of the form:
	\begin{equation}\label{eq:lambda:s}
		\Lambda_S = b_S + Q_S [Q_F^T (r-b_F )] ,
	\end{equation}
	combining Eq.~\eqref{eq:into:room} and~\eqref{eq:movements}. Analogously, the relative motion of $\Omega$ in $\bm{x}_S$ is
	\begin{align}
		\bm{x}_S = b_S + Q_S [Q_\Omega^T (\bm{x}_\Omega -b_\Omega )] = \underbrace{b_S - Q_S Q_\Omega^T b_\Omega}_{b_*} + \underbrace{Q_S Q_\Omega^T}_{Q_*}\bm{x}_\Omega .\label{eq:omega:s}
	\end{align}
	With these considerations, we can imagine the scanner $\bm{x}_S$ and the FoV $\bm{x}_F$ to be fixed, but scanning along the generalized Lissajous trajectory $\Lambda_S$ in Eq.~\eqref{eq:lambda:s}, and $\bm{x}_\Omega$ moving with respect to $\bm{x}_S$ according to Eq.~\eqref{eq:omega:s}. Therefore, we can use \emph{iii)} in Lemma~\ref{lemma:transf} to transform the data collected by the scanner into $\bm{x}_\Omega$. Indeed, if $\hat{s}(t)$ is the data collected by the scanner, we obtain the transformed data $s(t)$ in $\bm{x}_{\Omega}$ by $s(t) = Q_*^T [\hat{s}(t)-b_* ]$, and the positions and velocities $\Lambda (t), \Lambda ' (t)$ in $\bm{x}_{\Omega}$ using Eq.~\eqref{eq:lemma:iii} with $\hat{r}(t) = \Lambda_S (t)$, $Q_\alpha (t) = Q_* (t)$, $b(t)=b_* (t)$ and their derivatives.
	
\end{proof}

Due to Theorem~\ref{theorem:transformations}, any multi-patching modality can be mathematically seen from the point of view of our generalized multi-patching framework upon transformation. We formulate this as a Corollary:
\vspace{\baselineskip}
\begin{corollary}
	Given any MPI scan involving rigid motions of the FoV due to a focus field, of the scanner and of the specimen, it is always possible reformulate it as a scan with fixed scanner and fixed specimen, but general scanning trajectory.
\end{corollary}
\vspace{\baselineskip}
\begin{rmk}
	We would like to stress once more that the data points describing the scan (position, velocities and the signal collected by the scanner) must be expressed in the f.o.r. $\bm{x}_\Omega$ of $\Omega$ in order to operate the first stage of the algorithm in Section~\ref{sec:rec:algorithm}. In the proof of Theorem~\ref{theorem:transformations} we provide transformation formulas to convert the data $\hat{s}(t),\hat{r}(t)$ and $\hat{v}(t)$ collected into $\bm{x}_\Omega$. For this reason we will assume that the data collected is already converted into $\bm{x}_\Omega$ when talking about the first stage of the algorithm and we will denote with $s(t)$,$r(t)$ and $v(t)$ the signal, the scan trajectory and the velocities in $\bm{x}_\Omega$.
\end{rmk}

\section{Reconstruction Algorithm for the Multi-Patching Setting}\label{sec:rec:algorithm}

\subsection{Reconstruction Formulae}
In Section~\ref{sec:mathematical:modeling} we have seen that the relationship between the data $s(t)$, the information about the points at which the signal is collected $r(t)$, $v(t)$ and the distribution $\rho$ is mediated by the MPI Core Operator. In particular, from  Eq.~\eqref{eq:MPI:core} and~\eqref{eq:trace:convolution} it is possible to derive the following two-stage reconstruction method~\cite{marz2016model,marz2022amee,math10183278,maerz2022icnaam,gapyak2022icnaam}: the first stage consists of reconstructing the MPI Core Operator from the data, imposing the relationship stated in Eq.~\eqref{eq:signal:dimensionless}; in the second stage the particle distribution $\rho$ is reconstructed solving the deconvolution problem in Eq.~\eqref{eq:trace:convolution} with input data the trace $u=\trace (A_h [\rho ])$ of the MPI Core Operator $A_h [\rho ]$ reconstructed in the first stage.

\paragraph{\textbf{Stage 1: Reconstruction of the MPI Core Operator}}

In the MPI setting, the data $s(t)$ is collected along a discrete time series in $[0,T]$ for some $T>0$, let $t_k = \frac{T}{L}k$ for $k\in\lbrace 1, \dots ,L\rbrace$ and $L$ the number of samples. At each time $t_k$, the position of the FFP point will be denoted with $r_k = r(t_k )$ and its velocity as $v_k = v (t_k )$, leading to the discretized version of the relationship expressed by Eq.~\eqref{eq:signal:dimensionless}:
\begin{equation}\label{eq:core:discrete}
	s_k = A_h [\rho ](r_k )v_k\quad\text{for }k\in\lbrace 1,\dots ,L\rbrace .
\end{equation}
From Eq.~\eqref{eq:core:discrete} the independence of the reconstruction method from the scanning trajectories is even more clear since the given data points are samples where the causality and connectivity is not important in phase-space of the form Eq.~\ref{eq:data:points}. In the first stage of the reconstruction algorithm, we aim at using the data coming from a scan, i.e., from the scanning points $s_k$ at positions $r_k$ and velocities $v_k$ to retrieve the MPI Core Operator which mediates the relationship (see Eq.~\ref{eq:MPI:core}) between the distribution $\rho$ and the data points. The output of the first stage is the trace of the MPI Core Operator, which is the distribution $\rho$ convolved with a known kernel (see Eq.~\eqref{eq:second:stage}). The trace of the MPI Core Operator will be the target of a deconvolution to retrieve $\rho$ in the second stage.
To reconstruct the MPI Core Operator the authors have proposed two methods so far: the LLSq (Local Least Square method)~\cite{marz2016model} and a variational approach introduced by the authors in~\cite{marz2022amee}. The latter has the advantage of not imposing any restriction on the grid used to discretize the region $\Omega$. In this paper we will use the variational approach: if $N$ is the number of voxels in the chosen grid, one obtains an estimation of the MPI Core Operator minimizing the functional $J[A ]$ defined as follows:
\begin{align}\label{eq:first:step}
	A &= \arg \min\limits_{\hat{A}} \; J\left[\hat{A}\right], & \quad\text{where}\quad
	J\left[\hat{A}\right] &= \underbrace{\frac{\lambda}{N} \lVert D \hat{A}\rVert_{2}^2}_{\eqcolon R\left[\hat{A}\right]}
	\; + \; \underbrace{\frac{1}{L} \sum\limits_{k=1}^L \left| s_k - I\left[\hat{A}\right](r_k) v_k \right|^2}_{\eqcolon F\left[\hat{A}\right]}
	,
\end{align}
where $D$ is the matrix obtained by discretizing the gradient using first differences, $\lambda >0$ is the regularization parameter controlling the strength of the regularization and $I[A]$ denotes the usage of an interpolation scheme to evaluate $I[A]$ in the specific point $r_k$.
A more precise description of the discretization procedure can be found in the numerical treatment and reconstruction algorithms further in this section.
\vspace{\baselineskip}
\begin{rmk}\label{remark:no:trajectory}
	We point out here a fundamental fact for multi-patching stamming from Eq.~\eqref{eq:core:discrete}: the MPI Core Operator is a matrix-valued operator dependent only on the single points of evaluation $r_k =r(t_k )$ and a direction vector $v_k = v(t_k )$ and not on the specific properties of the chosen trajectory $t\mapsto r(t)$ itself. Indeed, any set of triplets $(s_k ,r_k ,v_k )$ can function as an input for the first stage and precisely this property allows us to merge different scans into a denser or enlarged input data set for the reconstruction. Moreover, the only data coming from the scanner is the signal $s_k$, while the positions and velocities $r_k$, $v_k$ are computed from the analytic expression of the trajectory. In particular, if any type of distortion of the chosen trajectory is known, potential correcting transformations of the data points $r_k$, $v_k$ can be applied.
\end{rmk}

\vspace{\baselineskip}
\paragraph{\textbf{Stage 2: Regularized Deconvolution}}

Once an approximation of the MPI Core Operator is obtained using Eq. \eqref{eq:first:step}, its trace $u=\trace (A )$ serves as input data for the deconvolution needed to retrieve $\rho$ as prescribed by Eq.~\eqref{eq:trace:convolution}. Deconvolution is an ill-posed problem, and regularization methods are necessary~\cite{bertero2021introduction,hansen2010discrete}.

In regularized deconvolution the spatial particle distribution $\rho$ is obtained by minimizing the energy functional introduced by the authors in~\cite{marz2022amee}:
\begin{equation}\label{eq:second:stage}
	\rho = \arg\min_{\hat{\rho}} E[\hat{\rho}] , \quad\text{where}\quad E[\hat{\rho}]=\lVert\kappa_h * \hat{\rho} - u\rVert_2^2 + \mu R[\hat{\rho}].
\end{equation}
In Eq.~\eqref{eq:second:stage}, the term $R[\rho ]$ is a regularizer weighted by the parameter $\mu >0$. In particular, in this paper we assume that the distribution is an element of the Sobolev space $\rho\in W^{1,1}(\Omega )\subset \text{BV}(\Omega)$ \footnote{Here $W^{k,p}(\Omega )$ is the Sobolev space of differentiation order $k$ and integration order $p$ and $\text{BV}(\Omega )$ denotes the space of functions with bounded variations~\cite{ambrosio2000fbv}} (and in particular with finite Total Variation) and consider as regularizer the so-called TV Smooth~\cite{David_Dobson_1996,vogel1996iterative,vogel1998fast} approximation of the Total Variation (TV), defined as follows
\begin{equation}\label{eq:tv:smooth:reg}
	R_\delta [\rho ] = \int_{\Omega}\sqrt{\delta + \lvert\nabla\rho (x)\rvert^2}\, dx
\end{equation}
approaching $\lvert\rho\rvert_{\text{TV}}$ for $\delta\to 0$. Of course, many other approximations of the TV norm can be considered. However, the choice of this specific approximation of the Total Variation will be motivated in the discussion on the numerical treatment of the reconstruction later.

\paragraph{\textbf{Non-negative Fused Lasso for Regularized Deconvolution}}

The minimization of the energy functional in Eq.~\eqref{eq:second:stage} leads to good results~\cite{math10183278} with both the TV Smooth approximation and a more standard Tikhonov regularizer. However, it can be extended to incorporate further \emph{a priori} MPI-specific knowledge. Because $\rho$ represents a particle distribution, it is a natural and common practice to impose positivity constraints, which also lead to better reconstruction quality~\cite{Weizenecker_etal2009,KnoppBiederer_etal2010}. Moreover, it has been shown~\cite{storath2016edge} that, because the tracer is usually concentrated in specific spots inside the field of view and in order to enhance the quality of the images, better preserve edges and suppress noise artifacts, one could also enforce sparsity-promoting priors. In particular, in this paper we extend the energy proposed in~\cite{math10183278} considering the following constrained minimization problem:
\begin{equation}\label{eq:pos:sparsity:constraint}
	\rho = \arg\min_{\hat{\rho}\geq 0}\bigl\lbrace\lVert \kappa_h * \hat{\rho} - u \rVert_2^2 + \mu R[\hat{\rho}] +\beta\lVert\hat{\rho}\rVert_1\bigr\rbrace .
\end{equation}
The positivity constraint $\hat{\rho}\geq 0$ is easily dealt with introducing the indicator function
\begin{equation*}
	\iota_+ (\hat{\rho}) = \begin{cases}
		0, & \text{if }\hat{\rho}\geq 0 \\
		+\infty & \text{otherwise}.
	\end{cases}
\end{equation*}
Indeed, with the definition of the indicator function above, it is equivalent to consider the following unconstrained minimization problem:
\begin{align}\label{eq:second:pos:sparsity}
	\rho = \arg\min_{\hat{\rho}} \bigl\lbrace \underbrace{\lVert \kappa_h *\hat{\rho}-u\rVert_2^2 + \mu R[\hat{\rho}]}_{\eqcolon F[\hat{\rho}]} + \underbrace{\beta\lVert\hat{\rho}\rVert_1}_{\eqcolon G_1(\hat{\rho})} + \, \iota_+ (\hat{\rho}) \bigr\rbrace .
\end{align}
In general the minimization problem in Eq.~\eqref{eq:second:pos:sparsity} involves terms of different nature: some terms are differentiable and some are not, but are convex. The differentiability of the term $F$ is however not true in general, as it depends on the choice of the regularizer: for instance, if one considers the near-isotropic discretizations of the TV norm~\cite{storath2016edge,storath2015joint,storath2017fast}, the only differentiable term in Eq.~\eqref{eq:second:pos:sparsity} is $\lVert \kappa_h *\hat{\rho}-u\rVert_2^2$. We note that various different regularizers $R[\hat{\rho}]$ can be employed in Eq.\eqref{eq:second:pos:sparsity}, for example the Tikhonov regularizer used in~\cite{math10183278}. The choice of the TV type regularizer is motivated by its edge-preserving properties.
However, if $R[\hat{\rho}]$ is the TV Smooth regularizer in Eq.~\eqref{eq:tv:smooth:reg}, the term $F[\hat{\rho}]$ in Eq.~\eqref{eq:second:pos:sparsity} is differentiable. With different regularizers other minimization algorithms are also possible and they potentially can be learning based. Such algorithm have already been employed in MPI dealing with the system-matrix approach, for example using Plug-and-Play~\cite{askin2022pnp}.

Furthermore, the sparsity-promoting term $G_1 (\hat{\rho})$ and the indicator function $\iota_+$ are both convex and can be dealt with using the theory of convex optimization~\cite{bauschkeconvex,ekelandconvex}, in particular with the use of proximal mappings and proximal algorithms~\cite{OPT-003}. We recall here that if $f\colon\mathbb{R}^n\to\mathbb{R}\cup\lbrace+\infty\rbrace$ is a closed proper convex function, the associated proximal mapping is well defined as:
\begin{align}\label{eq:proximal}
	\prox_{f}\colon  \mathbb{R}^n & \to \mathbb{R}^n \nonumber\\
	v & \mapsto  \prox_f (v)\coloneq\arg\min_x \bigl ( f(x)+\frac{1}{2}\lVert x-v\rVert_2^2\bigr ) .
\end{align}

To deal with the different nature of the terms in the energy functional in Eq.~\eqref{eq:second:pos:sparsity}, in this paper we use the generalized forward-backwards splitting scheme by Raguet et al.~\cite{doi:10.1137/120872802}. More explicitly, the iterative splitting scheme applied is:
\begin{align}\label{eq:alg:lasso}
	z_1^{(k+1)} & = z_1^{(k)} + \lambda_k \biggl [ \prox_{2\gamma G_1} \bigl ( 2\rho^{(k)} - z_1^{(k)} - \gamma\nabla F(\rho^{(k)})\bigr ) - \rho^{(k)} \biggr ]	\\
	z_2^{(k+1)} & = z_2^{(k)} + \lambda_k \biggl [ \prox_{2\gamma \iota_+} \bigl ( 2\rho^{(k)} - z_2^{(k)} - \gamma\nabla F(\rho^{(k)})\bigr ) - \rho^{(k)} \biggr ]	\notag\\
	\rho^{(k+1)} & = \frac{1}{2}(z_1^{(k+1)}+z_2^{(k+1)}) .\notag
\end{align}
Convergence conditions for the algorithm depend on $\Lip (\nabla F)$, which in turn depends on the discretization chosen. Therefore convergence conditions and a convergence proof will be provided in Section~\ref{sec:convergence:alg} after discussing in more detail the numerical treatment of the two stages.

\paragraph{\textbf{Numerical Treatment and Reconstruction Algorithms}}

We now provide a detailed discussion on the discretization techniques chosen, pseudocodes for each stage of the algorithm.

From now on, we will consider the 2D scenario, as it is similar to the 3D case but has simpler notation. We consider the case where $\Omega = [a,b]\times [c,d]$, discretized by an $N_x\times N_y$ Cartesian grid. Let $h_x = \frac{b-a}{N_x}$ be the discretization step size along the $x$ axis and $h_y = \frac{d-c}{N_y}$ the step size along the $y$ axis. The pixels are then identified with their central points $(x_i ,y_j )\in\Omega$ with coordinates:
\begin{equation}\label{eq:pixel:indeces}
	x_i = a + \biggl (i+\frac{1}{2}\biggr ) h_x ,\quad y_j = c + \biggl (j+\frac{1}{2}\biggr ) h_y ,
\end{equation}
for $i\in\lbrace 0,\dots , N_x -1\rbrace$ and $j\in\lbrace 0,\dots , N_y -1\rbrace$.

\textbf{First Stage}. Once the discretization of the set $\Omega$ is chosen, we aim at reconstructing on it the MPI Core Operator, i.e., we are looking for a tensor $A\in\mathbb{R}^{N_x\times N_y\times 2\times 2}$, $\bigl ( A_{i,j}^{p,q}\bigr )$ where $i,j$ refer to the pixel and $p,q$ are the indices of the matrix, minimizing the functional in Eq.~\eqref{eq:first:step}, i.e., that solve the corresponding gradient system $\nabla J [A] = 0$.

The regularizing term $\lVert D\hat{A}\rVert$ in Eq.~\eqref{eq:first:step} will be represented as follows:
\begin{equation}\label{eq:first:discrete:reg}
	\lVert D\hat{A}\rVert_2^2 = \sum_{i=0}^{N_x -1}\sum_{j=0}^{N_y -1}\biggl ( \biggl \lVert \frac{\hat{A}_{i+1,j} - \hat{A}_{i,j}}{h_x}\biggr \rVert_F^2 + \biggl \lVert \frac{\hat{A}_{i,j+1} - \hat{A}_{i,j}}{h_y}\biggr \rVert_F^2\biggr )
\end{equation}
As computed in~\cite{marz2022amee}, the partial derivative of the data fidelity term $F$ in Eq.~\eqref{eq:first:step} with respect to $A$ are:
\begin{align}\label{eq:first:reg:partials}
	\frac{\partial F[A]}{\partial A_{i,j}^{p,q}} &= -\frac{2}{L}\sum_{k=1}^L
	\left\langle s_k - I[A](r_k )v_k , I\bigl [ \delta_{i,j}^{p,q}\bigr ] (r_k )v_k\right\rangle ,
\end{align}
where $\delta_{i,j}^{p,q}$ refers to the Kronecker symbol, and the partial derivative of the regularizing term in Eq.~\eqref{eq:first:discrete:reg} are:
\begin{align}
	\frac{\partial R[A]}{\partial A_{i,j}^{p,q}} &=-\frac{2\lambda}{N_x N_y}\left( \frac{A_{i-1,j}^{p,q}-2A_{i,j}^{p,q}+A_{i+1,j}^{p,q}}{h_x^2}+\frac{A_{i,j}^{p,q}-2A_{i,j}^{p,q}+A_{i,j+1}^{p,q}}{h_y^2}\right).
\end{align}
Equations~\eqref{eq:first:step} and~\eqref{eq:first:reg:partials} always involve an interpolation scheme represented by the symbol $I[A] (r)$ in order to evaluate the MPI Core Operator in the points $r_k$, that in general differ from the pixel centers $(x_i ,y_j )$. In the experiments performed in this paper we use bicubic interpolation, that is, if the point $(r_x ,r_y)$ is contained in the $(i,j)$-th pixel with center $(x_i ,y_j )$, then the interpolation $I[A]$ in $(r_x ,r_y)$ is performed as follows:
\begin{align}\label{eqn:interpolation}
	I[A^{p,q}] (r_x , r_y ) =&
	\begin{pmatrix}
		L_{-1}(s_x) \\ L_{0}(s_x) \\ L_{+1}(s_x) \\ L_{+2}(s_x)
	\end{pmatrix}^T
	\begin{pmatrix}
		A_{i-1,j-1}^{p,q} & A_{i-1,j}^{p,q} & A_{i-1,j+1}^{p,q} & A_{i-1,j+2}^{p,q} \\
		A_{i  ,j-1}^{p,q} & A_{i  ,j}^{p,q} & A_{i  ,j+1}^{p,q} & A_{i  ,j+2}^{p,q} \\
		A_{i+1,j-1}^{p,q} & A_{i+1,j}^{p,q} & A_{i+1,j+1}^{p,q} & A_{i+1,j+2}^{p,q} \\
		A_{i+2,j-1}^{p,q} & A_{i+2,j}^{p,q} & A_{i+2,j+1}^{p,q} & A_{i+2,j+2}^{p,q} \\
	\end{pmatrix}
	\begin{pmatrix}
		L_{-1}(s_y) \\ L_{0}(s_y) \\ L_{+1}(s_y) \\ L_{+2}(s_y)
	\end{pmatrix},
\end{align}
where we define the Lagrange polynomials
\begin{align}\label{eq:lagrange}
	L_{-1}(s)&=-\frac{1}{6}s(s-1)(s-2), & L_{0}(s) &=\frac{1}{2}(s+1)(s-1)(s-2),\\
	L_{+1}(s)&=-\frac{1}{2}s(s+1)(s-2), & L_{+2}(s) &=\frac{1}{6}s(s+1)(s-1), \notag
\end{align}
with values $s_x = \frac{r_x - x_i}{h_x}$ and $s_y = \frac{r_y - y_j}{h_y}$. Finally, the gradient system $\nabla J[A]=0$ can be rewritten as a non-homogeneous linear system of the form $GA=b$ where
\begin{align}
	(GA)_{i,j}^{p,q} &= \frac{2}{L}\sum_{k=1}^L \left\langle  I[A](r_k )v_k , I\left [ \delta_{i,j}^{p,q}\right ] (r_k )v_k\right\rangle +\frac{\partial R[A]}{\partial A_{i,j}^{p,q}}, \label{eq:1:stage:system:G}\\
	b_{i,j}^{p,q} &=\frac{2}{L}\sum_{k=1}^L \left\langle s_k , I\left [\delta_{i,j}^{p,q}\right ] (r_k )v_k\right\rangle ,\label{eq:1:stage:system:b}
\end{align}
and $GA=b$ is solved with respect to $A$ with the conjugate gradient method~\cite{golub2013matrix}, for $G$ is symmetric and positive definite, being the Hessian of a convex quadratic function. A pseudo-code of the first stage algorithm can be seen in Alg.~\ref{alg:first}. The algorithm in Alg.~\ref{alg:first} is mainly divided into three sections: (i) the computation of the Lagrange polynomials that involves first associating each point to the center of the pixel where it falls in and to compute the quantities in Eq.~\eqref{eq:lagrange}(both operations require $\mathcal{O}(L)$ computations); (ii) the setup of the right hand side in Eq.~\eqref{eq:1:stage:system:b} and the evaluation of $GA$ in each iteration of the CG method in Eq.~\eqref{eq:1:stage:system:G} is of complexity $\mathcal{O} (L\cdot N_x \cdot N_y )$, because the evaluations in Eq.~\eqref{eqn:interpolation} for the bicubic interpolation are in constant time (of course the complexity of the interpolation plays a role in the complexity of the algorithm for global interpolation schemes).

\begin{algorithm}
	\caption{Pseudocode for the first stage of the algorithm}\label{alg:first}
	\textbf{Input}: scan data $\lbrace (s_k , r_k ,v_k )\colon k=1,\dots L\rbrace$ in the frame of reference of $\Omega = [a,b]\times [c,d]$;\\
	\textbf{Output}: reconstructed trace $u\in\mathbb{R}^{N_x\times N_y}$ of the MPI Core Operator;
	
	\begin{algorithmic}
		\For{$k=1$ \textbf{to} $L$}
		\State $L_x (k) \gets [L_{-1}(s_x^{(k)}),L_{0}(s_x^{(k)}),L_{+1}(s_x^{(k)}),L_{+2}(s_x^{(k)})]$; \Comment{Compute and store the Lagrange polynomials in Eq.~\eqref{eq:lagrange}}
		\State $L_y (k) \gets [L_{-1}(s_y^{(k)}),L_{0}(s_y^{(k)}),L_{+1}(s_y^{(k)}),L_{+2}(s_y^{(k)})]$;    
		\EndFor
		\For{\textbf{each} $p,q,i,j$}
		\State $b_{i,j}^{p,q} \gets \frac{2}{L}\sum_{k=1}^L \left\langle s_k , I\left [\delta_{i,j}^{p,q}\right ] (r_k )v_k\right\rangle$; \Comment{Setup the right hand side  fo $GA=b$}
		\EndFor
		\Procedure{ApplyMatrixG}{$A$}
		\For{\textbf{each} $p,q,i,j$}
		\State $(GA)_{i,j}^{p,q} \gets \frac{2}{L}\sum_{k=1}^L \left\langle  I[A](r_k )v_k , I\left [ \delta_{i,j}^{p,q}\right ] (r_k )v_k\right\rangle +\frac{\partial R[A]}{\partial A_{i,j}^{p,q}}$; \Comment{$L_x$,$L_y$ are used to compute $I[A]$ with Eq.~\eqref{eqn:interpolation}}
		\EndFor
		\State \textbf{return} $GA$;
		\EndProcedure
		\State $\tilde{A} \gets$ ConjGrad$($\Call{ApplyMatrixG}{$\bullet$}$,b)$;\Comment{Apply the CG method};
		\State u $\gets$ Trace($\tilde{A}$);
	\end{algorithmic}
\end{algorithm}

\textbf{Second Stage} The problem to solve in the second stage is a regularized deconvolution, as the relationship between the distribution of particles $\rho $ and the output of the first stage is expressed in the model by Eq.~\eqref{eq:second:stage}. Again, let us consider an $N_x\times N_y$ grid  discretization of $\Omega = [a,b]\times [c,d]$ and reconstruct an approximation of the original distribution of the form $\rho\in\mathbb{R}^{N_x\times N_y}$, indexed as $\rho_{i,j}\approx \rho (x_i ,y_j )$. Because of the modeling of the scanner introduced in Section~\ref{sec:mathematical:modeling}, the support of $\rho$ is completely contained in $\Omega$, and the convolution in the energy functional $E[\rho ]$ in Eq.~\eqref{eq:second:stage} is equivalent to an integral over $\Omega$ instead of $\mathbb{R}^2$, that is:
\begin{equation}\label{eq:convo:omega}
	(\kappa_h * \rho ) (x) = \int_{\Omega} \rho (x)\kappa_h (x-y) \, dy ,
\end{equation}
which will be approximated using the midpoint rule. The discretization (the sampling) of the convolution kernel $\kappa_h$ is performed on the same discretization grid as the target solution $\rho_{i,j}$, leading to a matrix-vector product $K_h \rho$ which approximates the convolution in Eq.~\eqref{eq:convo:omega} and can be easily handled using discrete convolution based on FFT. We point out that the matrix $K_h$ is symmetric because $\kappa_h (x)=\kappa_h (-x)$.

To discretize the TV Smooth regularizer in Eq.~\eqref{eq:tv:smooth:reg} we proceed as follows: first, consider the forward and backward first order difference operators of the discretized target distribution $\rho_{i,j}$:
\begin{equation}\label{eq:firs:ord:diff}
	D_x^+ \rho_{i,j} = \frac{\rho_{i+1 ,j}-\rho_{i,j}}{h_x}, \quad D_x^- \rho_{i,j} = \frac{\rho_{i ,j}-\rho_{i-1,j}}{h_x}
\end{equation}
where $h_x = \frac{b-a}{N_x}$ is the discretization step size (the operators $D_y^+$ and $D_y^-$ are defined analogously with step size $h_y=\frac{d-c}{N_y}$). The operators in Eq.~\eqref{eq:firs:ord:diff} are such that
\begin{equation}
	\frac{(D_x^+ \rho_{i,j})^2 + (D_x^- \rho_{i,j})^2}{2} = \partial_x \rho (x_i , y_j )^2 + O (h_x^2 ).
\end{equation}
This means that we have the averages $W_{i,j}$ of the form
\begin{equation}\label{eq:def:wij}
	W_{i,j}\coloneq \frac{(D_x^+ \rho_{i,j})^2 + (D_x^- \rho_{i,j})^2}{2} + \frac{(D_y^+ \rho_{i,j})^2 + (D_y^- \rho_{i,j})^2}{2} ,
\end{equation}
which give a second order accurate approximation of $\lvert \nabla\rho\rvert^2$ at the cell centers. Because we assume that the support of $\rho$ is contained in $\Omega$, we impose Dirichlet zero boundary conditions which are implemented with zero-padding of $\rho$. Then, the TV Smooth regularizer of Eq.~\eqref{eq:tv:smooth:reg} can be discretized as follows:
\begin{equation}\label{eq:tv:smooth:discrete}
	R_\delta [\rho ] = h_x h_y \sum_{i=0}^{N_x -1}\sum_{j=0}^{N_y -1}\sqrt{W_{i,j}+\delta}.
\end{equation}
\begin{rmk}
	Any choice of a discretization of TV raises questions on its isotropicity. As discussed in~\cite{hosseini2023secondorder}, even the classic isotropic TV model (ROF~\cite{rudin1992nonlinear}) is not literally isotropic, as it lacks rotational invariance with respect to $90^\circ$ rotations. The discretization in Eq.~\eqref{eq:tv:smooth:discrete} however, considers in each pixel the differences with all neighboring (adjacent but not diagonally adjacent) pixels, guaranteeing isotropicity at least with respect to  $90^\circ$ rotations.
\end{rmk}

In particular, we have the discrete version of the minimization problem in Eq.~\eqref{eq:second:pos:sparsity}:
\begin{align}\label{eq:second:pos:sparsity:discrete}
	\rho = \arg\min_{\hat{\rho}} \biggl\lbrace \lVert K_h \hat{\rho}-u\rVert_2^2 + \mu h_x h_y \sum_{i=0}^{N_x -1}\sum_{j=0}^{N_y -1}\sqrt{W_{i,j}+\delta } + \beta\lVert\hat{\rho}\rVert_1 +  \iota_+ (\hat{\rho}) \biggr\rbrace .
\end{align}
With the discretization in Eq.~\eqref{eq:tv:smooth:discrete}, the partial derivatives of the regularizer are:
\begin{equation}\label{eq:partial:of:reg}
	\frac{\partial R_\delta}{\partial \rho_{i,j}} [\rho ] = -\frac{A_x^+ g_{i,j} D_x^+ \rho_{i,j}-A_x^- g_{i,j}D_x^- \rho_{i,j}}{h_x} - \frac{A_y^+ g_{i,j} D_y^+ \rho_{i,j}-A_y^- g_{i,j}D_y^- \rho_{i,j}}{h_y}
\end{equation}
where $g_{i,j}=\frac{1}{\sqrt{\delta + W_{i,j}}}$ and 
\begin{equation}\label{eq:operators:a:differences}
	A_x^+ g_{i,j} \coloneq \frac{g_{i+1,j}+g_{i,j}}{2},\quad A_x^- g_{i,j} \coloneq \frac{g_{i,j}+g_{i-1,j}}{2}
\end{equation}
are the forward and backward average operators applied to $g$ ($A_y^+$ and $A_y^-$ are defined analogously).

The discretizations performed so far can be more compactly represented in matrix form. Indeed, if $D$ denotes the first order discretization matrix of the gradient and $M[\rho ] = \text{diag}(g_{1,1},\dots ,g_{N_x N_y ,N_x N_y})$ the diagonal matrix representing the multiplication with $g=\frac{1}{\sqrt{\delta + \lvert\nabla\rho\rvert^2}}$, and $A$ the matrix representing the application of the averages $A_x^+$,$A_x^-$,$A_y^+$ and $A_y^-$ in Eq.~\eqref{eq:operators:a:differences}, then the gradient of the regularizer is 
\begin{equation}\label{eq:grad:tv:disc}
	\nabla R_\delta [\rho ] = -D^T A M[\rho ]D\rho ,
\end{equation}
and the matrix form of the Euler-Lagrange equation $\nabla E [\rho ]=0$ with the energy function $E$ in Eq.~\eqref{eq:second:stage} is:
\begin{equation}\label{eq:EL:tv:smooth}
	-\mu D^T A M[\rho ]D\rho + K_h^T (K_h \rho - u) = 0.
\end{equation}
For the non-negative fused Lasso scheme, the left-hand side in the Euler-Lagrange Eq.~\eqref{eq:EL:tv:smooth} is actually the gradient of the differentiable part $F[\rho ]$ in Eq.~\eqref{eq:second:pos:sparsity}, i.e.,
\begin{equation}\label{eq:grad:f}
	\nabla F[\rho ] = -\mu D^T A M[\rho ]D\rho + K_h^T (K_h \rho - u).
\end{equation}

Finally, one can explicitly compute the proximal mappings needed in Eq.~\eqref{eq:alg:lasso}. In particular, they are separable with respect to the grid centers, so it is enough to compute the action of the proximal mappings on each pixel:
\begin{align}
	\biggl ( \prox_{2\gamma G_1} (\rho )\biggr )_{i,j} & = \arg\min_{r\in\mathbb{R}}\bigl\lbrace 2\gamma G_1 (r) + \frac{1}{2}\lVert r - \rho_{i,j}\rVert_2^2\bigr\rbrace = \sign (\rho_{i,j})\cdot\max\bigl (\lvert \rho_{i,j}\rvert - 2\gamma\beta ,0\bigr ) \\
	\biggl (\prox_{2\gamma \iota_+} (\rho ) \biggr )_{i,j} & = \arg\min_{r\in\mathbb{R}}\bigl\lbrace 2\gamma\iota_+ (r)+\frac{1}{2}\lVert r - \rho_{i,j}\rVert_2^2\bigr\rbrace = \max (0, \rho_{i,j}) .
\end{align}
\begin{algorithm}
	\caption{Pseudocode for the second stage of the algorithm}\label{alg:lasso}
	\textbf{Input}: trace $u\in\mathbb{R}^{N_x\times N_y}$ of the MPI Core Operator from the first stage and parameters $\mu>0$, $\beta >0$,$\gamma >0$;\\
	\textbf{Output}: reconstructed distribution $\rho\in\mathbb{R}^{N_x\times N_y}$;
	
	\begin{algorithmic}
		\State $z_1 \gets u$;
		\State $z_2 \gets u$;
		\State $\rho \gets u$;
		\Repeat
		\State $\xi \gets 2\rho + \gamma \bigl( - \mu D^T A M[\rho ]D\rho + K_h^T (K_h\rho -u)\bigl ) $;
		\For{$i=0$ \textbf{to} $N_x-1$}
		\For{$j=0$ \textbf{to} $N_y-1$}
		\State $(z_1)_{i,j} \gets (z_1 )_{i,j} + \sign(\xi_{i,j})\cdot \max \bigl ( \lvert\xi_{i,j}\rvert - 2\gamma\beta , 0\bigr ) - \rho_{i,j}$;
		\State $(z_2 )_{i,j} \gets (z_2 )_{i,j } + \max (0, \xi_{i,j}-(z_2 )_{i,j})-\rho_{i,j}$;
		\EndFor
		\EndFor
		\State $\rho \gets \frac{1}{2}(z_1 + z_2 )$;
		\Until{the stopping criterion is met;}
	\end{algorithmic}
\end{algorithm}

A pseudocode for the deconvolution algorithm can be seen in Alg.~\ref{alg:lasso}, where we initialize  $\rho$,$z_1$ and $z_2$ with the trace $u$ of the MPI Core Operator reconstructed in the first stage of the algorithm. The parameter $\lambda_k$ in Eq.~\eqref{eq:alg:lasso} were chosen for simplicity to be $\lambda_k =1$ for every $k$. In Alg.~\ref{alg:lasso} the computation of the gradient (Eq.~\eqref{eq:grad:tv:disc}) can be computed using Eq.~\eqref{eq:partial:of:reg} and \eqref{eq:operators:a:differences}, yielding complexity of order $\mathcal{O}(N_x\cdot N_y)$, while the term $K_h^T K_h\rho$ can be computed using the FFT algorithm in $\mathcal{O}(N_x\cdot N_y \cdot \log N_x \cdot \log N_y )$. From this it follows that in Alg.~\ref{alg:lasso} $\mathcal{O}(N_x\cdot N_y \cdot \log N_x \cdot \log N_y )$ operations are performed in each iteration.

\paragraph{\textbf{Summary of the Approach}} For better clarity we summarize the two stage algorithm. The input data consist of the scan data $s_k$ and of the positions $r_k$ and velocities $v_k$ of the scanning trajectory at time $t_k$, for given indices $k$. These data are linked by Eq.~\eqref{eq:core:discrete} via the MPI Core Operator $A$, containing the information about the target distribution $\rho$. The MPI Core Operator is reconstructed at each point of the chosen reconstruction grid as follows: we minimize the functional $J[A]$ in Eq.~\eqref{eq:first:step}, which imposes the relation in Eq.~\eqref{eq:core:discrete} in the regularized least square sense, where the chosen regularizer imposes smoothness of the solution, as motivated by theory. The minimization of the functional $J[A]$ is performed by using the CG method applied to the linear system of equations obtained from the Euler-Lagrange equations relative to $J[A]$ (see Alg.~\ref{alg:first}). Because the trace of the MPI Core Operator is the convolution of $\rho$ with a specific kernel, we consider the scalar-valued field $u$ given by the traces of the matrix-valued field $A$ at each pixel, obtained in the first stage. We then perform a regularized deconvolution with the known kernel, choosing as regularizer the TV-Smooth approximation of the TV norm and adding positivity constraints and a sparsity promoting prior. This deconvolution is performed by minimizing the energy $E$ in Eq.~\eqref{eq:second:pos:sparsity}. Because the energy $E$ involves both smooth terms and non-smooth but convex terms, the algorithm that performs the deconvolution uses an iterative forward-backward splitting scheme which uses proximal mappings to deal with the convex non-smooth priors. This specific deconvolution algorithm is shown in Alg.~\ref{alg:lasso}, whose convergence is proved in Section~\ref{sec:convergence:alg}.

\subsection{Convergence of the Algorithm}\label{sec:convergence:alg}

We now provide a proof of convergence of Alg.~\ref{alg:lasso} to a minimizer of the energy functional in Eq.~\eqref{eq:second:pos:sparsity:discrete} with the chosen parameters. The proof will be mainly based on applying Theorem 2.1 in Raguet \emph{et al.}\cite{doi:10.1137/120872802} after proving that $\nabla F$ in Eq.~\eqref{eq:grad:f} is Lipschitz.

\begin{theorem}\label{thm:convergence}
	Let  $\eta = \frac{1}{\Lip (\nabla F)}$, $\gamma\in (0,2\eta)$ and $\lambda_k \in \biggl (0, \min \bigl (\frac{3}{2},\frac{1}{2}+\frac{\eta}{\gamma}\bigr ) \biggr )$ for $k\in\mathbb{N}$ and $\nabla F$ as in Eq.~\eqref{eq:grad:f}. Then the sequence $\rho^{(k)}$ defined in Alg.~\ref{alg:lasso} converges to a minimizer of the functional in Eq.~\eqref{eq:second:pos:sparsity:discrete}.
\end{theorem}
\vspace{\baselineskip}
\begin{proof}
	If $\nabla F$ is Lipschitz, then assumptions (A1),(A2) in Theorem 2.1 in~\cite{doi:10.1137/120872802} are satisfied. Consequently, weak convergence of the sequence $\rho^{(k)}$ to a global minimizer of Eq.~\eqref{eq:second:pos:sparsity} is guaranteed. Of course, because in the discretized functional in Eq.~\eqref{eq:second:pos:sparsity:discrete} everything is in a Euclidean space, weak converge is equivalent to strong convergence.
	
	We now prove that the gradient $\nabla F$ is indeed Lipschitz and we give an estimation of the Lipschitz constant.
	To this aim, let $\rho_1$,$\rho_2\in\mathbb{R}^{N_x\cdot N_y}$ and we separate the study of the Lipschitz constant of $\nabla F$ by splitting the estimate using the triangle inequality and the fact that the operator $K_h$ is linear:
	\begin{align}\label{eq:lipschitz:constant:separation}
		\lVert \nabla F (\rho_1 )-\nabla F(\rho_2 )\rVert_2 & = \bigl\lVert -\mu D^T A M[\rho_1 ]D\rho_1 + K_h^* K_h (\rho_1 -u) +\mu D^T A M[\rho_2 ]D\rho_2 - K_h^* K_h (\rho_2 -u) \bigr\rVert_2 \notag\\
		& \leq \mu \bigl\lVert D^T A M[\rho_1 ]D\rho_1 - D^T A M[\rho_2 ]D\rho_2 \bigr\rVert_2 + \bigl\lVert K_h^* K_h \bigr\rVert_{2} \lVert \rho_1 - \rho_2 \rVert_2
	\end{align}
	Now, we need to estimate the Lipschitz constant of the operator $\nabla R_{\delta}[\rho ]=-D^T A M[\rho ]D\rho $. The idea is to use the Mean Value Theorem~\cite{edwardsadvancedcalculus} and obtain a Lipschitz constant by taking an upper bound of the Frobenius norm of the Hessian matrix of $R_{\delta}$. This means that we need to estimate a bound for the entries of the Hessian matrix of the regularizer $R_{\delta}$. We compute the action of the Hessian using the calculus of variations: consider a variation $\psi\in\mathbb{R}^{N_x\cdot N_y}$. We evaluate the gradient $\nabla R_{\delta}$ in $\rho + t\psi$ and consider the directional derivative along $\psi$:
	\begin{align}
		\lim_{t\to 0}\frac{\nabla R_{\delta}[\rho + t\psi ]-\nabla R_{\delta}[\rho ]}{t} & = \lim_{t\to 0}\frac{-D^T A M[\rho +t\psi]D(\rho + t\psi )+D^T A M[\rho]D\rho}{t} \notag \\
		& = -D^T A M[\rho]D\psi  + \lim_{t\to 0} - D^T A\biggl [ \frac{M[\rho + t\psi]-M[\rho ]}{t}\biggr ]D\rho . \label{eq:frechet:grad:R:2}
	\end{align}
	Now, the matrix $M[\rho + t\psi ]$ is a diagonal $(N_x\cdot N_y )\times (N_x \cdot N_y )$ matrix\footnote{We use here a double indexing notation, writing $(ij,mn)$ in place of $\bigl ( (i,j)\, , (m,n)\bigr )$ for the $\bigl ( (i,j)\, , (m,n)\bigr )$-th entry of an $N_x N_y \times N_x N_y$ matrix.} whose diagonal entries are of the form:
	\begin{equation}\label{eq:M:diag:plus:var:entry}
		\begin{split}
			\biggr ( M[\rho + t\psi ]\biggr )_{ij,ij} = \biggl ( \delta + \frac{1}{2} \bigl  ( D_x^+ \rho_{i,j} + t D_x^+ \psi_{i,j} \bigr )^2 & +\frac{1}{2} \bigl ( D_x^- \rho_{i,j} + t D_x^- \psi_{i,j}\bigr )^2 + \\
			 & + \frac{1}{2} \bigl ( D_y^+ \rho_{i,j} + t D_y^+ \psi_{i,j}\bigr )^2 + \frac{1}{2} \bigl ( D_y^- \rho_{i,j} + t D_y^- \psi_{i,j}\bigr )^2\biggr )^{-\frac{1}{2}} .
		\end{split}
	\end{equation}
	The idea is now to expand the squares in Eq.~\eqref{eq:M:diag:plus:var:entry}. Indeed, if we recall the definition of $W_{i,j}$ in Eq.~\eqref{eq:def:wij} and consider a short-hand notation for the mixed terms:
	\begin{equation}
		B_1^{i,j} = D_x^+ \rho_{i,j}D_x^+ \psi_{i,j},\quad
		B_2^{i,j} = D_x^- \rho_{i,j}D_x^- \psi_{i,j},\quad
		B_3^{i,j} = D_y^+ \rho_{i,j}D_y^+ \psi_{i,j},\quad
		B_4^{i,j} = D_y^- \rho_{i,j}D_y^- \psi_{i,j} ,
	\end{equation}
	then we can factor out the term $(\delta + W_{i,j})^{-\frac{1}{2}}$ in Eq.~\eqref{eq:M:diag:plus:var:entry}:
	\begin{equation}\label{eq:M:diag:plus:var:entry:factored}
		\biggr ( M[\rho + t\psi ]\biggr )_{ij,ij} =\frac{1}{\sqrt{\delta +W_{i,j}}} \biggl ( 1 + \frac{B_1^{i,j} + B_2^{i,j} + B_3^{i,j} + B_4^{i,j}}{\delta + W_{i,j}} t + O(t^2 )\biggr )^{-\frac{1}{2}}
	\end{equation}
	We aim to use the Taylor expansion $\frac{1}{\sqrt{1+x}}=1 - \frac{1}{2}x + O (x^2 )$ if $\lvert x\rvert <1$. It is indeed applicable, as for $t$ small enough, we have eventually:
	\begin{equation}
		\biggl \lvert \frac{t\sum_{\alpha =1}^{4}B_{\alpha}^{i,j}}{\delta + W_{i,j}}\biggr\rvert <1 .
	\end{equation}
	In particular, we obtain the following estimate for the $(i,j)$-th diagonal entry of $M[\rho + t\psi ]$:
	\begin{equation}
		\biggr ( M[\rho + t\psi ]\biggr )_{ij,ij} = \frac{1}{\sqrt{\delta + W_{i,j}}}\biggl (1- \frac{t}{2}\biggl ( \frac{B_1^{i,j} + B_2^{i,j} + B_3^{i,j} + B_4^{i,j}}{\delta + W_{i,j}} \biggr ) + O (t^2 ) \biggr ) .
	\end{equation}
	A Taylor expansion for the diagonal entries of the difference in Eq.~\eqref{eq:frechet:grad:R:2} immediately follows:
	\begin{equation}
		\biggr ( M[\rho + t\psi ] - M[\rho ]\biggr )_{ij,ij}= -\frac{t}{2}\frac{( B_1^{i,j} + B_2^{i,j} + B_3^{i,j} + B_4^{i,j})}{\bigl ( \delta + W_{i,j}\bigr )^{\frac{3}{2}}} + O(t^2 ) ,
	\end{equation}
	and consequently, taking the limit for $t\to 0$, the action of the Hessian matrix of $R$, which we denote by $H$, on $\psi$ can be written in matrix form as:
	\begin{equation}\label{eq:hessian:action}
		H(\rho )\psi = - D^T A B[\rho,\psi ]D\rho -  D^T A M[\rho ]D\psi ,
	\end{equation}
	where $B[\rho ,\psi ]$ is the diagonal matrix whose $(i,j)$-th diagonal entry is 
	\begin{equation}\label{eq:b:rho:psi:entry}
		\biggl ( B[\rho ,\psi ] \biggr )_{ij,ij} = -\frac{1}{2}\frac{( B_1^{i,j} + B_2^{i,j} + B_3^{i,j} + B_4^{i,j})}{\bigl ( \delta + W_{i,j}\bigr )^{\frac{3}{2}}} .
	\end{equation}
	Now, we aim at a bound on the Frobenius norm of the Hessian independent of $\rho$, as it would give us a Lipschitz constant for the non-linear operator $\nabla R_{\delta}$. From now on, $e_{i,j}$ will denote the $(i,j)$-th vector of the canonical basis of $\mathbb{R}^{N_x\cdot N_y}$, and $\langle\cdot  \, ,\cdot\rangle$ the standard scalar product in $\mathbb{R}^{N_x\cdot N_y}$. With this notation, the $(mn ,ij)$-th entry of the Hessian matrix $H$ in a point $\rho$ is:
	\begin{equation}\label{eq:hessian:matrix:entry}
		\bigl \langle e_{m,n}\, , H(\rho )e_{i,j}\bigr \rangle = - \underbrace{\bigl\langle e_{m,n}\, , D^T A B[\rho , e_{i,j}]D\rho\bigr\rangle }_{(\rom{1} )} - \underbrace{\bigl\langle e_{m,n}\, , D^T A M[\rho ]De_{i,j}\bigr\rangle}_{(\rom{2} )} .
	\end{equation}
	Before separately evaluating $(\rom{1})$ and $(\rom{2} )$, we recall - see Eq.~\eqref{eq:partial:of:reg} - that the operator $D^T A B D $ acts for a general diagonal matrix $B = \text{diag}(b_{1,1},\dots ,b_{N_x ,N_y})$ on a point $\rho$ according to the following formula:
	\begin{equation}\label{eq:action:dbd}
		\bigl\langle e_{m,n}\, , D^T A BD\rho\bigr\rangle = \frac{A_x^+ b_{m,n} D_x^+ \rho_{m,n}-A_x^- b_{m,n}D_x^- \rho_{m,n}}{h_x} + \frac{A_y^+ b_{m,n} D_y^+ \rho_{m,n}-A_y^- b_{m,n}D_y^- \rho_{m,n}}{h_y} ,
	\end{equation}
	where the operators $A_x^+$,$A_x^-$,$A_y^+$ and $A_y^-$ are defined in Eq.~\eqref{eq:operators:a:differences}. We denote with $\mathcal{D}$ a generic difference operator among $D_x^+$,$D_x^-$,$D_y^+$ and $D_y^-$. Analogously, we will write $\mathcal{A}$ when we talk about any of the operators of the form $A_x^+$,$A_x^-$,$A_y^+$ or $A_y^-$.
	
	Using the definition of $\mathcal{A}$ and the triangular inequality, we have a first estimate of the absolute value of the general term in Eq.~\eqref{eq:action:dbd}:
	\begin{align}\label{eq:triangular:big:expansion}
		\biggl\lvert\bigl\langle e_{m,n}\, , D^T A BD\rho\bigr\rangle\biggr\rvert \leq &
		\biggl\lvert \frac{b_{m+1,n}D_{x}^+\rho_{m,n}}{2h_x}\biggr\rvert + 
		\biggl\lvert \frac{b_{m,n}D_{x}^+\rho_{m,n}}{2h_x}\biggr\rvert + 
		\biggl\lvert \frac{b_{m,n}D_{x}^-\rho_{m,n}}{2h_x}\biggr\rvert +
		\biggl\lvert \frac{b_{m-1,n}D_{x}^-\rho_{m,n}}{2h_x}\biggr\rvert + \dots \notag\\
		& \dots +
		\biggl\lvert \frac{b_{m,n+1}D_{y}^+\rho_{m,n}}{2h_y}\biggr\rvert +
		\biggl\lvert \frac{b_{m,n}D_{y}^+\rho_{m,n}}{2h_y}\biggr\rvert +
		\biggl\lvert \frac{b_{m,n}D_{y}^+\rho_{m,n}}{2h_y}\biggr\rvert +
		\biggl\lvert \frac{b_{m,n-1}D_{y}^+\rho_{m,n}}{2h_y}\biggr\rvert .
	\end{align}
	We now study $(\rom{1} )$ alone. The fact that we are dealing with the canonical basis vector $e_{i,j}$ in place of $\psi$ simplifies the estimation. Indeed, $\mathcal{D} e_{i,j}$ is always in the set $\lbrace -1, 0 , 1\rbrace$ and in particular the sum of the mixed terms $B_{\alpha }^{i,j}$ is of the form: 
	\begin{equation}
		\sum_{\alpha = 1}^{4}B_{\alpha}^{m,n} = \frac{\omega_1}{h_x}D_x^+\rho_{m,n}+\frac{\omega_2}{h_x}D_x^-\rho_{m,n}+\frac{\omega_3}{h_y}D_y^+\rho_{m,n}+\frac{\omega_4}{h_y}D_y^-\rho_{m,n} ,
	\end{equation}
	with the coefficients $\omega_\alpha\in\lbrace -1,0,1\rbrace$ for $\alpha = \lbrace 1,2,3,4\rbrace$.
	With this observation, we can use a simplified vector notation to carry on the computations. Indeed, defining the following vectors
	\begin{equation}\label{eq:vector:notation:differences}
		\overrightarrow{r}_{m,n}=\bigl (D_x^+\rho_{m,n},D_x^-\rho_{m,n},D_y^+\rho_{m,n},D_y^-\rho_{m,n}\bigr ) ,\quad\text{and}\quad \overrightarrow{\omega} = \biggl ( \frac{\omega_1}{h_x},\frac{\omega_2}{h_x},\frac{\omega_3}{h_y},\frac{\omega_4}{h_y}\biggr )
	\end{equation}
	the diagonal $(m,n)$-th entry of the diagonal matrix $B[\rho ,e_{i,j }]$ can be written as
	\begin{equation}
		b_{m,n} = -\frac{1}{2}\frac{\langle \overrightarrow{r}\, ,\overrightarrow{\omega}\rangle}{\bigl ( \delta + \frac{1}{2}\langle \overrightarrow{r}\, , \overrightarrow{r} \rangle \bigr )^{\frac{3}{2}}} .
	\end{equation}
	Moreover, the norm of $\overrightarrow{\omega}$ is bounded from above:
	\begin{equation}\label{eq:omega:norm:bound}
		\lVert\overrightarrow{\omega}\rVert_2 = \sqrt{\frac{\omega_1^2}{h_x^2} + \frac{\omega_2^2}{h_x^2}+\frac{\omega_3^2}{h_y^2}+\frac{\omega_4^2}{h_y^2}}\leq \sqrt{\frac{2}{h_x^2}+\frac{2}{h_y^2}}.
	\end{equation}
	We also observe that the following index shifts hold true:
	\begin{equation}\label{eq:shift:of:d}
		D_x^+ \rho_{m,n} = D_x^- \rho_{m+1 ,n}\quad,\quad D_x^- \rho_{m,n} = D_x^+ \rho_{m-1 ,n}\quad ,\quad  D_y^+ \rho_{m,n} = D_y^- \rho_{m ,n+1}\quad ,\quad D_y^- \rho_{m,n} = D_y^+ \rho_{m ,n-1}.
	\end{equation}
	In particular, it is clear that, up to a shift in Eq.~\eqref{eq:shift:of:d}, any term of the form $b_{k,l}\mathcal{D}$ in the expansion in Eq.~\eqref{eq:triangular:big:expansion} for $(k,l)$ in $\lbrace (m+1,n),(m,n),(m-1,n),(m,n+1),(m,n-1)\rbrace$ is of the form:
	\begin{equation}
		b_{k,l}\mathcal{D} = -\frac{1}{2}\frac{\langle \overrightarrow{r}\, ,\overrightarrow{\omega}\rangle}{\bigl ( \delta + \frac{1}{2}\langle \overrightarrow{r}\, , \overrightarrow{r} \rangle \bigr )^{\frac{3}{2}}} \langle\overrightarrow{r}\, , e_{\alpha}\rangle ,
	\end{equation}
	and $e_{\alpha}$ is a vector of the canonic basis of $\mathbb{R}^4$. Using the Cauchy-Schwarz inequality and the bound on the norm in Eq.~\eqref{eq:omega:norm:bound} we obtain the following estimate:
	\begin{equation}\label{eq:clubsuit:estimate}
		\lvert b_{k,l}\mathcal{D}\rvert \leq \frac{1}{2}\frac{\lVert \overrightarrow{r}\rVert_2^2 \lVert \overrightarrow{\omega}\rVert_2}{\bigl (\delta +\frac{1}{2}\lVert\overrightarrow{r}\rVert_2^2 \bigr )^{\frac{3}{2}}} \leq \sqrt{2}\biggl ( \frac{1}{h_x^2}+\frac{1}{h_y^2} \biggr )^{\frac{1}{2}} \max_{x\in\mathbb{R}^+}\lbrace\varphi (x)\rbrace , \quad\text{where }\varphi (x)\coloneq \frac{x}{\bigl ( \delta + x\bigr )^{\frac{3}{2}}}.
	\end{equation}
	Now, the function $\varphi$ defined in Eq.~\eqref{eq:clubsuit:estimate} achieves its maximum in $x=2\delta$, leading to the following estimate:
	\begin{equation}\label{eq:bkl:D:term:estimate:c1}
		\lvert b_{k,l}\mathcal{D}\rvert \leq \frac{2\sqrt{2}}{3^{\frac{3}{2}}\sqrt{\delta}}\sqrt { \frac{1}{h_x^2}+\frac{1}{h_y^2} } \eqcolon C_1 .
	\end{equation}
	Finally, applying the estimate of Eq.~\eqref{eq:bkl:D:term:estimate:c1} to each of the term in Eq.~\eqref{eq:triangular:big:expansion} we obtain
	\begin{equation}
		(\rom{1} )\leq 2\biggl (\frac{1}{h_x}+\frac{1}{h_y}\biggr ) C_1 .
	\end{equation}

	Next, we estimate $(\rom{2} )$, in this case the operator is $D^T A B D$ where $B$ is the diagonal matrix whose $(m,n)$-th diagonal entry is $b_{m,n}=\frac{1}{\sqrt{\delta + W_{m,n}}}$ and with the definition of $\overrightarrow{r}$ as in Eq.~\eqref{eq:vector:notation:differences} , any term $b_{k,l}\mathcal{D}$ for $(k,l)\in\lbrace (m+1,n),(m,n),(m-1,n),(m,n+1),(m,n-1)\rbrace$ is of the form:
	\begin{equation}
		b_{k,l}\mathcal{D} =\frac{1}{\sqrt{\delta +\frac{1}{2}\lVert \overrightarrow{r}\rVert_2^2}}\langle \overrightarrow{r}\, , e_\alpha \rangle ,
	\end{equation}
	for some vector $e_{\alpha}$ of the canonical basis of $\mathbb{R}^4$. Again, the Cauchy-Schwarz inequality allows us to find a bound on the absolute values of the term needed:
	\begin{equation}
		\biggl\lvert \frac{1}{\sqrt{\delta +\frac{1}{2}\lVert \overrightarrow{r}\rVert_2^2}}\langle \overrightarrow{r}\, , e_\alpha \rangle\biggr\rvert = \frac{\lvert\langle \overrightarrow{r}\, , e_{\alpha}\rangle\rvert}{\sqrt{\delta + \frac{1}{2}\lVert\overrightarrow{r}\rVert_2^2}}\leq \frac{\lVert\overrightarrow{r}\rVert}{\sqrt{\delta + \frac{1}{2}\lVert\overrightarrow{r}\rVert_2^2}} = \eta (\lVert\overrightarrow{r}\rVert_2)\quad\text{where}\quad \eta(x)=\frac{x}{\sqrt{\delta + \frac{x^2}{2}}}\text{ for }x\in\mathbb{R}.
	\end{equation}
	The newly defined function $\eta$ is bounded from above by $\sqrt{2}$, and applying this bound to each of the terms in Eq.~\eqref{eq:triangular:big:expansion} we have the following estimate of $(\rom{2} )$:
	\begin{equation}
		(\rom{2} )\leq 2\sqrt{2}\biggl (\frac{1}{h_x}+\frac{1}{h_y}\biggr ) .
	\end{equation}
	We therefore have a bound of the $(mm,ij)$-th entry of the Hessian matrix -- see Eq.~\eqref{eq:hessian:matrix:entry} -- and because the Hessian matrix is an $N_x\cdot N_y \times N_x\cdot N_y$ matrix, its Frobenius norm is bounded as well:
	\begin{align}
		\lVert H(\rho)\rVert_F = \sqrt{\sum_{m,i=1}^{N_x}\sum_{n,j=1}^{N_y}\bigl\lvert \bigl \langle e_{m,n}\, , H(\rho )e_{i,j}\bigr \rangle\bigr\rvert^2}  & \leq
		\sqrt{\sum_{k=1}^{N_x N_y}\sum_{l=1}^{N_x N_y}\biggl ( 2\biggl (\frac{1}{h_x}+\frac{1}{h_y}\biggr )C_1 + 2\sqrt{2} \biggl (\frac{1}{h_x}+\frac{1}{h_y}\biggr )\biggr )^2} \notag\\
		& =\sqrt{(N_x N_y)^2\biggl ( 2\biggl (\frac{1}{h_x}+\frac{1}{h_y}\biggr ) \bigl (\sqrt{2}+C_1\bigr )\biggr )^2} \notag\\
		& = 2N_x N_y \biggl (\frac{1}{h_x}+\frac{1}{h_y} \biggr ) (\sqrt{2}+C_1 ) \eqcolon L ,
	\end{align}
	and $L$ is a Lipschitz constant for $\nabla R_{\delta}$ and consequently, recalling Eq.~\eqref{eq:lipschitz:constant:separation}, 
	\begin{equation}\label{eq:lipschitz:constant}
		\tilde{L}\coloneq \mu L + \lVert K_h^* K_h \rVert_{2}
	\end{equation}
	is a Lipschitz constant for $\nabla F$.
\end{proof}

\section{Numerical Experiments}\label{sec:experiments}

In this section we describe the experiments performed. We use the two-stage algorithm described in Section~\ref{sec:rec:algorithm} and explain how the parameters in each of the stages have been optimized.

\subsection{Experimental Setup}

The algorithms described in this paper have been implemented in Python 3.9, using the packages Numpy, SciPy and PyTorch. The numerical experiments were performed on a Virtual Machine with 8 VCPUs, 16 GB of RAM and Ubuntu 22.04.1 LTS.

\paragraph{Simulation of the Data and Reconstruction Grid} We work with simulated data based on the simulation process in~\cite{marz2016model,math10183278}. We use different type of generalized Lissajous trajectories of the form in Eq.~\eqref{eq:gen:liss}, where the term $r(t)$ - the standard Lissajous trajectory in Eq.~\eqref{eq:lissajous} - mimics the setting in the Open MPI project~\cite{knopp2020openmpidata} with 1632 samples per scan period. The particular choices of offset trajectory $b(t)$ and of angle function $\alpha (t)$ for each experiment will be explained in the description of each experiment performed in what follows in this section. In general, given a ground truth distribution (also called phantom) $\rho_{\text{GT}}$, the data $s(t)$ are obtained using Eq.~\eqref{eq:signal:dimensionless} with a resolution parameter $h$ set to $h=0.01$ (see Remark~\ref{rmk:h}). To the data obtained we add $10\%$ Gaussian noise, i.e., the discrete data are simulated to be
\begin{equation}
	s_k = s(t_k )+\varepsilon N_k ,\quad \varepsilon = 0.1 \max_{k=1,\dots ,L}\lbrace \lvert s(t_k )\rvert\rbrace ,
\end{equation}
where $N_k\sim \mathcal{N}(0,1)$ are i.i.d. standard normally distributed random variables. This accounts for the inevitable noise of the measured data.

Each experiment below is performed on a domain $\Omega = [a,b]\times [c,d]$ with different choices of $a,b,c,d\in\mathbb{R}$. Particular choices are motivated when describing each experiment. In general, we employ an $N_x\times N_y$ grid discretization of $\Omega$ with $N_x =N_y \in\lbrace 100,200\rbrace $ for all experiments but Experiment 7 and 8, where we chose a $40\times 40$ grid. We point out that the choice of discretization grids finer than $40\times 40$ is to be already considered a higher resolution setup, as in current setups the grid sizes are about $40\times 40$ as in the Open MPI Project~\cite{knopp2020openmpidata}.

\paragraph{Quantitative Assessment}

To quantitatively assess the quality of the reconstructions obtained with the algorithm proposed, we use the Peak Signal-to-Noise Ratio (PSNR). We recall its definition: if $\rho_{\text{GT}}$ is the ground truth distribution and $\rho_{\text{rec}}$ is the reconstruction, the PSNR is defined as follows:
\begin{equation}\label{eq:psnr}
	\mathrm{PSNR} = 10\log_{10} \biggl (\frac{(\max \rho_{\text{GT}})^2}{\mathrm{MSE}}\biggr ) \quad\text{with }\mathrm{MSE}=\frac{1}{N_x N_y}\sum_{i=0}^{N_x -1}\sum_{j=0}^{N_y -1}\bigl ( \rho_{\text{GT},i,j}-\rho_{\text{rec},i,j} \bigr)^2 ,
\end{equation}
where MSE stand for the Mean Square Error, and $\rho_{\text{GT},i,j}$ (resp. $\rho_{\text{rec},i,j}$) denote the $(i,j)$-th pixel of the image. Because of its frequent employment for the measure of image quality, we have also computed the Structural Similarity Index Measure (SSIM)~\cite{Wang2004ImageQA} for all reconstructions.

\paragraph{Choice of Model and Regularization Parameters}

\def\imratio{0.35}
\begin{figure}[t]
	\centering
	\begin{subfigure}[t]{\imratio\linewidth}
		\includegraphics[width=\linewidth]{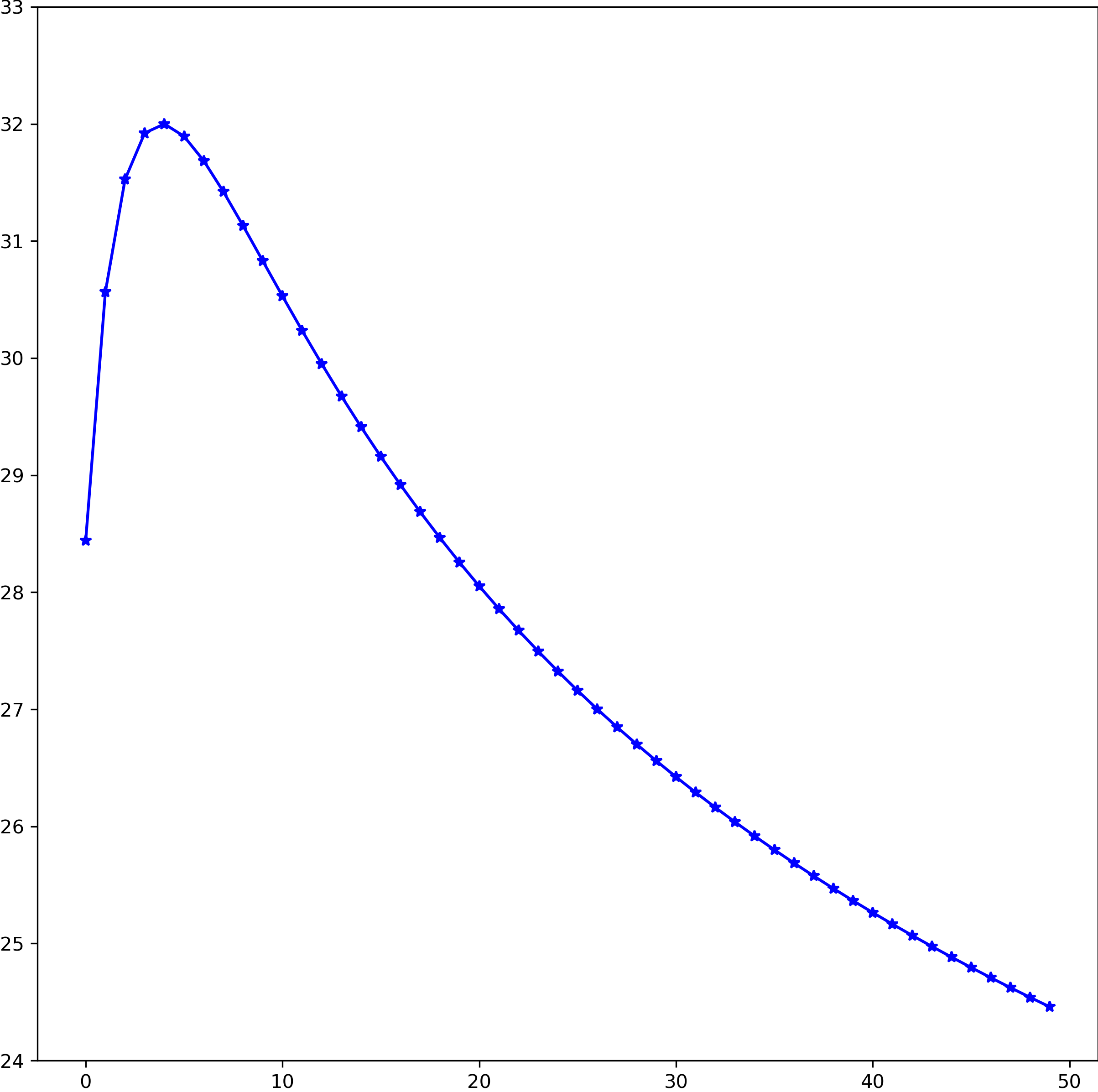}
		\caption{\centering\scriptsize PSNR vs. $\lambda$ for the first stage.}
		\label{subfig:plots:psnr:u}
	\end{subfigure}
	\hfil
	\begin{subfigure}[t]{\imratio\linewidth}
		\includegraphics[width=\linewidth]{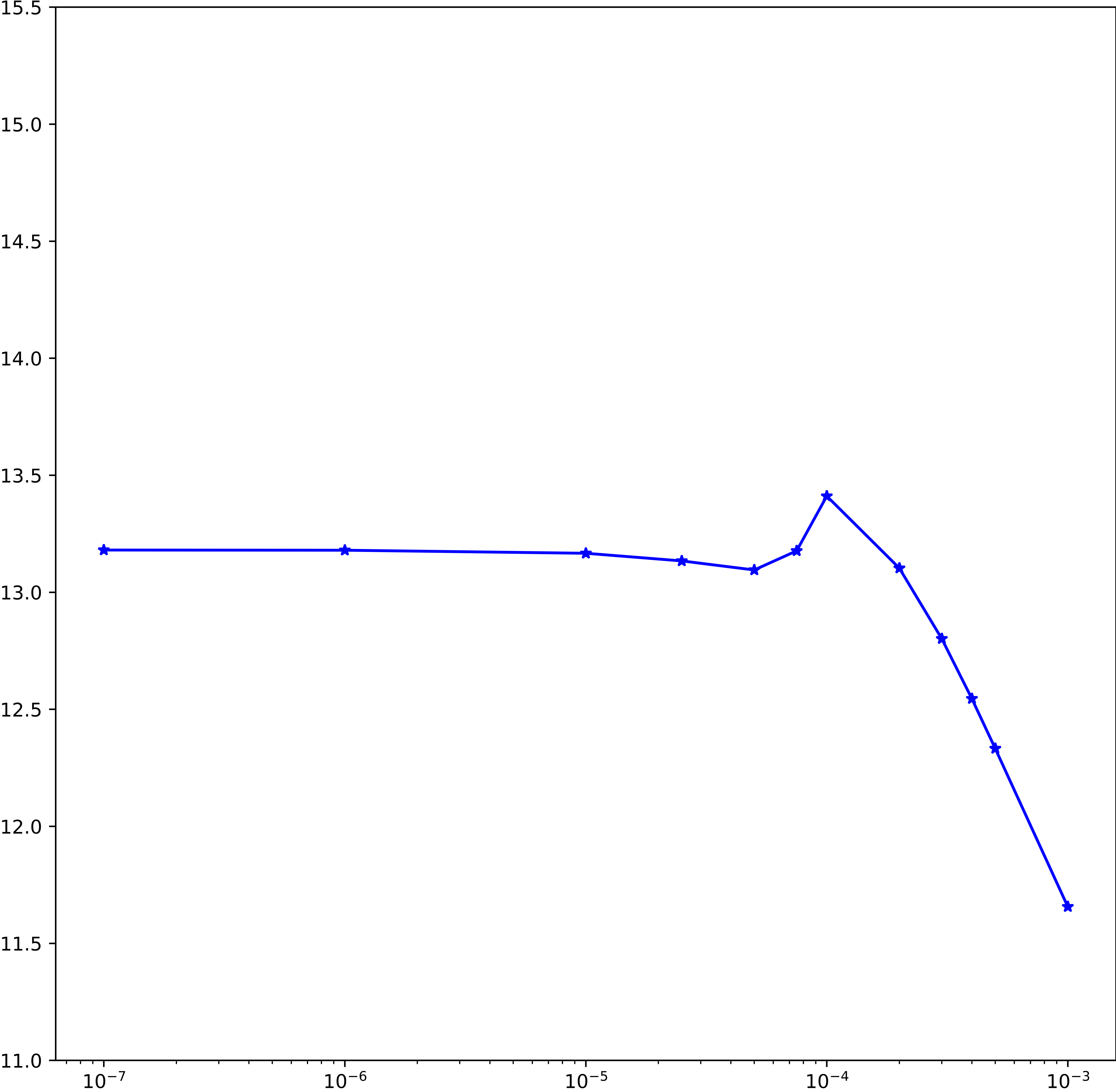}
		\caption{\centering\scriptsize PSNR vs. $\log (\mu )$ for the second stage.}
		\label{subfig:plots:psnr:rho}
	\end{subfigure}
	\caption{Plots of the behavior of the PSNR as the regularization parameters $\lambda$ for the first stage in Fig.~\ref{subfig:plots:psnr:u} and $\mu$ for the second stage in Fig.~\ref{subfig:plots:psnr:rho} vary. These curves refer to Experiment 1 (Fig.~\ref{subfig:exp:vess:multi:10:u},~\ref{subfig:exp:vess:multi:10:rho}) with 10x10 multi-patch scans (Fig.~\ref{subfig:scans:1010}). The parameter $\lambda$ in Fig.~\ref{subfig:plots:psnr:u} varies between 1 and 50 with step size 1, while in Fig.~\ref{subfig:plots:psnr:rho}, $\mu$ has been first optimized on $\lbrace 10^{-7},10^{-6},10^{-5},10^{-4},10^{-3}\rbrace$, and once the correct scale of $\mu$ has been spot, more points have been added around the presumed maximum value. The values of $\lambda$ and $\mu$ that maximize the PSNR have been chosen as the optimal reconstruction parameters of the first and second stages, respectively.}
	\label{plots:psnr}
\end{figure}

Both stages of the proposed reconstruction algorithm presented are variational minimization problems that involve regularization terms whose strengths are controlled by parameters which have to be chosen adequately. In particular, we need to find the parameter $\lambda$ in Eq.~\eqref{eq:first:step} for the first stage of the algorithm, as well as the parameters $\mu$ and $\beta$ in Eq.~\eqref{eq:second:pos:sparsity} for the second stage. In our experiments, we aim at presenting the potential of the proposed mathematical framework for multi-patching. To that end, we show results for in a sense optimal parameter choices (w.r.t. PSNR) being aware of the underlying ground truth. We point out that we do not employ a method of parameter selection in the sense of data analysis. Instead, we discuss the effect of the parameter on the final reconstruction later on; cf. the corresponding paragraph in Section~\ref{sec:disc:reg:params}. We point out that the automatic selection of the parameters is an important line of future research which is however not dealt with here.

Coming back, we choose parameters more precisely as follows: we perform reconstruction on a range of parameters and compute the PSNR defined in Eq.~\eqref{eq:psnr} w.r.t. the ground truth. In particular, the ground truth $\rho_{\text{GT}}$ is used to compute the PSNR value of the second stage, while in the first stage, the PSNR values of the traces are computed considering as ground truth the actual ground truth $\rho_{\text{GT}}$ convolved with the kernel $\kappa_h$ in Eq.~\eqref{eq:trace:convolution}. Concretely, for the first stage we let $\lambda$ range between 1 and 50 with step 1. For the second stage, we first performed reconstructions with $\mu\in\lbrace 10^{-7},10^{-6},10^{-5},10^{-4},10^{-3}\rbrace$ for every experiment (but experiment 5 in Fig.~\ref{fig:1000:recs} where it was necessary to decrease the magnitude and try values of magnitude $10^{-n}$ for $n\in\lbrace 8,\dots ,13\rbrace$.) After the correct magnitude has been pinpointed, denote it with $10^{-n}$, we fine-tune by considering $t \cdot 10^{-n-1}$ and $s\cdot 10^{-n}$ for $t\in\lbrace 2.5,5,7.5\rbrace$ and $s\in\lbrace 2,3,4,5\rbrace$. Because in the second stage the regularization of the deconvolution is mainly taken care of by the TV Smooth regularization term weighted by $\mu$ and the parameter $\beta$ is in charge of weighting the sparsity enforcing prior, we decided to optimize for $\mu$ while setting $\beta = 1$ for all experiments but the concentration phantom in Fig.~\ref{subfig:concentration:gt} and the Experiment in Fig.~\ref{fig:1000:recs}, where we set it to $\beta = 0.1$. The reason for this choice is explained in Sec.~\ref{sec:exp:MPI} and Sec.~\ref{sec:scan:move}. (Not fine-tuning for both $\mu$ and $\beta$ in the second stage of the algorithm allowed us to reduce the number of reconstructions to be performed.) In Fig.\ref{plots:psnr} we plotted an example of the general behavior of the PSNR as the parameter increases for both stages. In all cases, the PSNR graph of the first stage (e.g., Fig.~\ref{subfig:plots:psnr:u}) is a (non-symmetrically) bell-shaped function with a clear global maximum. The value for $\lambda$ (resp. $\mu$) corresponding to the maximum of the respective graph, and consequently, the reconstruction with said value, is chosen. In particular, the trace $u$ (highest PSNR) of the MPI Core Operator reconstructed in the first stage serves as input for the second stage, and the resulting reconstructed distribution (best w.r.t. PSNR) is display in the figures. The choice of the parameters with the baseline method in Experiment 7 and 8 are explained in the respective sections \ref{sec:exp:sm:22} and \ref{sec:exp:sm:lasso}.
\begin{rmk}\label{rem:inpainting}
	The parameter $\lambda$ in Eq.~\eqref{eq:first:step} is in charge of the strength of the inpainting effect of the Tikhonov regularizer. This means that a higher number of sample points which also corresponds in our experiments to a denser sampling (as increased density of the sampling is one of our goals), needs a lower inpainting strength, i.e., a lower value of $\lambda$. This intuitive explanation is supported by the decrease of the optimal parameter as the number of sample points increase (see Tab.\ref{tab:exp}).
\end{rmk}
\vspace{\baselineskip}

In all the experiments, the trace is computed solving the system $GA = b$ with $G$ in Eq.~\eqref{eq:1:stage:system:G} and $b$ in Eq.~\eqref{eq:1:stage:system:b} with CG with 1000 maximum iterations if a tolerance on the residuals of $10^{-12}$ is not reached before.

For Stage 2, we set $\delta = 10^{-16}$ and as a stopping criterion we considered the relative residuals, i.e., $\frac{\lVert \rho^{(k+1)}-\rho^{k}\rVert_2}{\lVert \rho^{k}\rVert_2}$. As a tolerance on the relative residuals we chose $\text{toll.} = 5\cdot 10^{-6}$ with $100,000$ maximum number of iterations. The choice for the value $5\cdot 10^{-6}$ is motivated by the fact that it is a trade off between less time-consuming but early stopped solutions when a tolerance of $10^{-5}$ is reached and the very time-consuming solutions obtained with a tolerance of $10^{-6}$. Finally, the convergence of Alg.~\ref{alg:lasso}, is guaranteed for a parameter choice that depends on the Lipschitz constant of the gradient of the differentiable part $F$ in Eq.~\eqref{eq:second:pos:sparsity}. The Lipschitz constant has been estimated from above in the proof of Theorem~\ref{thm:convergence}. This estimation however, leads in general to a big Lipschitz constant, which in turn produces a very small step size $\eta$ in the gradient descent part in Alg.~\ref{alg:lasso}. Indeed, in the proof of Theorem~\ref{thm:convergence}, we found out -- see Eq.~\eqref{eq:lipschitz:constant} -- that a Lipschitz constant for $\nabla F$ is
\begin{equation}
	\tilde{L} = 2\mu N_x N_y \biggl (\frac{1}{h_x}+\frac{1}{h_y} \biggr ) \biggl (\sqrt{2}+C_1 \biggr ) + \lVert K_h^* K_h \rVert_{2}, 
\end{equation}
where $C_1$ is the constant in Eq.~\eqref{eq:bkl:D:term:estimate:c1}. For our choice of $\delta = 10^{-16}$, $\eta\propto \frac{1}{\tilde{L}}$ is of order $10^{-11}$ for the experiments on a $200\times 200$ grid. If we increase $\delta$ up to $10^{-7}$ we get a constant $\eta$ of order $10^{-7}$, which is still very small. Moreover, further increasing $\delta$ compromises the quality of the approximation of the TV norm by the TV Smooth regularizer. Therefore, it is necessary to find by hand parameters $\gamma$ for which Alg.~\ref{alg:lasso} converges and that it is big enough. After performing many experiments we set the parameter $\gamma = 10^{-3}$ with which all reconstructions converged.

\def\imratio{0.32}
\begin{figure}[t]
	\centering
	\begin{subfigure}[t]{\imratio\linewidth}
		\includegraphics[width=\linewidth]{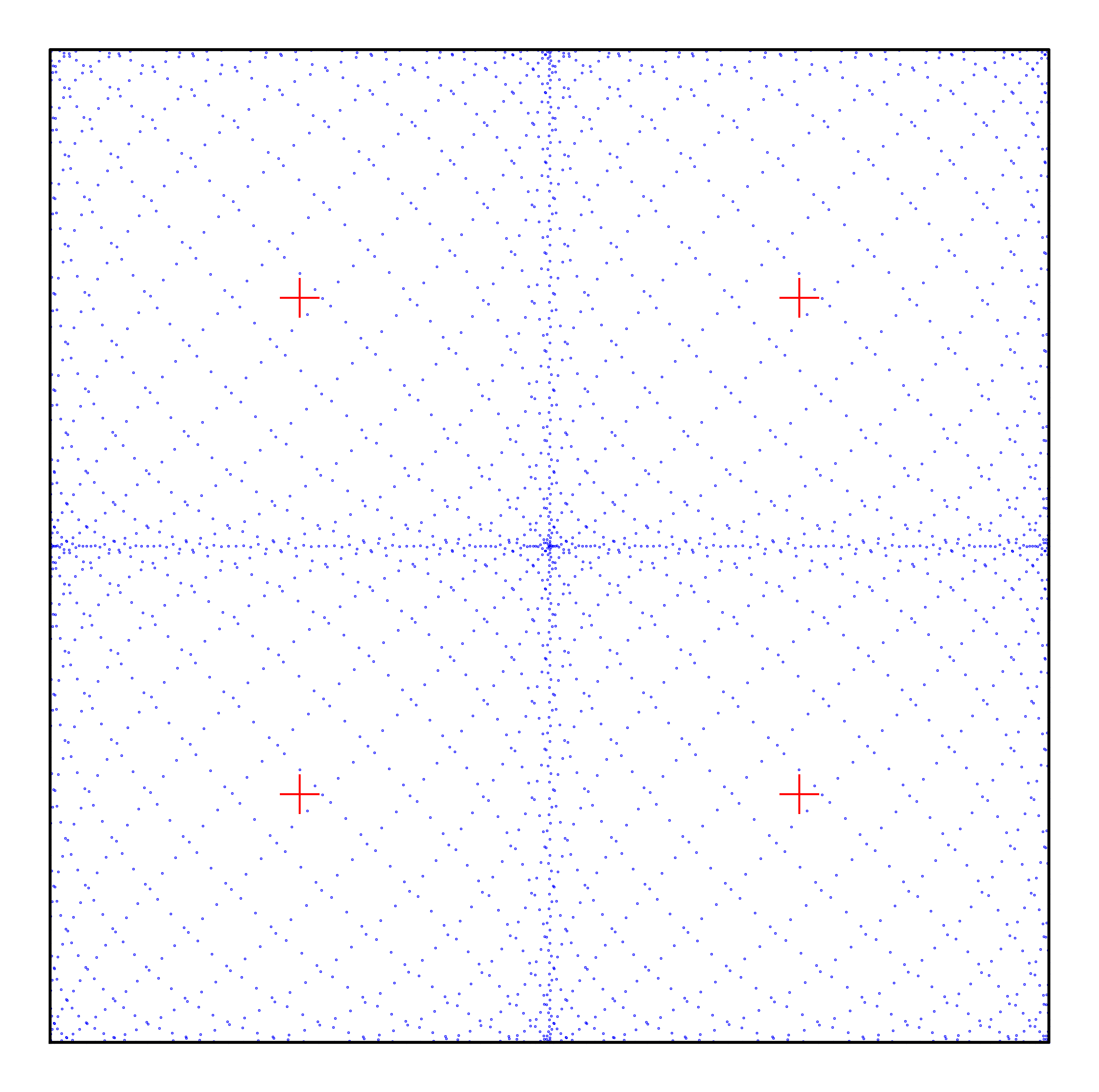}
		\caption{\centering\scriptsize $2\times 2$}
		\label{subfig:scans:22}
	\end{subfigure}
	\begin{subfigure}[t]{\imratio\linewidth}
		\includegraphics[width=\linewidth]{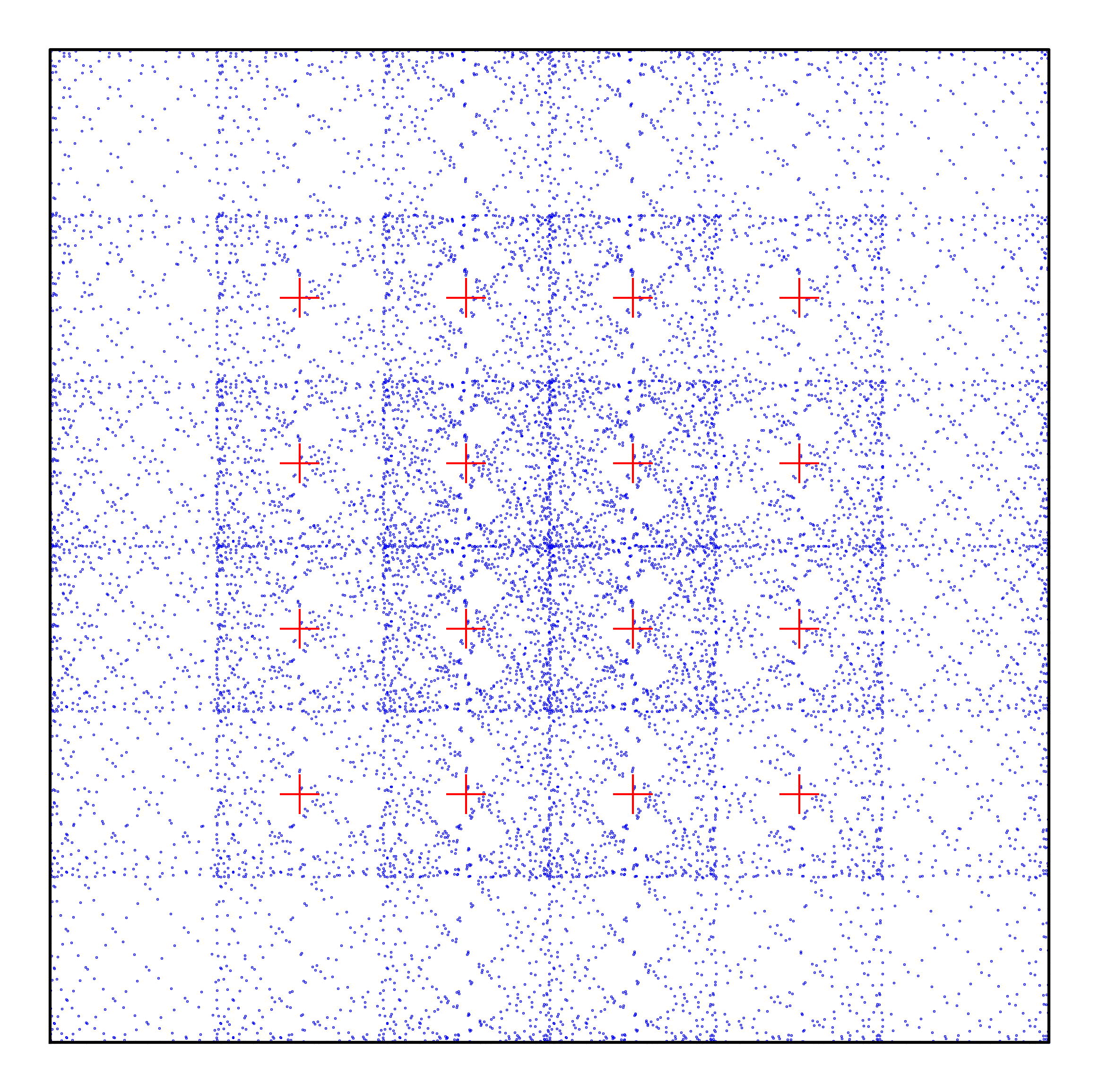}
		\caption{\centering\scriptsize $4\times 4$}
		\label{subfig:scans:44}
	\end{subfigure}
	\begin{subfigure}[t]{\imratio\linewidth}
		\includegraphics[width=\linewidth]{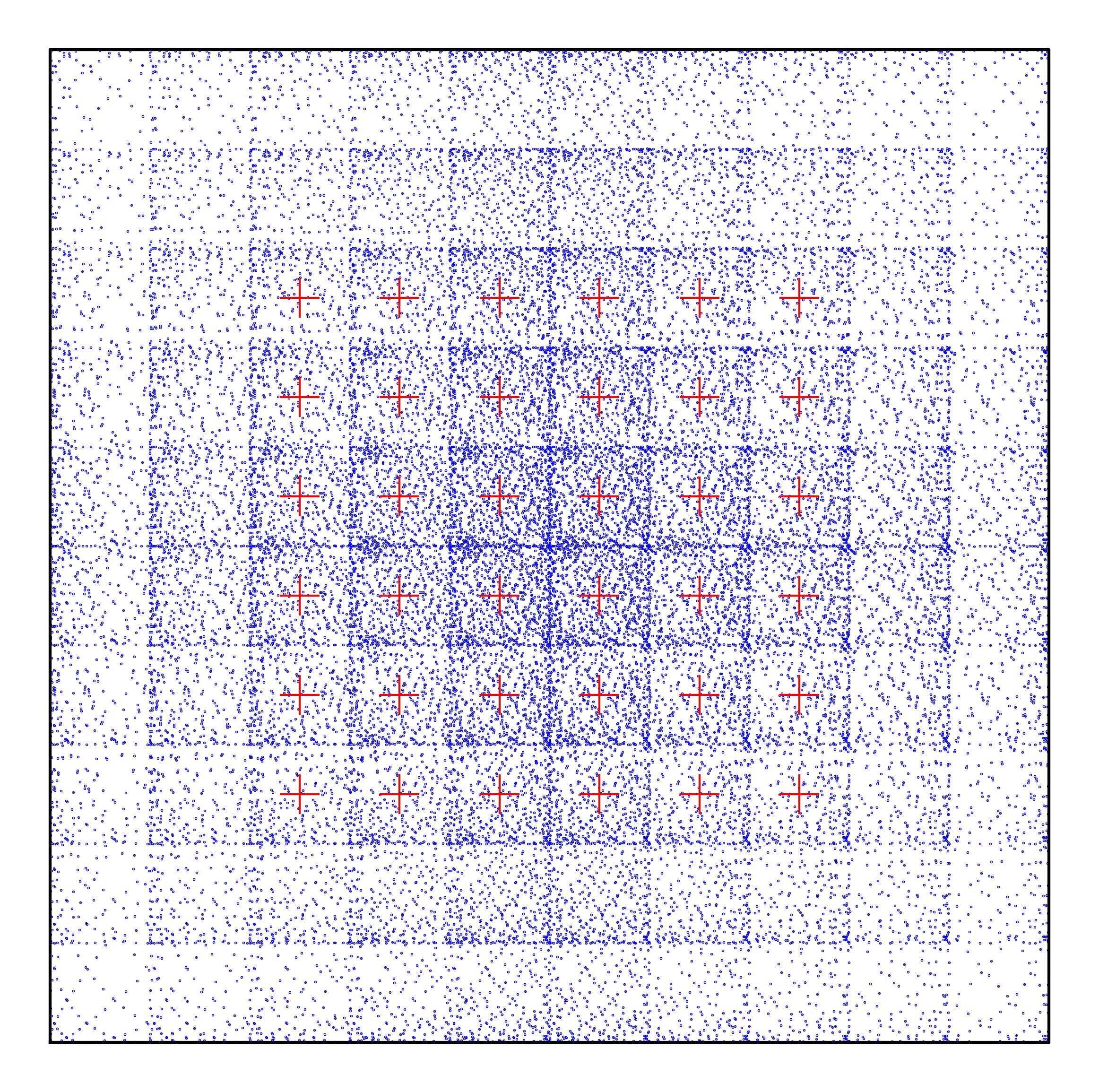}
		\caption{\centering\scriptsize $6\times 6$}
		\label{subfig:scans:66}
	\end{subfigure}
	\par\medskip
	\begin{subfigure}[t]{\imratio\linewidth}
		\includegraphics[width=\linewidth]{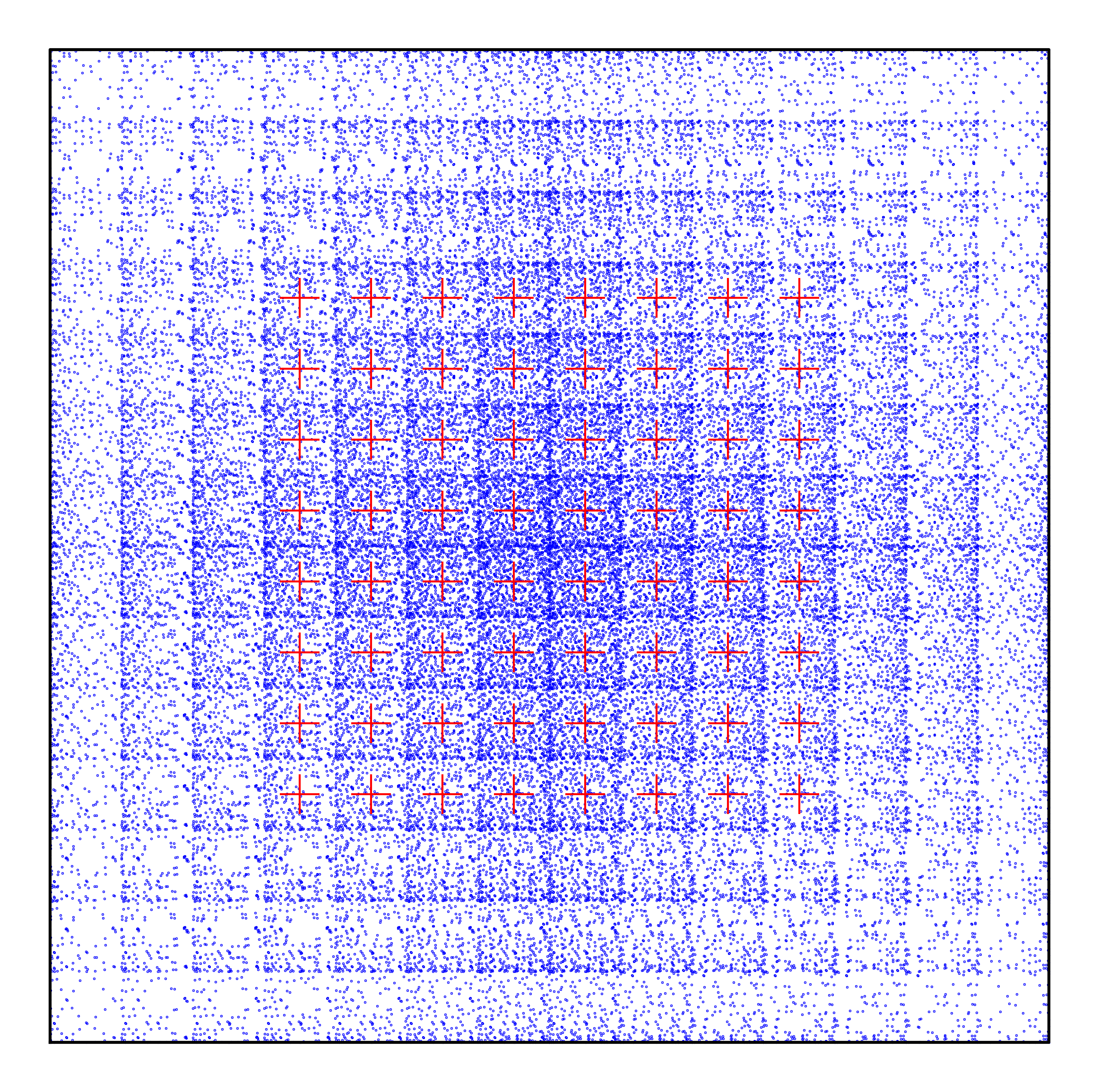}
		\caption{\centering\scriptsize $8\times 8$}
		\label{subfig:scans:88}
	\end{subfigure}
	\begin{subfigure}[t]{\imratio\linewidth}
		\includegraphics[width=\linewidth]{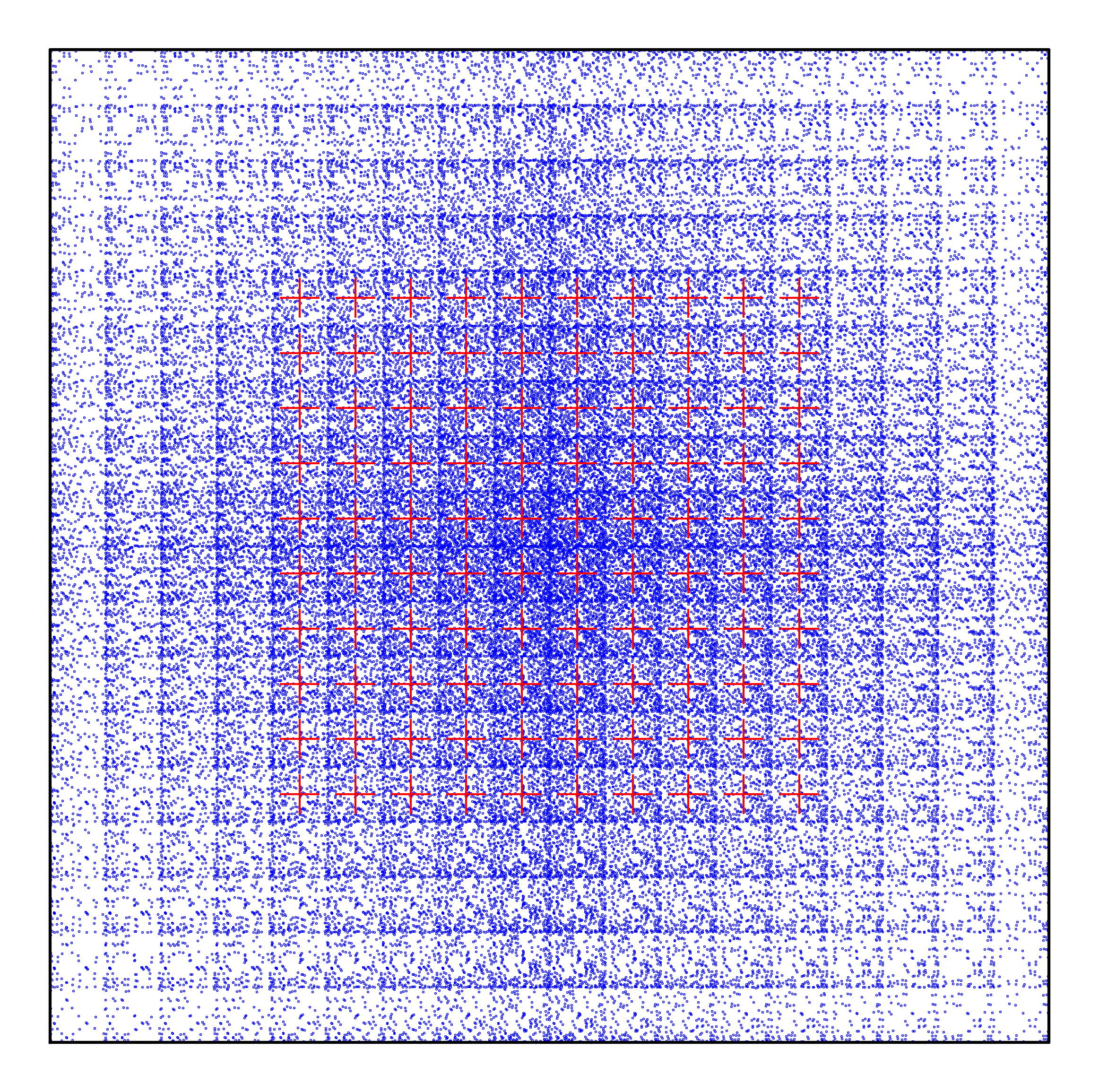}
		\caption{\centering\scriptsize $10\times 10$}
		\label{subfig:scans:1010}
	\end{subfigure}
	\begin{subfigure}[t]{\imratio\linewidth}
		\includegraphics[width=\linewidth]{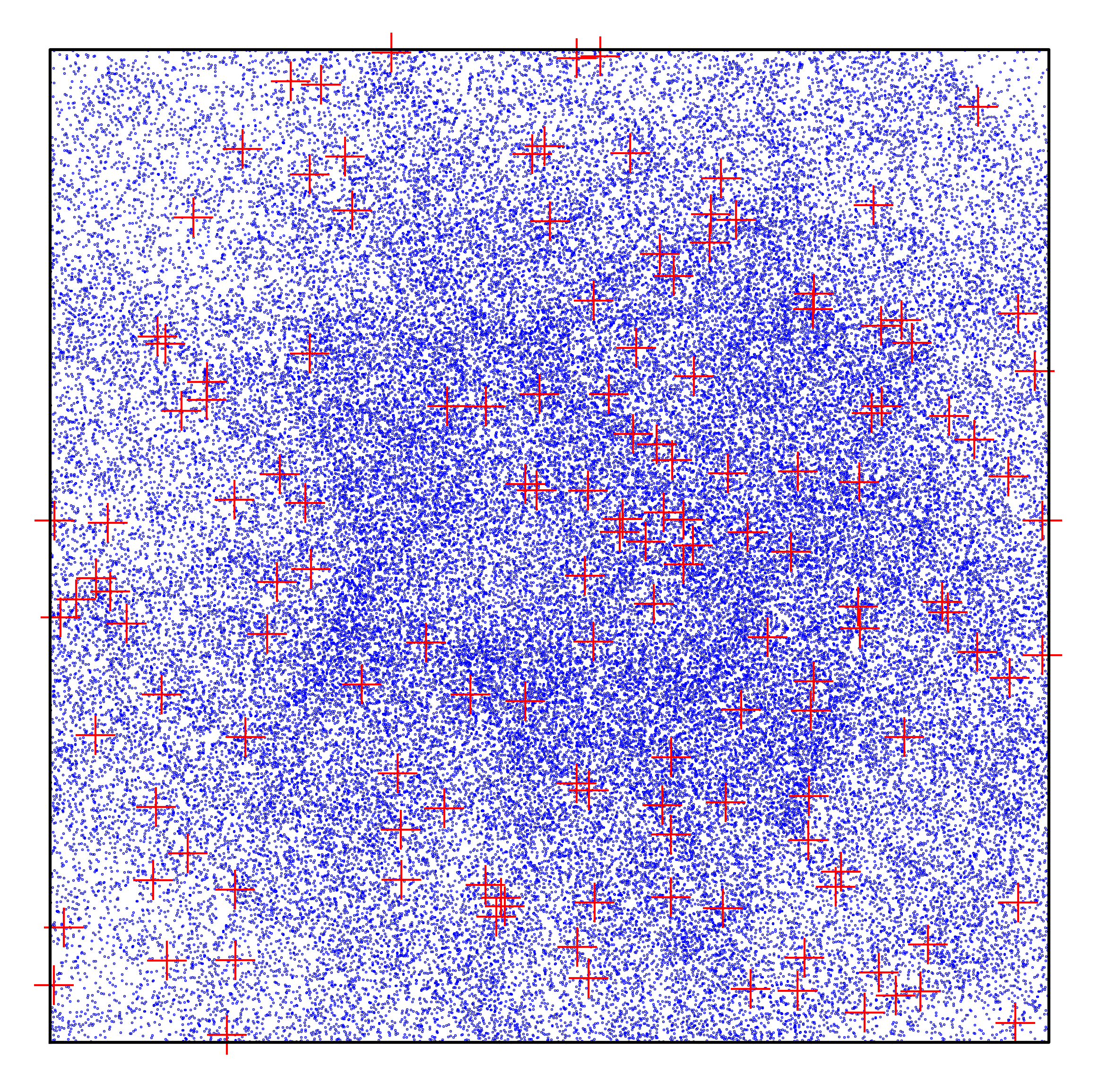}
		\caption{\centering\scriptsize Random}
		\label{subfig:scans:r}
	\end{subfigure}
	\caption{Sampling distribution of the scans in Experiment 1 and 3. Figures~\ref{subfig:scans:22},~\ref{subfig:scans:44},~\ref{subfig:scans:66},~\ref{subfig:scans:88} and~\ref{subfig:scans:1010} show the increased density in the sampling locations by performing standard multi-patching on the same region $\Omega=[-2,2]^2$ with FoV of amplitudes $A_x=A_y=1$. In particular, we have represented in blue the data points $r_k$ for each scan $\xi$, and in red the centers $b_\xi$ of the FoVs in each scan $\xi$. In Fig.~\ref{subfig:scans:r} we show the sampling locations obtained using 143 scans with random offset vectors and rotation angles. Random sampling has the property that the points tend to be more uniformly distributed, compared with standard multi-patching, where the region close to the boundary is always less sampled than the center of $\Omega$. Comparing Fig.~\ref{subfig:exp:vess:multi:10:rho} and~\ref{subfig:exp:vess:multi:r:rho}, illustrates the benefit of random sampling, as a more uniform distribution of the sampling allow to have more data close to the boundary. Indeed, increased sampling of the boundary is reflected in the reconstruction by a better reconstruction of the phantom close to the boundary.}
	\label{fig:scans}
\end{figure}

\subsection{Experiment 1: Standard Multi-Patching}

\def\imratio{0.19}
\begin{figure}[t]
	\centering
	\begin{subfigure}[t]{\imratio\linewidth}
		\includegraphics[width=\linewidth]{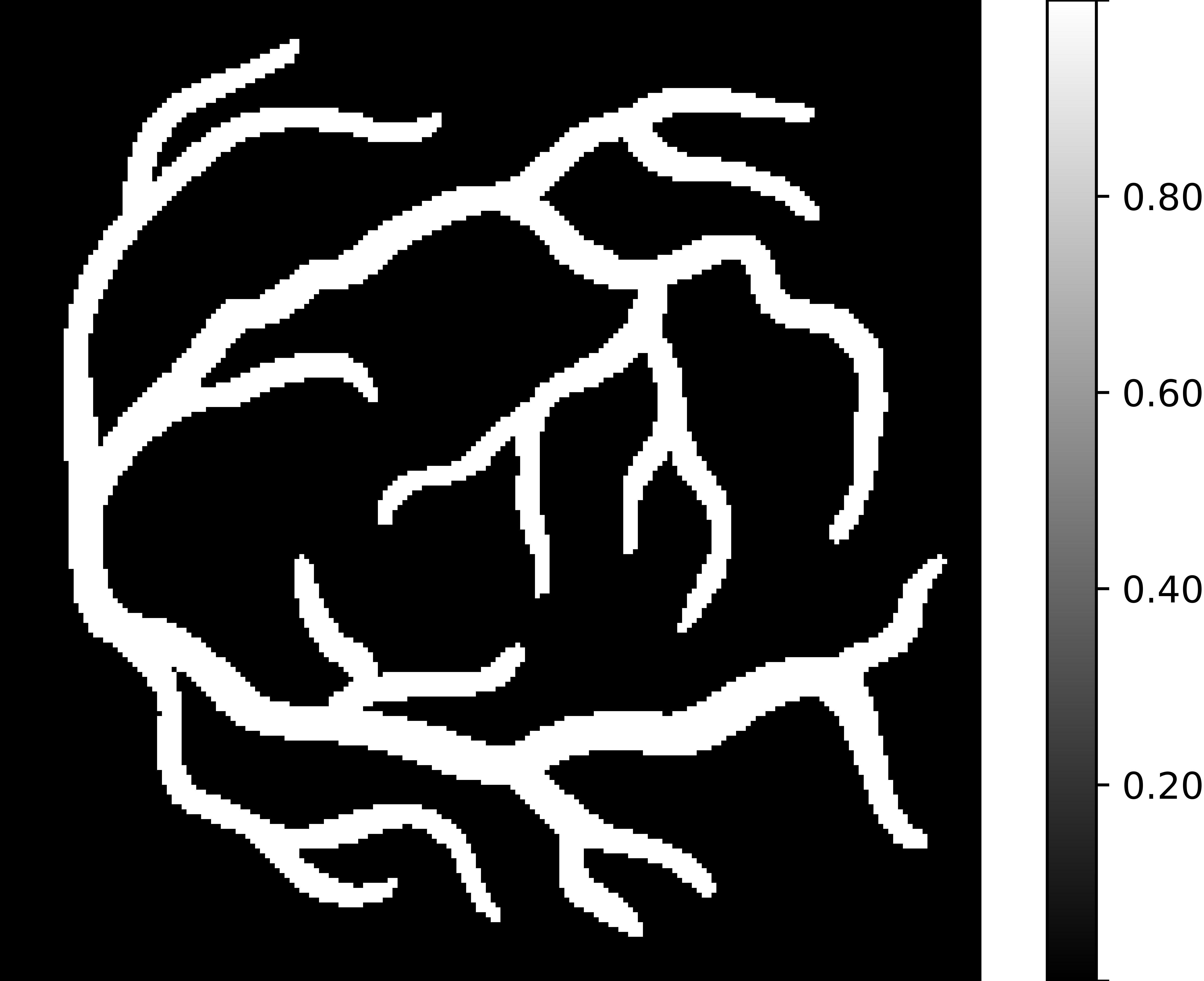}
		\caption{\centering\scriptsize $\rho_{\text{GT}}$.}
		\label{subfig:exp:vess:multi:gt}
	\end{subfigure}
	\hfil
	\begin{subfigure}[t]{\imratio\linewidth}
		\includegraphics[width=\linewidth]{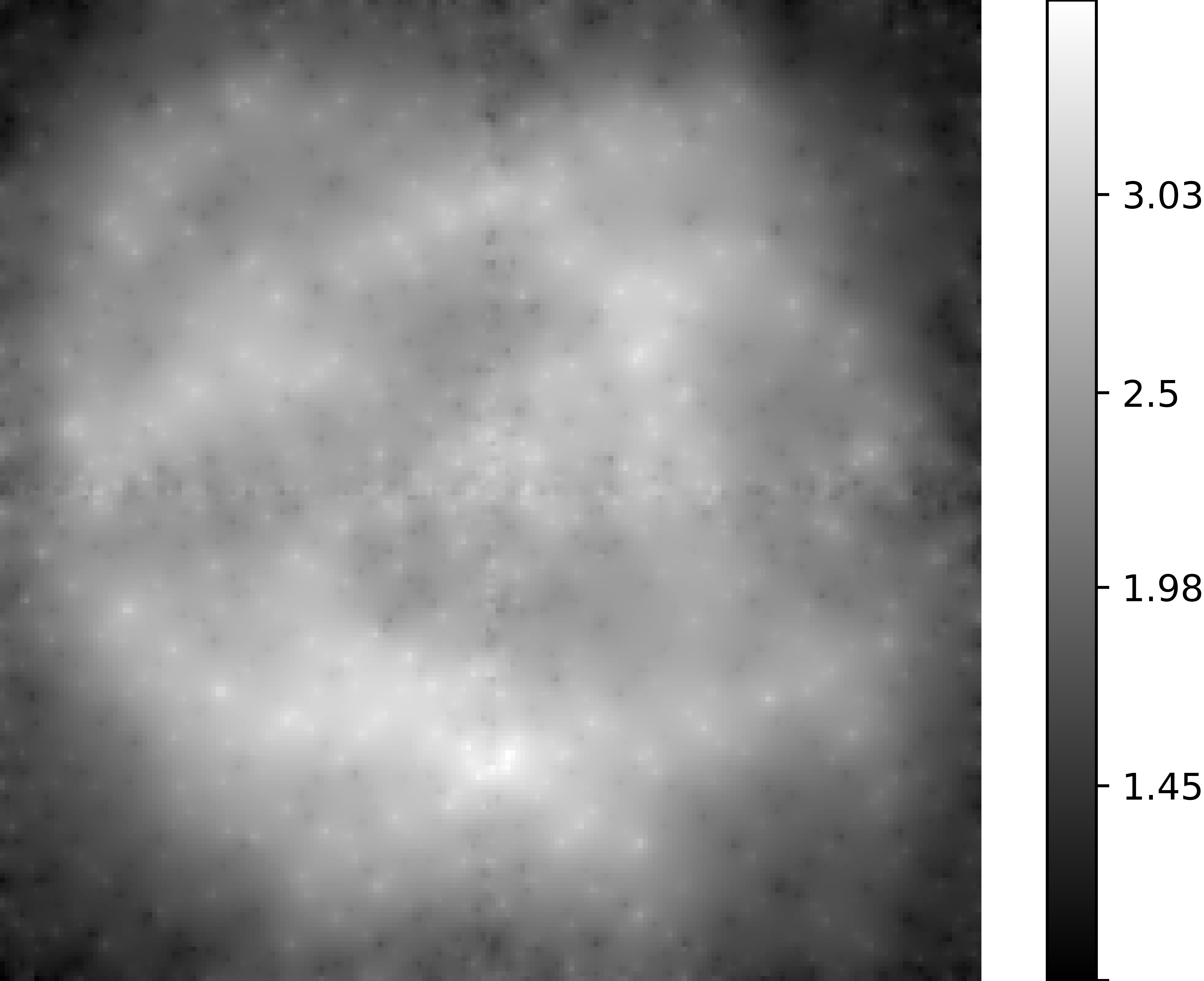}
		\caption{\centering\scriptsize $2\times 2 $, $u$ PSNR 26.21, SSIM 0.7853.}
		\label{subfig:exp:vess:multi:2:u}
	\end{subfigure}
	\hfil
	\begin{subfigure}[t]{\imratio\linewidth}
		\includegraphics[width={\linewidth}]{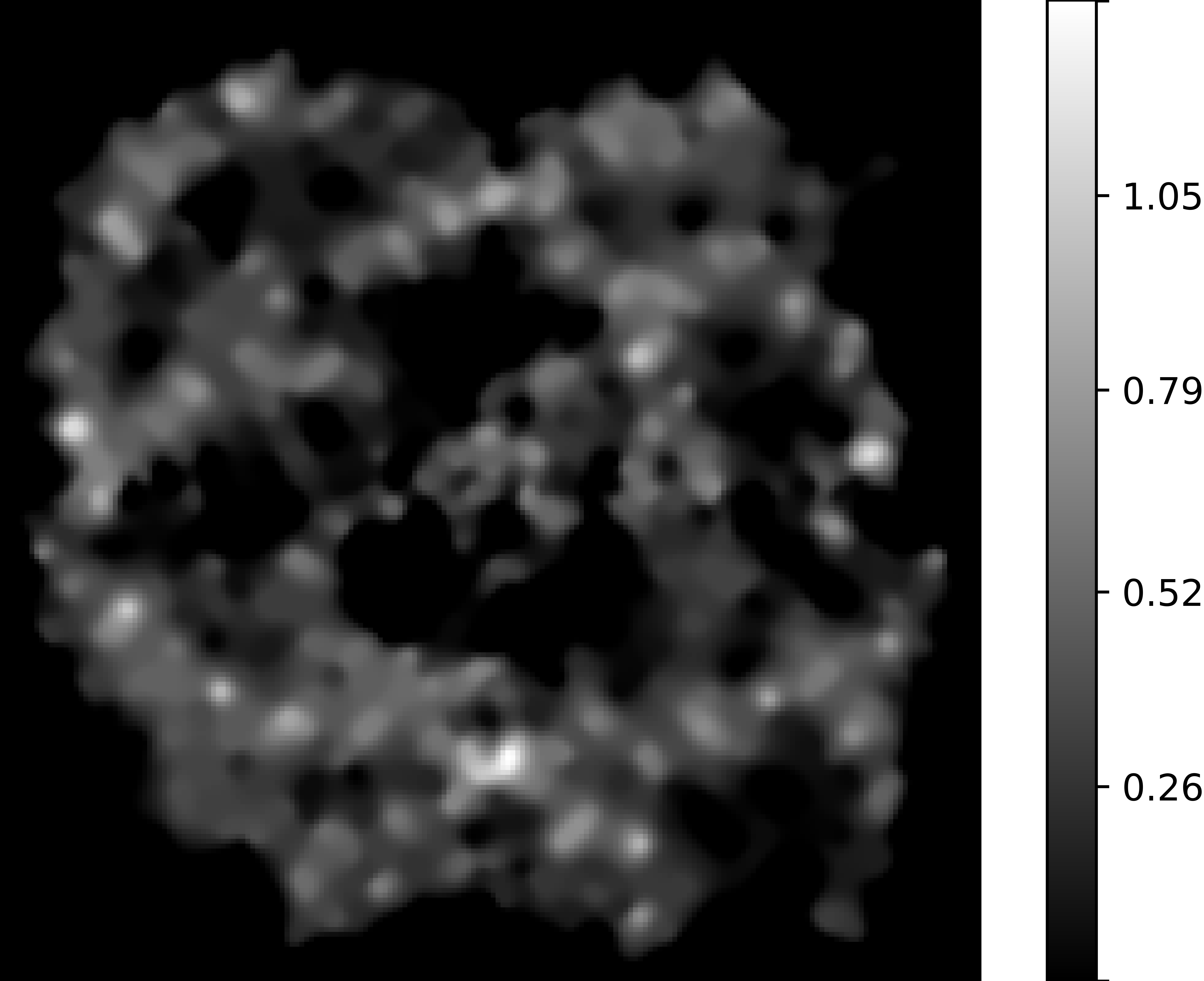}
		\caption{\centering\scriptsize $2\times 2$, $\rho$ PSNR 10.14, SSIM 0.3316.}
		\label{subfig:exp:vess:multi:2:rho}
	\end{subfigure}
	\hfil
	\begin{subfigure}[t]{\imratio\linewidth}
		\includegraphics[width=\linewidth]{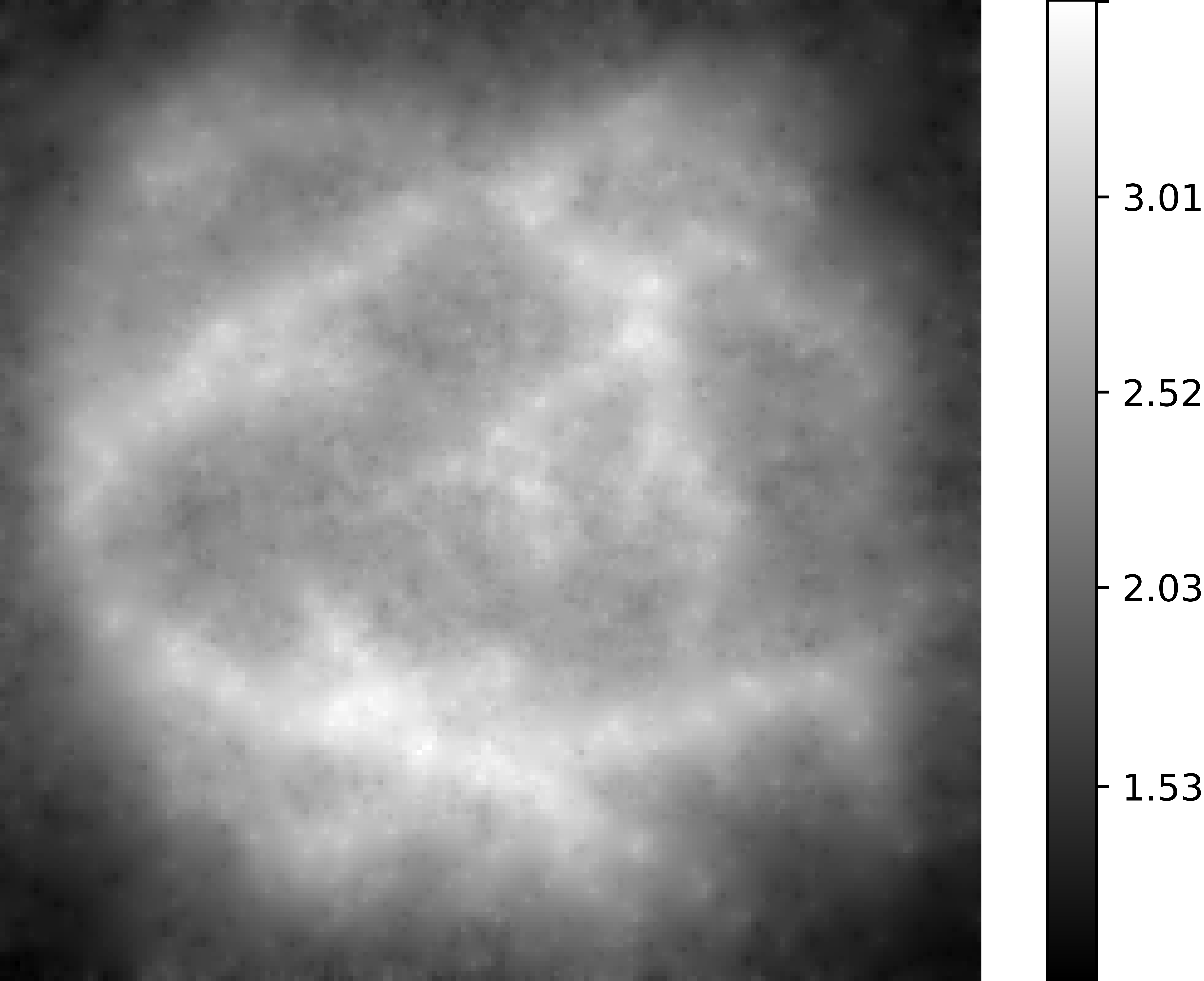}
		\caption{\centering\scriptsize $4\times 4 $, $u$ PSNR 28.99, SSIM 0.8307.}
		\label{subfig:exp:vess:multi:4:u}
	\end{subfigure}
	\hfil
	\begin{subfigure}[t]{\imratio\linewidth}
		\includegraphics[width={\linewidth}]{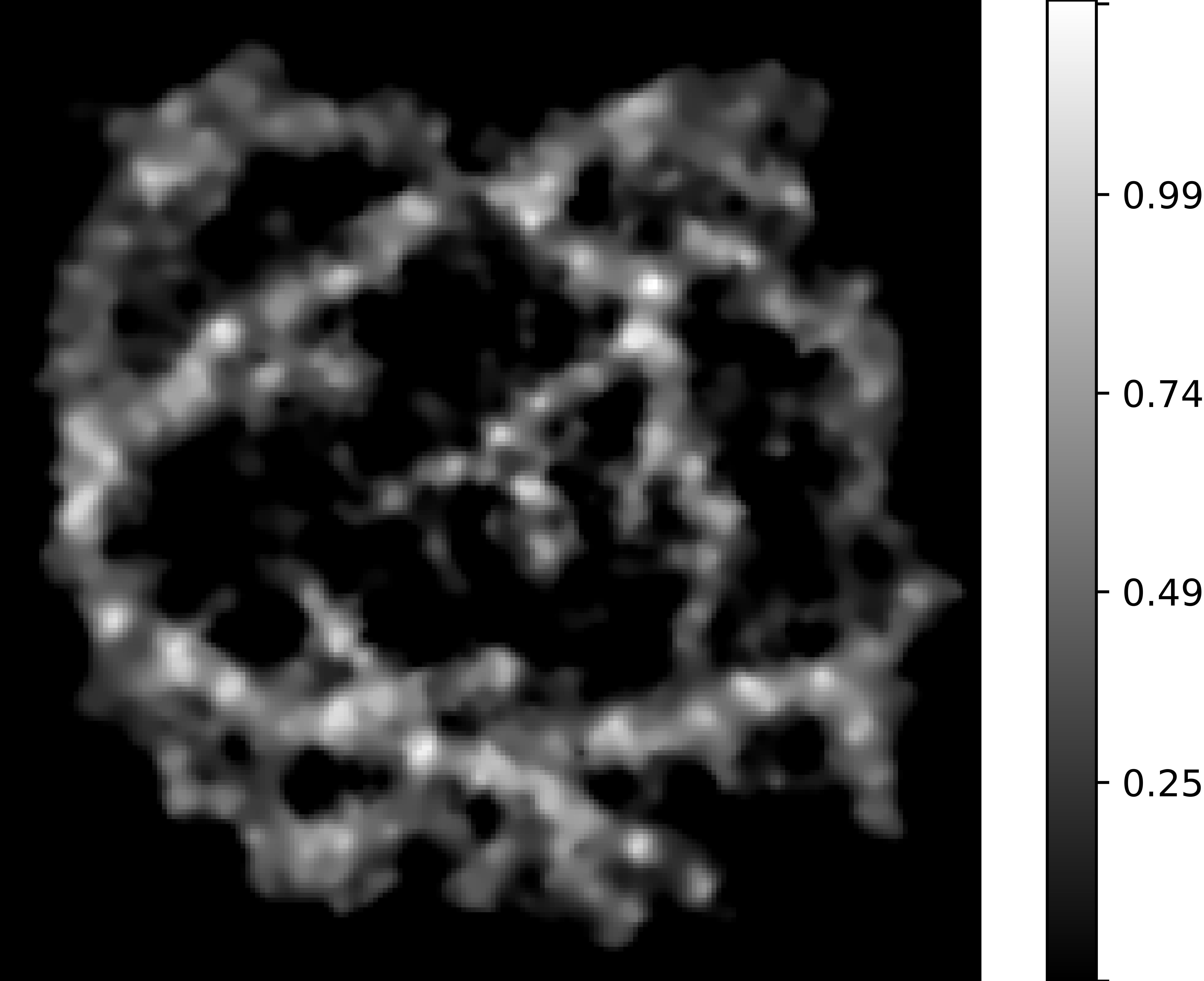}
		\caption{\centering\scriptsize $4\times 4$, $\rho$ PSNR 11.36, SSIM 0.4337.}
		\label{subfig:exp:vess:multi:4:rho}
	\end{subfigure}
	\par\medskip
	\begin{subfigure}[t]{\imratio\linewidth}
		\includegraphics[width=\linewidth]{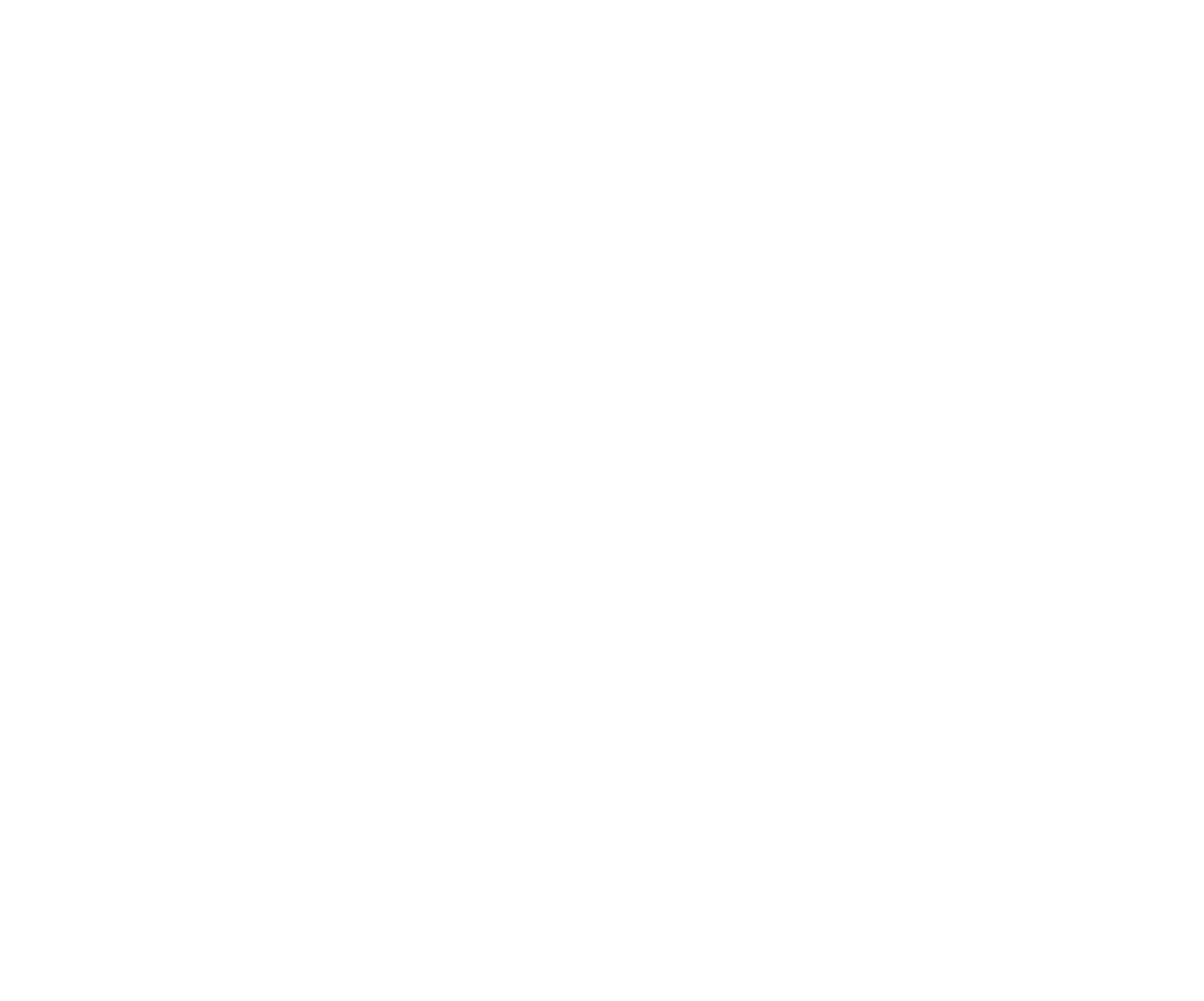}
	\end{subfigure}
	\hfil
	\begin{subfigure}[t]{\imratio\linewidth}
		\includegraphics[width=\linewidth]{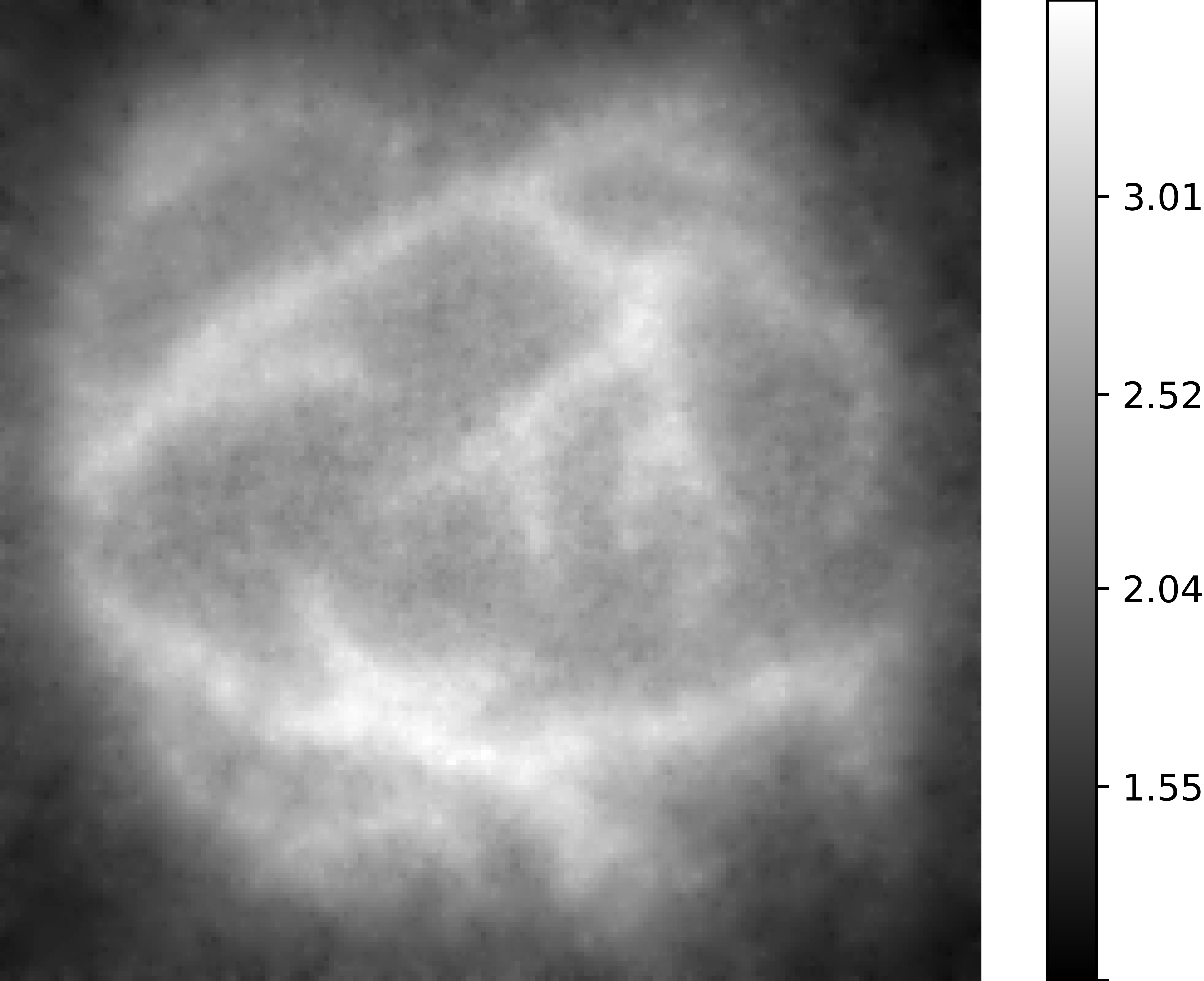}
		\caption{\centering\scriptsize $6\times 6 $, $u$ PSNR 30.10, SSIM 0.8521.}
		\label{subfig:exp:vess:multi:6:u}
	\end{subfigure}
	\hfil
	\begin{subfigure}[t]{\imratio\linewidth}
		\includegraphics[width={\linewidth}]{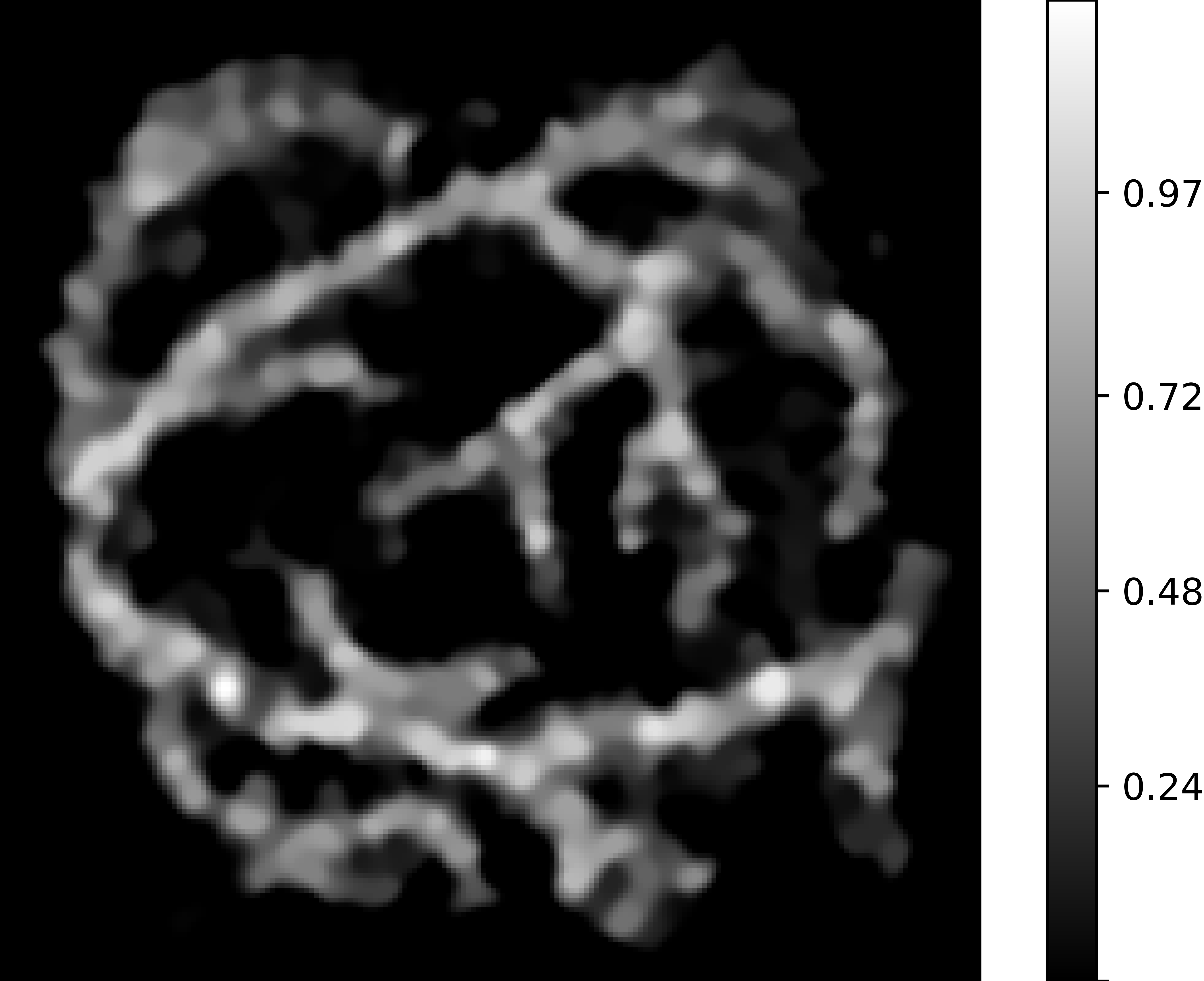}
		\caption{\centering\scriptsize $6\times 6$, $\rho$ PSNR 12.29, SSIM 0.5247.}
		\label{subfig:exp:vess:multi:6:rho}
	\end{subfigure}
	\hfil
	\begin{subfigure}[t]{\imratio\linewidth}
		\includegraphics[width=\linewidth]{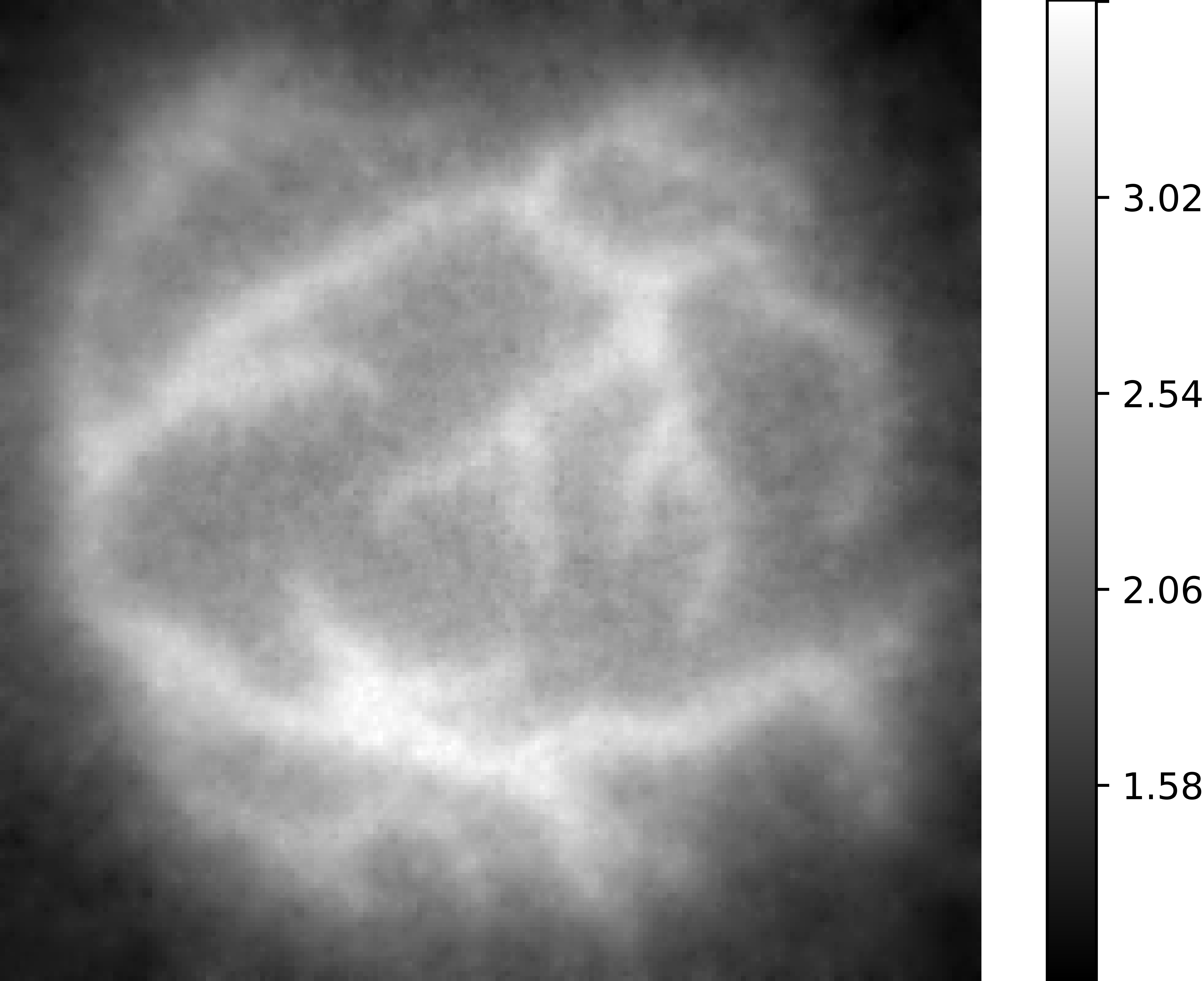}
		\caption{\centering\scriptsize $8\times 8 $, $u$ PSNR 31.51, 0.8703.}
		\label{subfig:exp:vess:multi:8:u}
	\end{subfigure}
	\hfil
	\begin{subfigure}[t]{\imratio\linewidth}
		\includegraphics[width={\linewidth}]{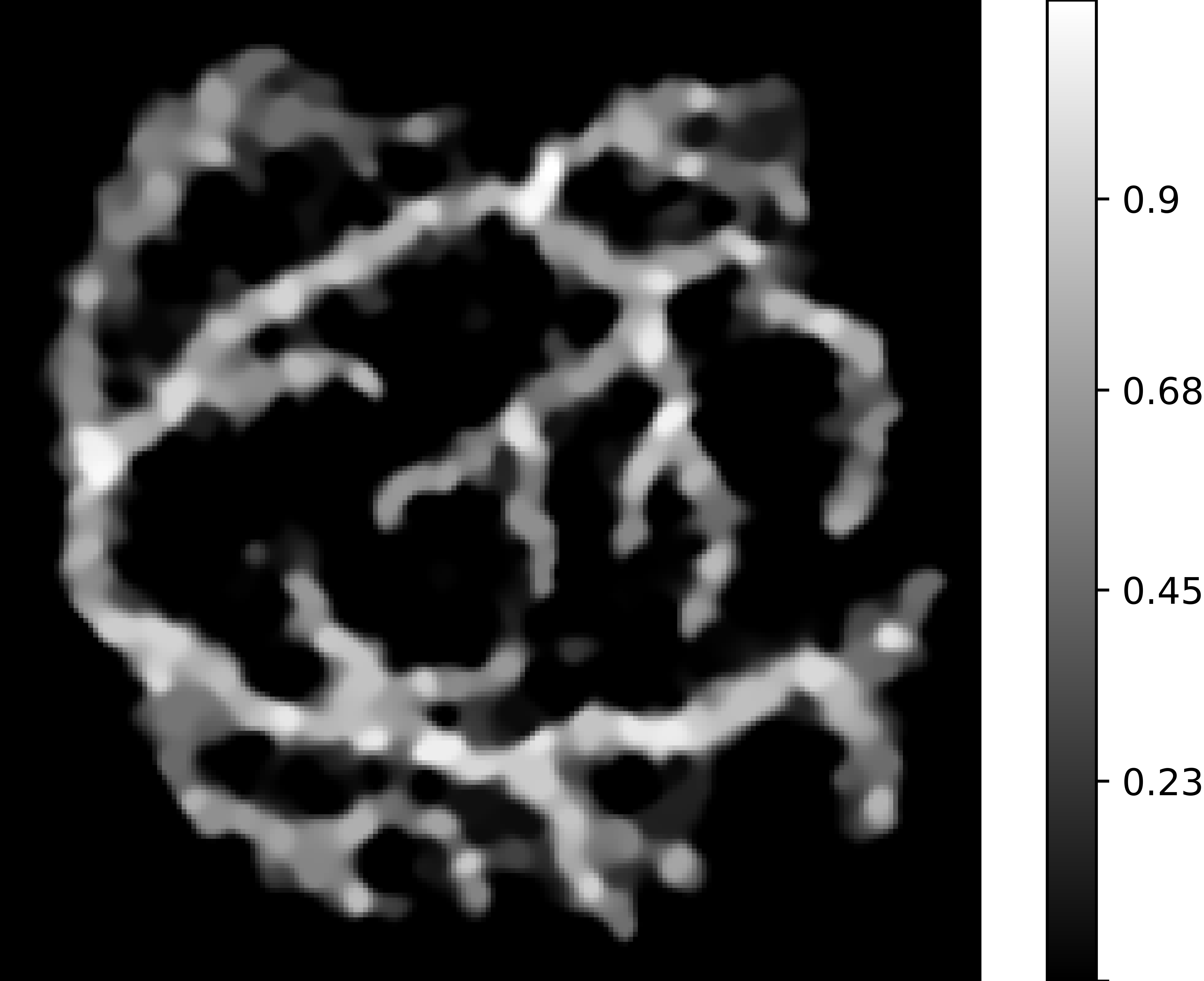}
		\caption{\centering\scriptsize $8\times 8$, $\rho$ PSNR 12.92, SSIM 0.5967.}
		\label{subfig:exp:vess:multi:8:rho}
	\end{subfigure}
	\par\medskip
	\begin{subfigure}[t]{\imratio\linewidth}
		\includegraphics[width=\linewidth]{Images/blank_dist.png}
	\end{subfigure}
	\hfil
	\begin{subfigure}[t]{\imratio\linewidth}
		\includegraphics[width=\linewidth]{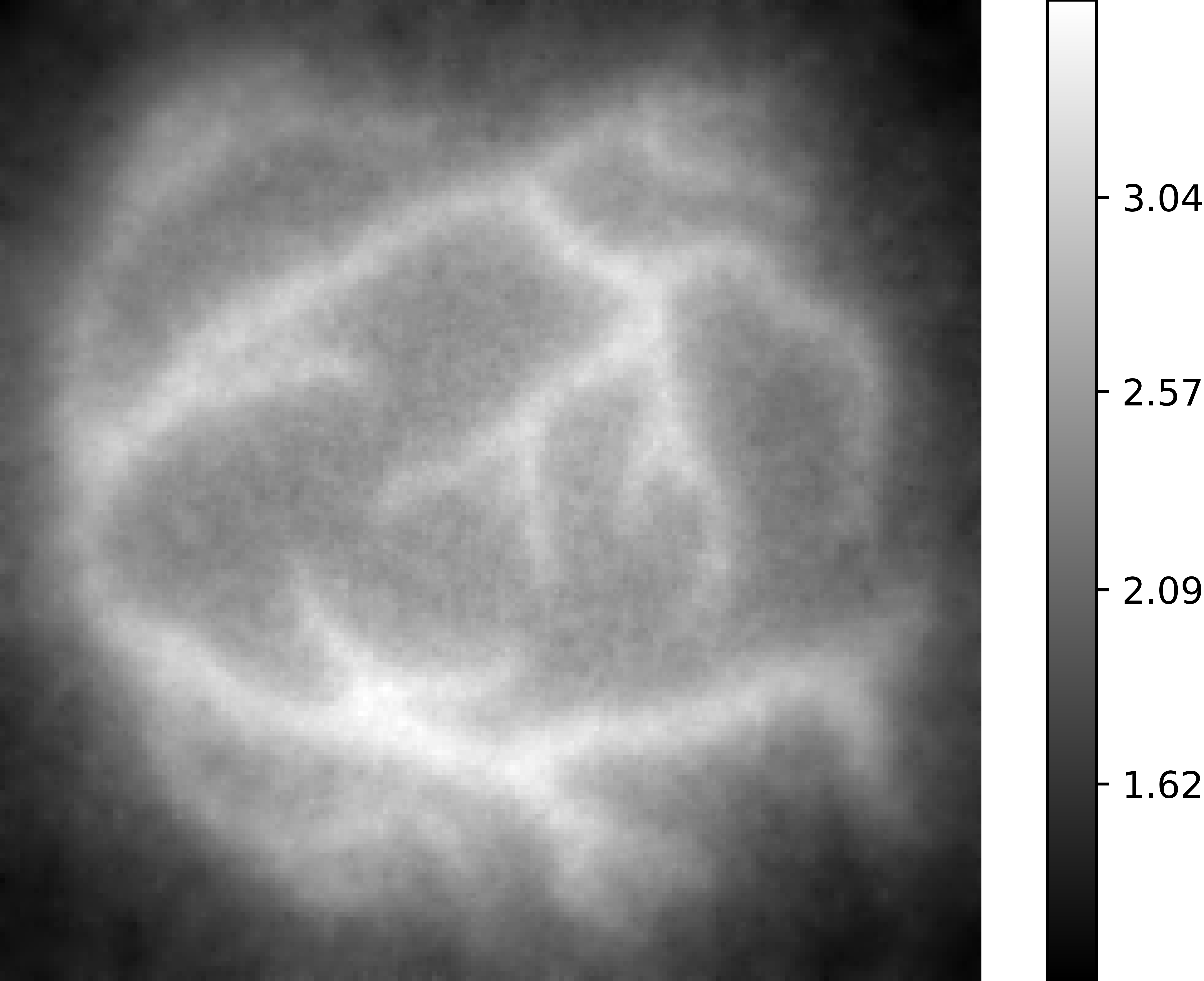}
		\caption{\centering\scriptsize $10\times 10 $, $u$ PSNR 32.00, SSIM 0.8857.}
		\label{subfig:exp:vess:multi:10:u}
	\end{subfigure}
	\hfil
	\begin{subfigure}[t]{\imratio\linewidth}
		\includegraphics[width={\linewidth}]{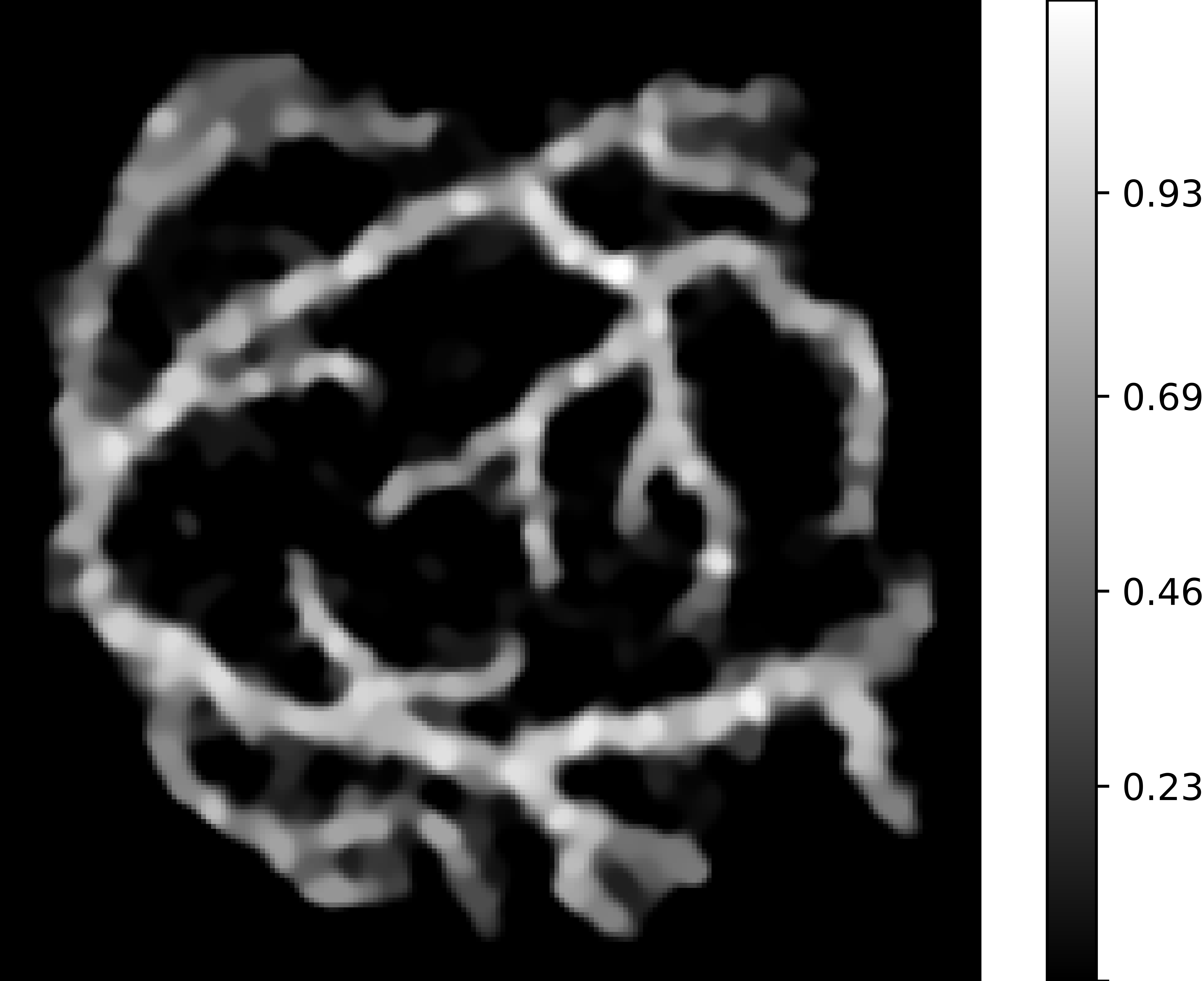}
		\caption{\centering\scriptsize $10\times 10$, $\rho$ PSNR 13.41, SSIM 0.6038.}
		\label{subfig:exp:vess:multi:10:rho}
	\end{subfigure}
	\begin{subfigure}[t]{\imratio\linewidth}
		\includegraphics[width=\linewidth]{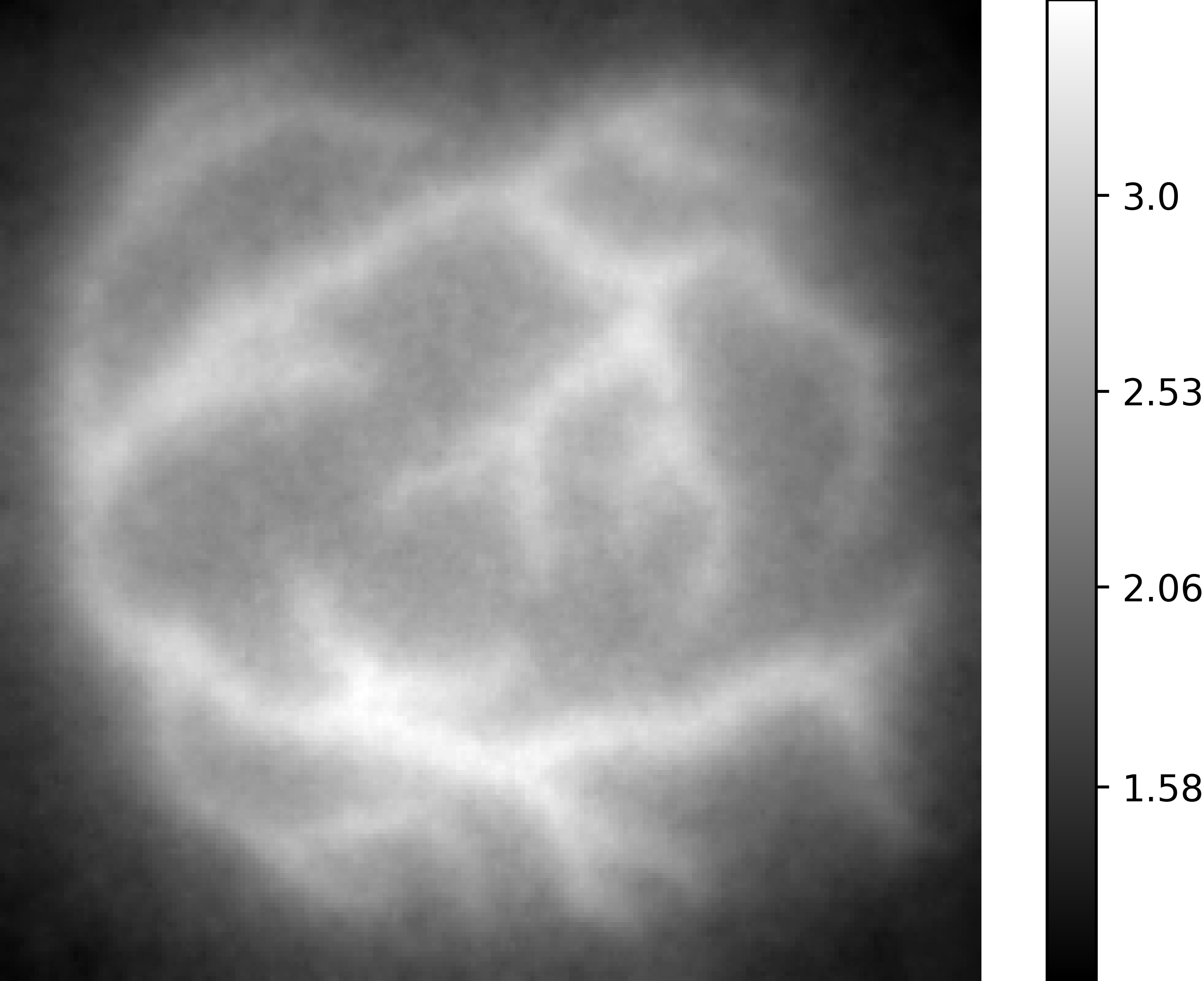}
		\caption{\centering\scriptsize Random, $u$ PSNR 34.53, SSIM 0.9238.}
		\label{subfig:exp:vess:multi:r:u}
	\end{subfigure}
	\hfil
	\begin{subfigure}[t]{\imratio\linewidth}
		\includegraphics[width={\linewidth}]{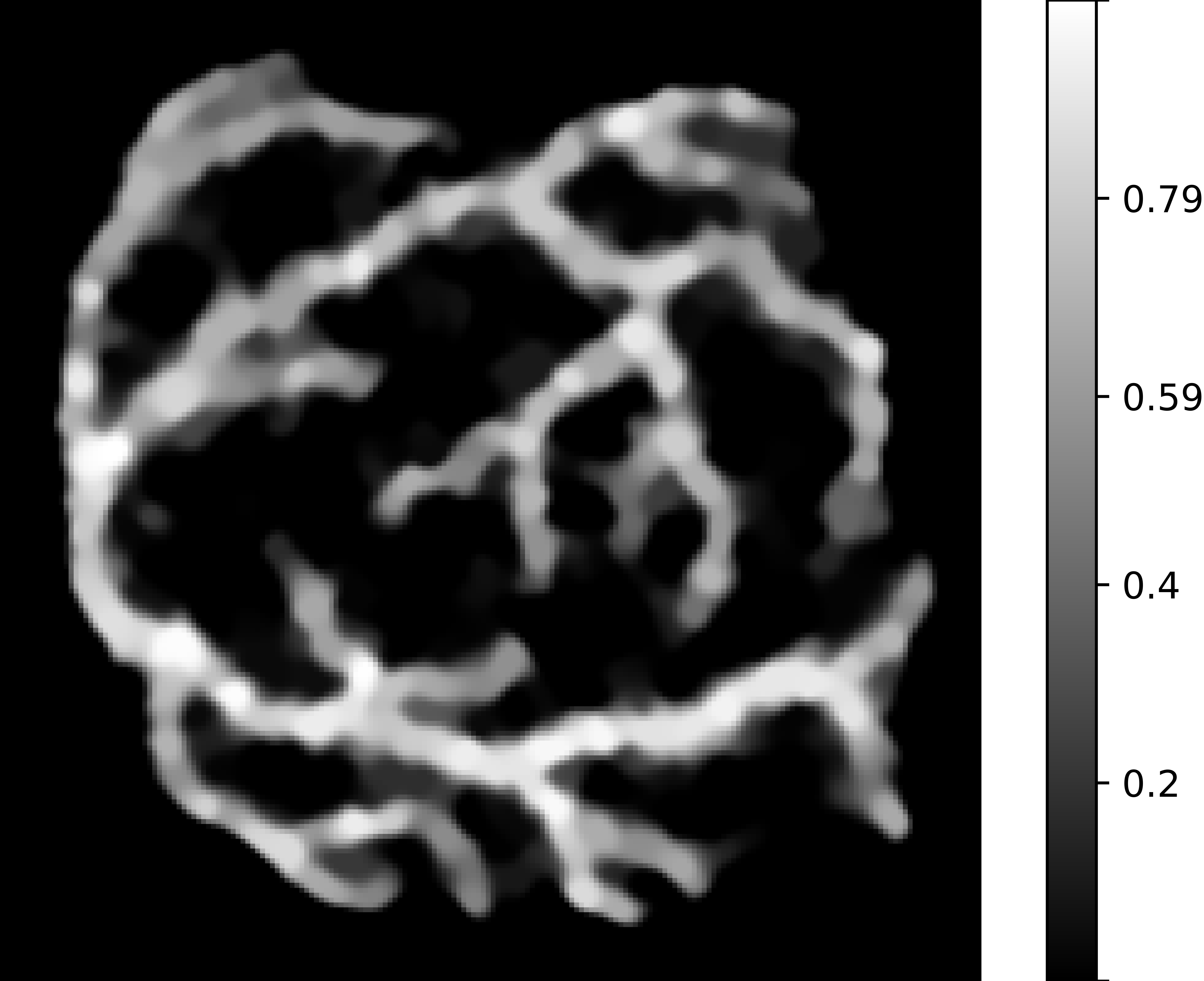}
		\caption{\centering\scriptsize Random, $\rho$ PSNR 13.33, SSIM 0.5903.}
		\label{subfig:exp:vess:multi:r:rho}
	\end{subfigure}
	
	\caption{An increased number of scans and denser samplings lead to better reconstructions. Figures~\ref{subfig:exp:vess:multi:2:u},~\ref{subfig:exp:vess:multi:4:u},~\ref{subfig:exp:vess:multi:6:u},~\ref{subfig:exp:vess:multi:8:u}and~\ref{subfig:exp:vess:multi:10:u} show the trace fields obtained with the first stage of the algorithm using the standard multi-patch scans in Fig.~\ref{subfig:scans:22},~\ref{subfig:scans:44},~\ref{subfig:scans:66},~\ref{subfig:scans:88} and~\ref{subfig:scans:1010}. Their reconstructions are shown in Figures~\ref{subfig:exp:vess:multi:2:rho},~\ref{subfig:exp:vess:multi:4:rho},~\ref{subfig:exp:vess:multi:6:rho},~\ref{subfig:exp:vess:multi:8:rho}and~\ref{subfig:exp:vess:multi:10:rho}. We observe that in case of a phantom that is bigger than the FoV, multi-patching can be used to produce a reconstruction after merging different scans that cover the region of interest (\ref{subfig:exp:vess:multi:2:u},~\ref{subfig:exp:vess:multi:2:rho}). Moreover, an overlap of scans, which produces a denser sampling by increasing the number of scans in the region (cf. Fig.~\ref{fig:scans}), improves the reconstruction quality. The last images (\ref{subfig:exp:vess:multi:r:u},~\ref{subfig:exp:vess:multi:r:rho}) show the result of Experiment 3, i.e., a reconstruction performed using the random samplings show in Fig.~\ref{subfig:scans:r}. For all reconstructions here we use a $200\times 200$ grid discretization of the region $\Omega = [-2,2]^2$.}
	\label{fig:exp:vess:multi:recs}
\end{figure}

In the first experiment we simulate multi-patch scans using the standard multi-patching with increasing number of scans (see Fig.~\ref{fig:scans}) of the vessel phantom in Fig.~\ref{subfig:exp:vess:multi:gt}. By standard multi-patching we mean that we cover the region $\Omega$ by equidistant overlapping scans. Consider $\Omega = [a,b]\times [c,d]$ and the field of view $[-A_x,A_x]\times [-A_y ,A_y]$, then if we want $I$ equidistant scans along the $x$ axis and $J$ equidistant scans along the $y$ axis, we will perform $\Xi = I\cdot J$ scans with offset vectors
\begin{equation}\label{eq:offset:std:multi}
	b_\xi \equiv b_{i,j} \coloneq (x_i , y_j )^T = (a + A_x + i\cdot d_x , c + A_y + j\cdot d_y )
\end{equation}
and angles $\alpha_{i,j} = 0$ for $i\in\lbrace 0,\dots , I-1\rbrace$ and $j\in\lbrace 0,\dots ,J-1\rbrace$, where $d_x = \frac{b-a-2A_x}{I-1}$ and $d_y = \frac{d-c-2A_x}{J-1}$ are the step sizes between consecutive centers of the FoVs along the $x$-axis and the $y$-axis, respectively. In our case we considered the vessel phantom in Fig.~\ref{subfig:exp:vess:multi:gt} in the region $\Omega=[-2,2]^2$ with FoV with amplitudes $A_x = A_y = 1$ and an increasing $I=J$ for $I\in\lbrace 2,4,6,8,10\rbrace$. We observe from the decreasing behavior of the optimal parameter $\lambda$ in Tab.~\ref{tab:exp}, that a higher sampling needs less inpainting from the regularizer in Eq.~\eqref{eq:first:discrete:reg}. This indeed confirms the intuitive explanation in Remark~\ref{rem:inpainting}. Moreover, the PSNR values increase for both stages as $I$ increases, supporting the use of standard multi-patching as a way to obtain improved reconstructions.

\subsection{Experiment 2: Ablation Study}

\def\imratio{0.24}
\begin{figure}[t]
	\centering
	\begin{subfigure}[t]{\imratio\linewidth}
		\includegraphics[width=\linewidth]{Images/vessel_gt.png}
		\caption{\centering\scriptsize $\rho_{\text{GT}}$.}
		\label{subfig:abl:gt}
	\end{subfigure}
	\hfil
	\begin{subfigure}[t]{\imratio\linewidth}
		\includegraphics[width=\linewidth]{Images/vessel_u_10.png}
		\caption{\centering\scriptsize $10\times 10 $, $u$ PSNR 32.00, SSIM 0.8857.}
		\label{subfig:abl:u}
	\end{subfigure}
	\hfil
	\begin{subfigure}[t]{\imratio\linewidth}
		\includegraphics[width={\linewidth}]{Images/vessel_rho_10.png}
		\caption{\centering\scriptsize $10\times 10$, $\rho$ with priors, PSNR 13.41, SSIM 0.6038.}
		\label{subfig:abl:rho}
	\end{subfigure}
	\hfil
	\begin{subfigure}[t]{\imratio\linewidth}
		\includegraphics[width={\linewidth}]{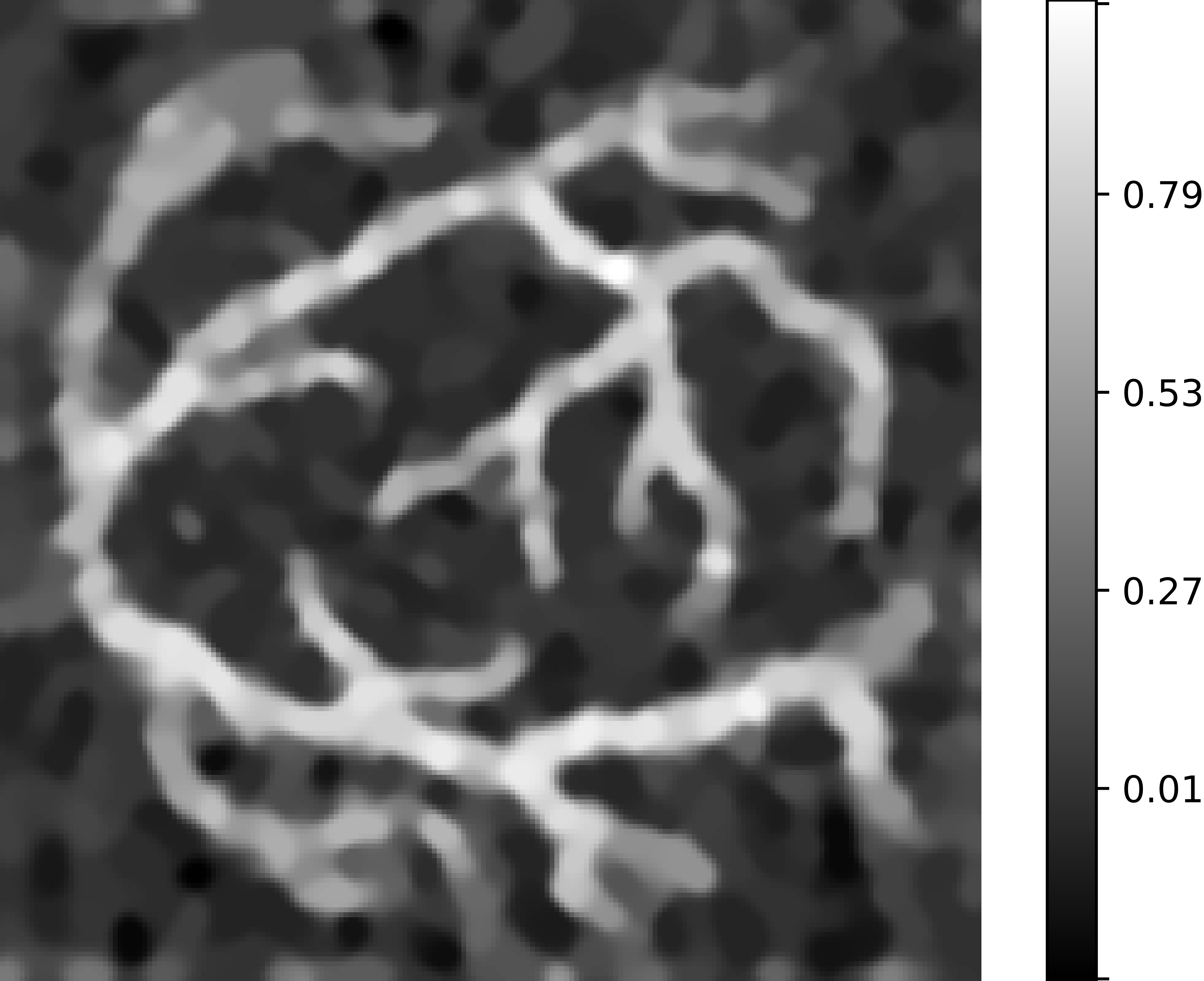}
		\caption{\centering\scriptsize $10\times 10$, $\rho$ without priors, PSNR 12.73, SSIM 0.2811.}
		\label{subfig:abl:rho:land}
	\end{subfigure}
	
	\caption{Comparison of the deconvolution in the second stage with the TV Smooth regularizer with and without (Fig.~\ref{subfig:abl:rho:land}) positivity constraints and sparsity promoting prior (Experiment 2). Both deconvolution are performed with the same trace field $u$ in Fig.~\ref{subfig:abl:u}, obtained performing the $10\times 10$ multi-patch scan in Fig.~\ref{subfig:scans:1010} of the vessel phantom in Fig.~\ref{subfig:abl:gt}. In particular, Fig.~\ref{subfig:abl:rho} is the reconstruction obtained using the functional with added priors in Eq.~\ref{eq:tv:smooth:discrete}, while Fig.~\ref{subfig:abl:rho:land} is obtained using the functional in Eq.~\ref{eq:second:stage}, without added priors. The quality enhancing properties of the priors are clear. }
	\label{fig:abl}
\end{figure}

In Experiment 2, we show the benefits of considering the functional in Eq.~\eqref{eq:second:pos:sparsity} with positivity constraint and sparsity enforcing prior, as opposed to the functional in Eq.~\eqref{eq:second:stage}. In particular, we consider the trace field $u$ obtained scanning the vessel phantom in Experiment 1, using a $10\times 10$ standard multi-patch scan (Fig.~\ref{subfig:scans:1010}) and performing the first stage of the algorithm. Then we compare the reconstructions obtained with Alg.\ref{alg:lasso} (Fig.~\ref{subfig:abl:rho}) and the reconstruction obtained by minimizing the functional in Eq.~\ref{eq:second:stage} (Fig.~\ref{subfig:abl:rho:land}). The input for both algorithms is the same reconstructed trace field $u$ with Alg.~\ref{alg:first}. In particular, we minimize the functional without further priors
\begin{equation}\label{eq:funcitonal:nospar}
	E[\rho ] = \lVert K_h \hat{\rho}-u\rVert_2^2 + \mu R[\hat{\rho}]
\end{equation}
with the TV Smooth regularizer $R$ in Eq.~\eqref{eq:tv:smooth:discrete}, using the Landweber iteration with $\rho^{(0)}=u$:
\begin{equation}
	\rho^{(k+1)} = \rho^{(k)}-\gamma\nabla F[\rho^{(k)} ]
\end{equation}
where $\nabla F[\rho]$ is as in Eq.~\eqref{eq:grad:f} and $\gamma = 10^{-3}$. The optimal parameter is $\mu = 1.75\cdot 10^{-4}$. Comparing the reconstruction obtained imposing positivity and promoting sparsity in Fig.~\ref{subfig:abl:rho} with the solution obtained minimizing the functional in Eq.~\eqref{eq:funcitonal:nospar}, it is visible the enhancement in the quality of the reconstruction, in particular with respect to compression of noise where the distribution $\rho$ is zero.

\subsection{Experiment 3: Random Samplings and Flexibility of the Generalized Multi-Patching}

In Experiment 3, we aim at pushing the boundaries with regards to what is possible with the generalized multi-patching approach. In particular, we opt for random scans, i.e., we perform scans using a FoV of amplitudes $A_x =A_y = 1$ on the region $\Omega = [-2,2]^2$ according to the generalized trajectories in Sec.~\ref{sec:gen:multipatching}, where the offset vector $b_{\xi}$ and the angles $\alpha_{\xi}$ in Eq.~\eqref{eq:piecewise:b:t} are random variables with uniform distribution, i.e.,
\begin{equation}
	b_{\xi}\sim\text{U}(\Omega ) ,\quad \alpha_{\xi}\sim\text{U}\bigl ([0,2\pi ) \bigr ) ,\quad\text{for }\xi\in\lbrace 0,\dots ,\Xi -1\rbrace .
\end{equation}
Because the offset vectors $b_{\xi}$ are taken in the region $\Omega$, if they are close enough to the boundary, there will be sample points outside $\Omega$. In the reconstructions, we consider only the data in $\Omega$, meaning that a higher number of scans is necessary to achieve a comparable number of samples in $\Omega$, compared with the standard multi-patching in Experiment 1, where each scan is fully contained in $\Omega$ (see Fig.~\ref{fig:scans}). In our case, we set $\Xi = 143$, which lead to the sampling in Fig.~\ref{subfig:scans:r}. In this experiment we see the first example of over-scanning methods. In particular, if the off-set vectors fall close enough to the boundary of $\Omega$, a portion of the data points collected are positioned outside of $\Omega$, producing \emph{de facto} an over-scan. Because only the data points in $\Omega$ can serve as input for the reconstruction algorithm, the data points ourside of the reconstruction region $\Omega$ are not considered. The amount of data points in the region $\Omega$ for this experiment is 173,260 which is equivalent to about $106$ scans with 1632 points per scan. The reconstructions with the random scans for the first (resp. second) stage can be seen in Fig.~\ref{subfig:exp:vess:multi:r:u} (resp. Fig.~\ref{subfig:exp:vess:multi:r:rho}). We observe that taking random scans allow for a more uniform sampling in the region $\Omega$, allowing for better reconstruction of those parts of the phantom that are close to the boundary (conf. Fig.~\ref{subfig:exp:vess:multi:10:rho} and Fig.~\ref{subfig:exp:vess:multi:r:rho}).

\subsection{Experiment 4: Phantoms Inspired by the Open MPI Dataset}\label{sec:exp:MPI}

\def\imratio{0.32}
\begin{figure}[t]
	\centering
	\begin{subfigure}[t]{\imratio\linewidth}
		\includegraphics[width=\linewidth]{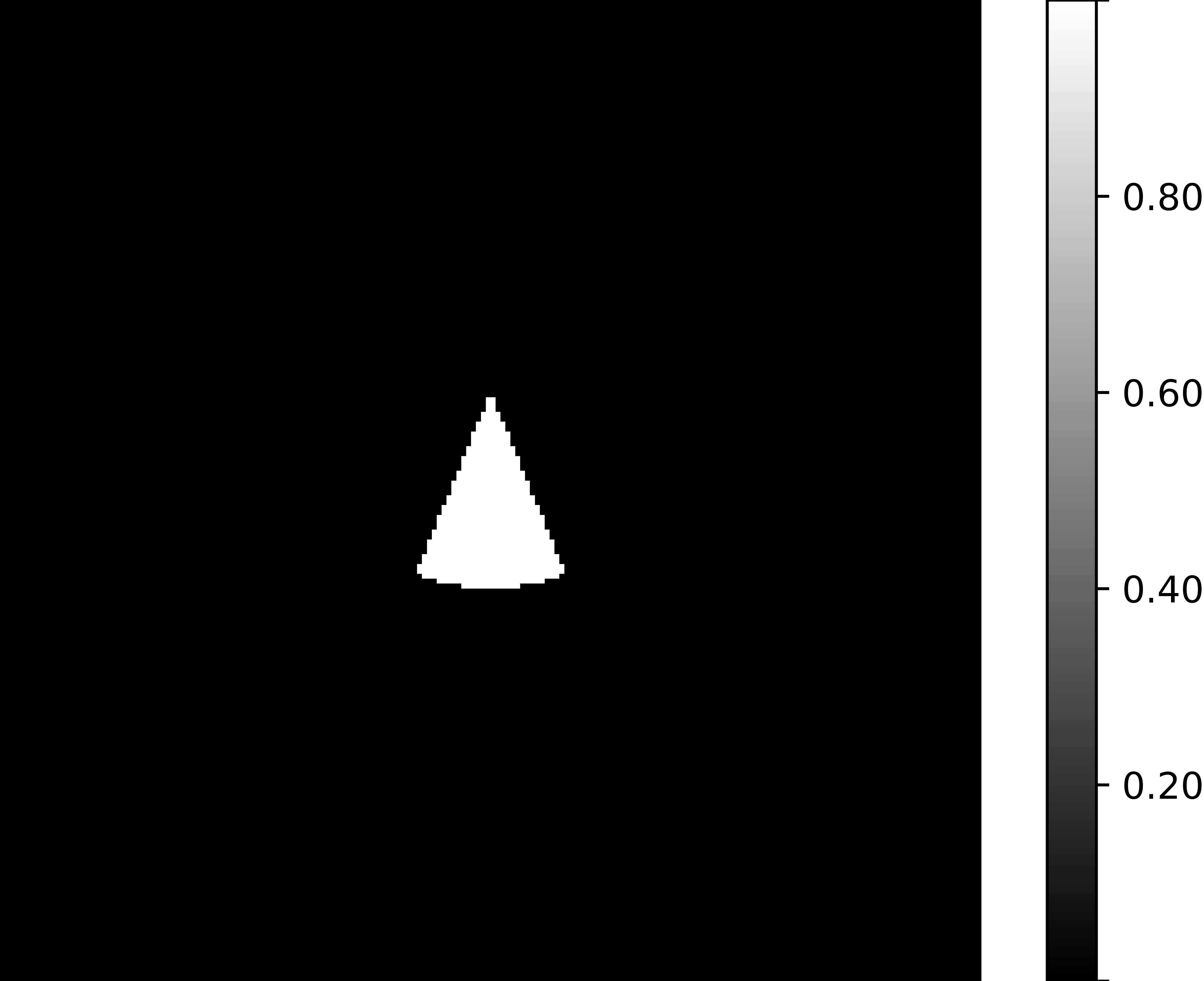}
		\caption{\centering\scriptsize Shape phantom $\rho_{\text{GT}}$.}
		\label{subfig:shape:gt}
	\end{subfigure}
	\hfil
	\begin{subfigure}[t]{\imratio\linewidth}
		\includegraphics[width=\linewidth]{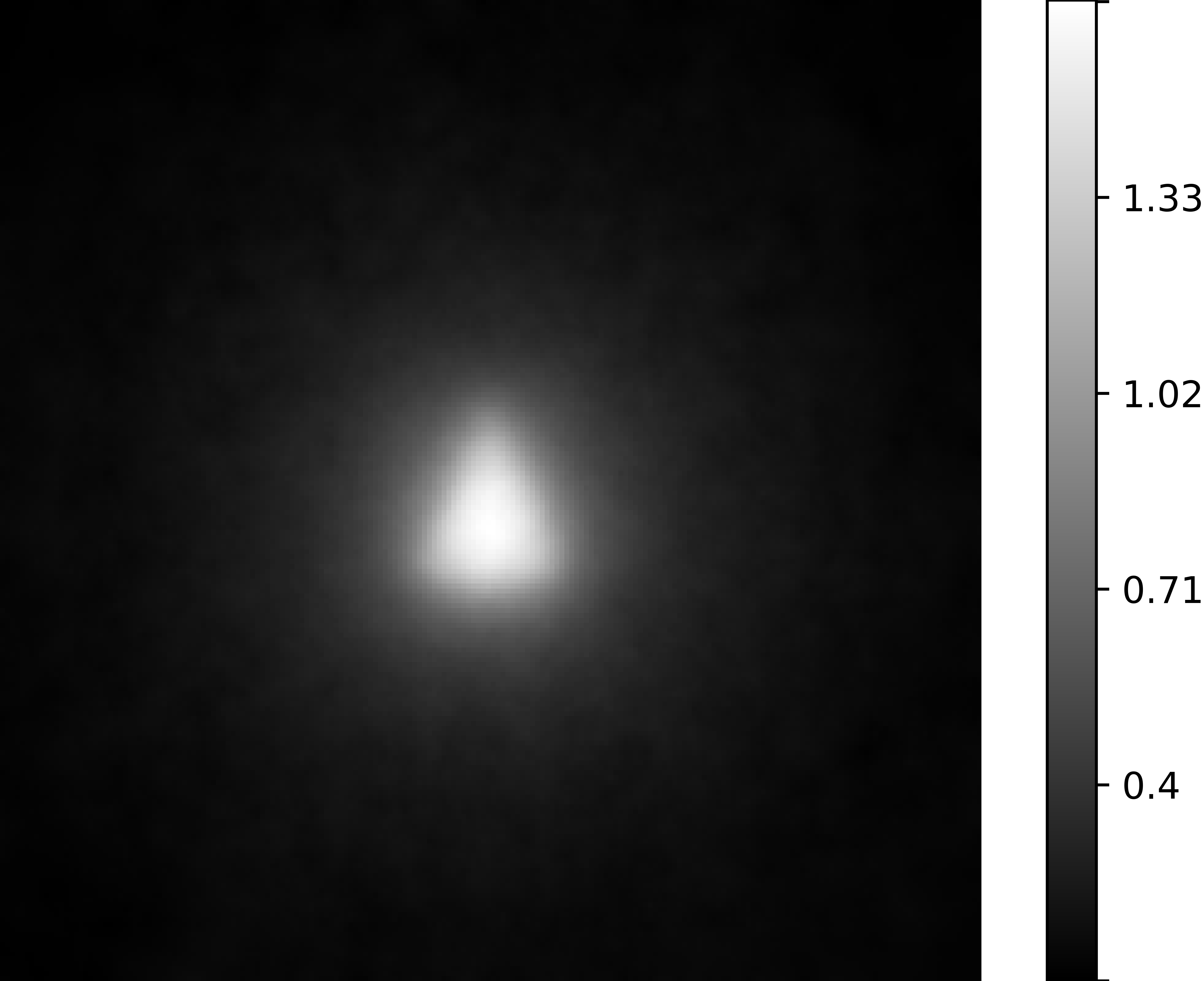}
		\caption{\centering\scriptsize $u$ PSNR 42.97, SSIM 0.9863.}
		\label{subfig:shape:u}
	\end{subfigure}
	\hfil
	\begin{subfigure}[t]{\imratio\linewidth}
		\includegraphics[width=\linewidth]{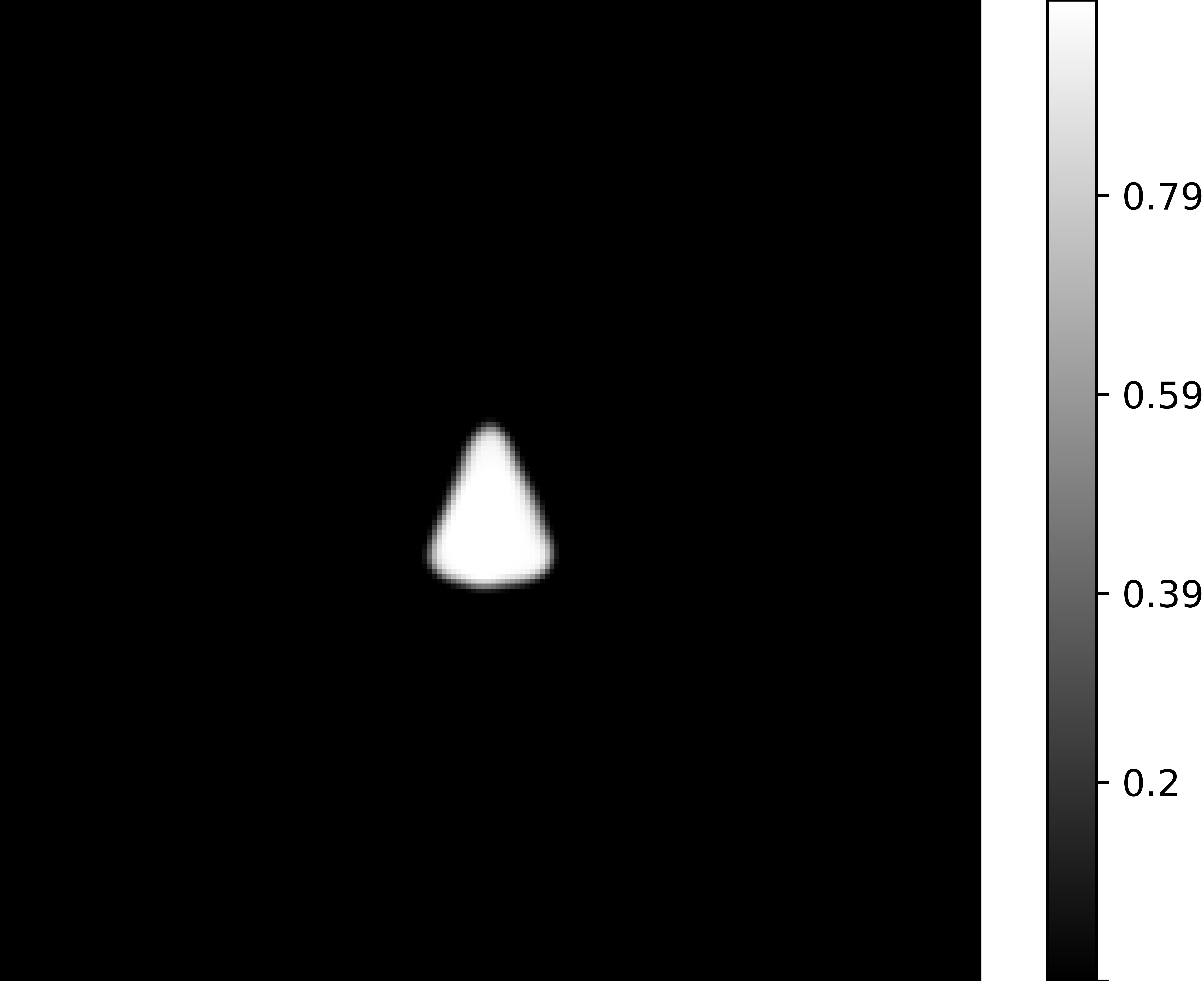}
		\caption{\centering\scriptsize $\rho$ PSNR 26.86, SSIM 0.9860.}
		\label{subfig:shape:rho}
	\end{subfigure}
	\par\medskip
	\begin{subfigure}[t]{\imratio\linewidth}
		\includegraphics[width=\linewidth]{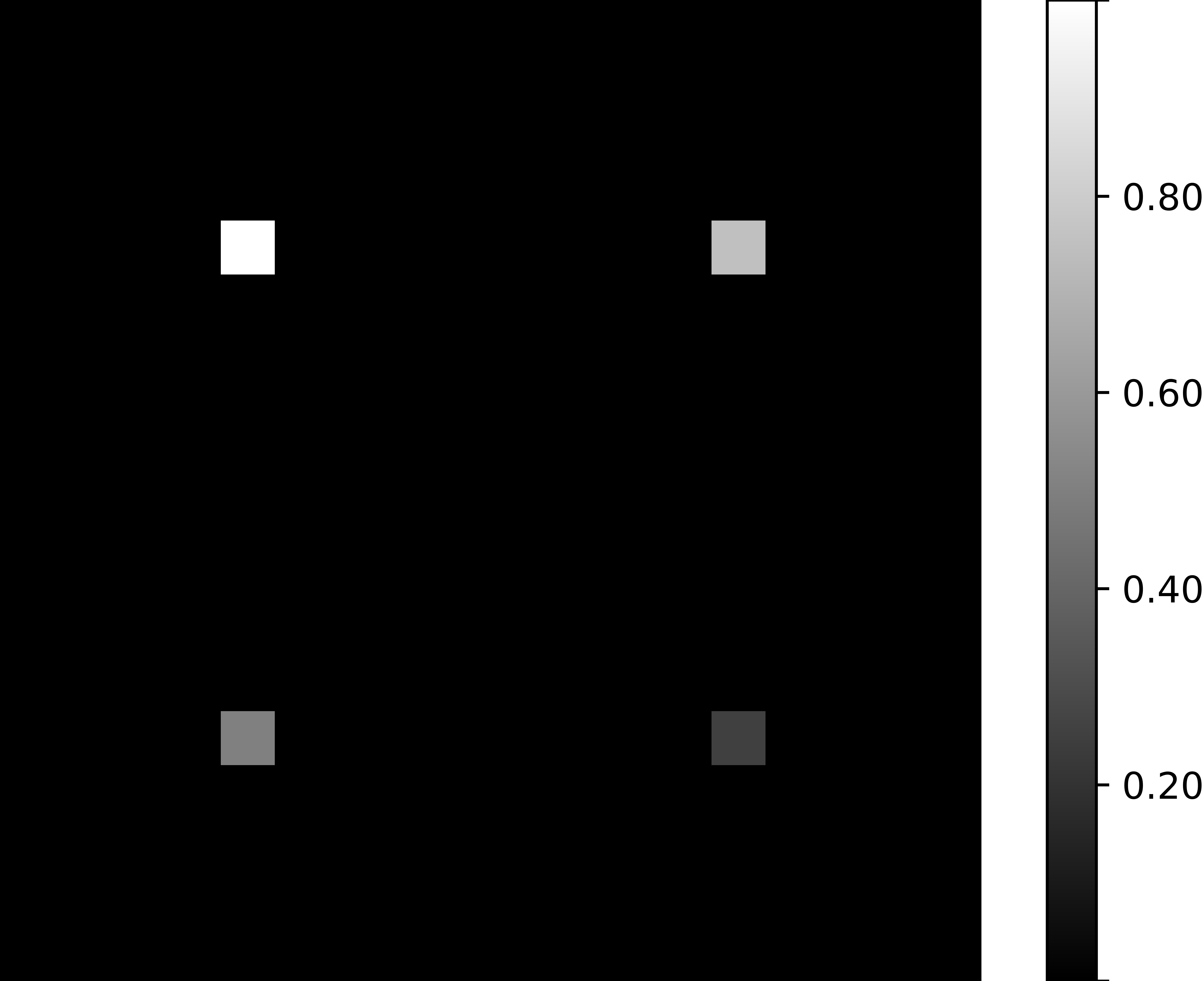}
		\caption{\centering\scriptsize Concentration phantom $\rho_{\text{GT}}$.}
		\label{subfig:concentration:gt}
	\end{subfigure}
	\hfil
	\begin{subfigure}[t]{\imratio\linewidth}
		\includegraphics[width=\linewidth]{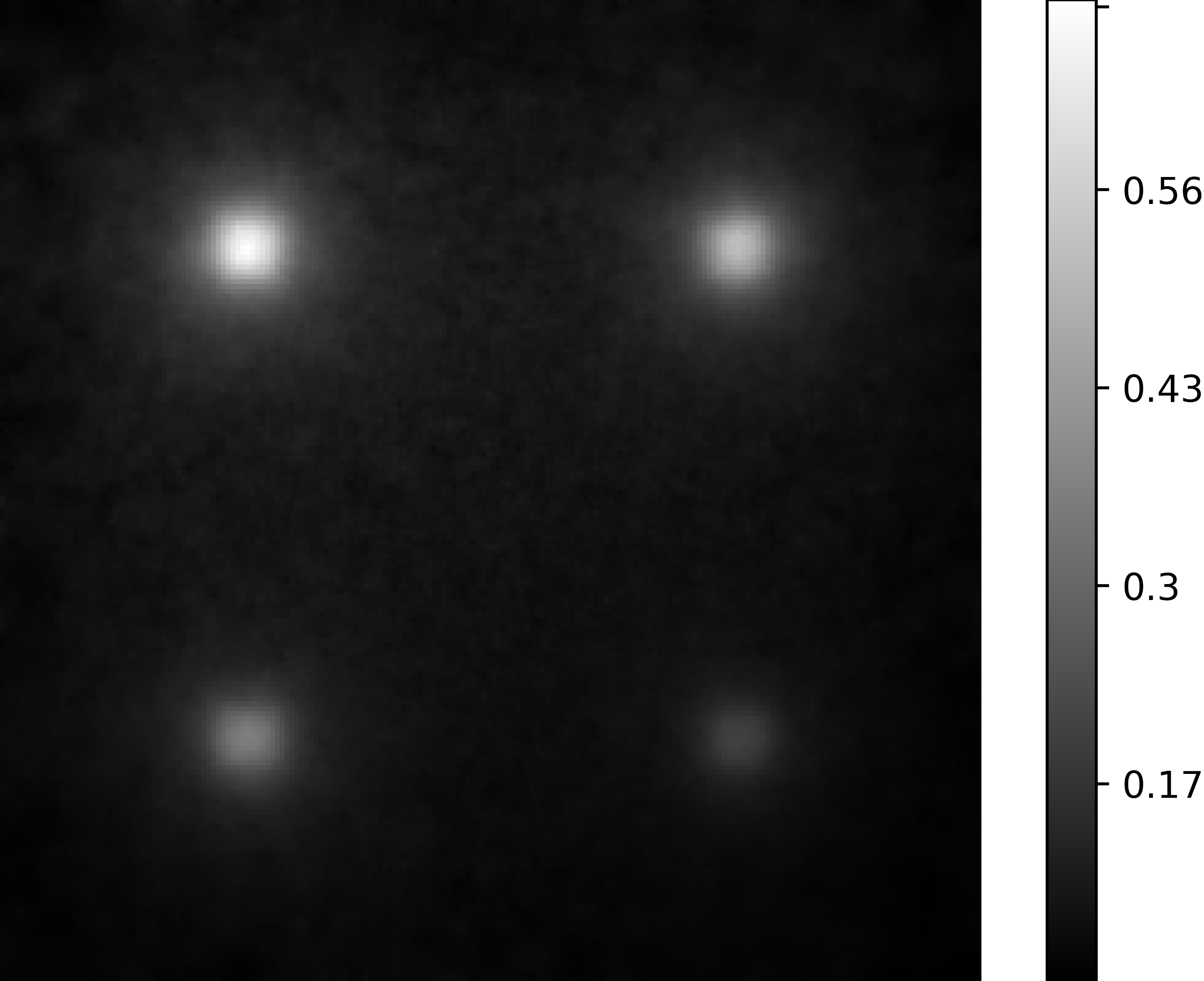}
		\caption{\centering\scriptsize $u$ PSNR 39.41, SSIM 0.9588.}
		\label{subfig:concentration:u}
	\end{subfigure}
	\hfil
	\begin{subfigure}[t]{\imratio\linewidth}
		\includegraphics[width={\linewidth}]{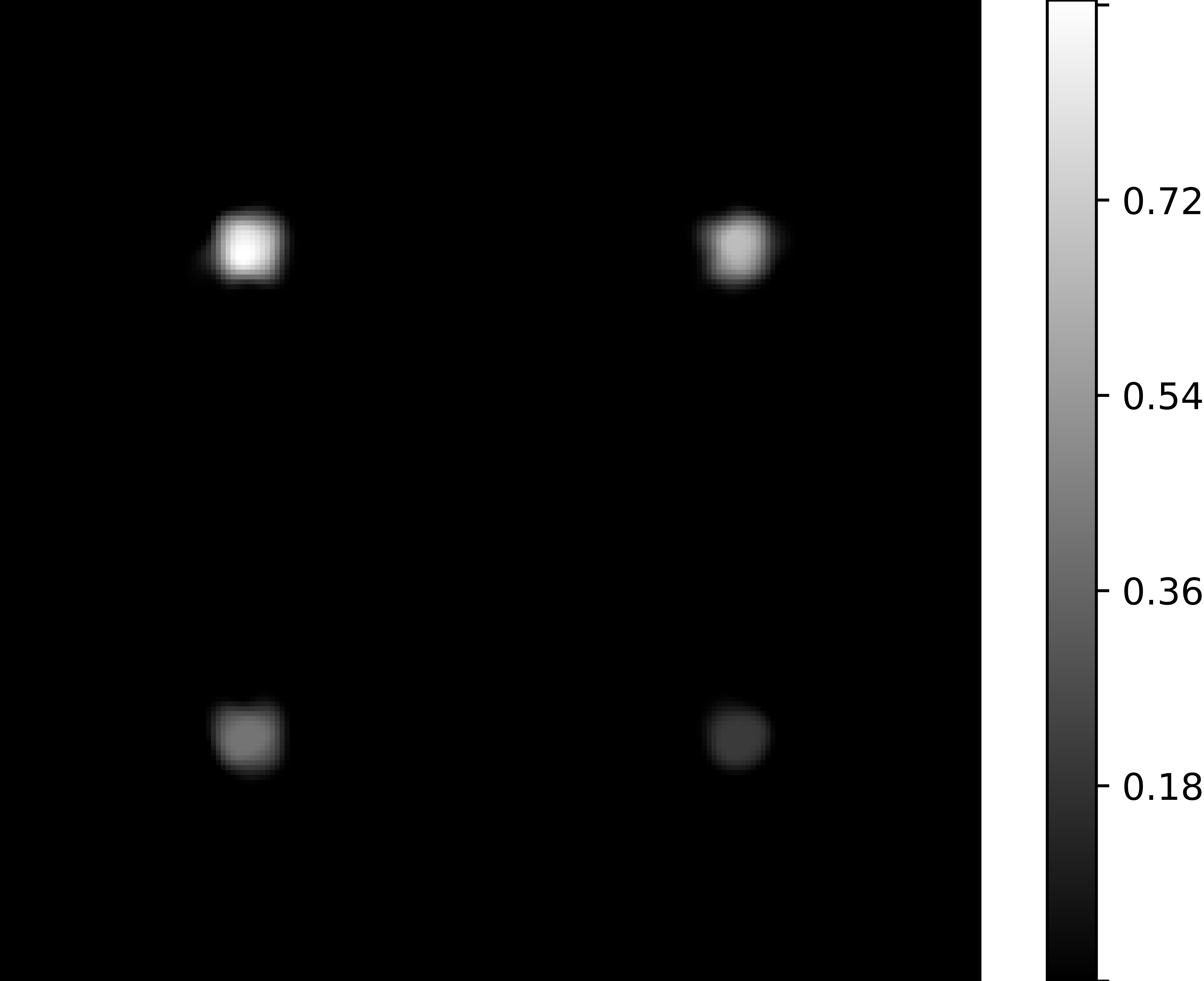}
		\caption{\centering\scriptsize $\rho$ PSNR 29.75, SSIM 0.9743.}
		\label{subfig:concentration:rho}
	\end{subfigure}
	
	\caption{Reconstructions of Open MPI Dataset-inspired shape and concentration phantoms of Experiment 4 with $10\times 10$ standard multi-patch scans (Fig.~\ref{subfig:scans:1010}) on a $200\times 200$ grid. The shapes of phantom~\ref{subfig:shape:gt} are well grasped by the reconstruction algorithm, although the edges are rounded off.}
	\label{fig:conc:shape:recs}
\end{figure}

To further test our algorithm, we performed the $10 \times 10$ scan on the region $\Omega = [-2,2]^2 $ with FoVs of amplitude $A_x = A_y =1$ of the phantoms inspire by the shape phantom and concentration phantom in the Open MPI Project~\cite{knopp2020openmpidata}. The ground truths can be seen in Fig.~\ref{subfig:shape:gt} and~\ref{subfig:concentration:gt}. The shape phantom in Fig.~\ref{subfig:shape:gt} is designed to test how much the reconstruction algorithm can preserve edges and pointy features and the level of smoothing out. The concentration phantom in Fig.~\ref{subfig:concentration:gt} is useful to test the ability of the reconstruction method to deal with four different concentration levels: 1, 0.75, 0.5 and 0.25 (see color bar in Fig.~\ref{subfig:concentration:gt}). In particular, it tests the robustness of the method with respect to low concentrations. Indeed, low concentrations are at risk of compression by the sparsity enforcing term $G_1$ in Eq.~\eqref{eq:second:pos:sparsity}. Moreover, the concentration phantom tests also the capability of the reconstruction method to keep the concentration ratios between different regions of the phantom. The reconstruction in Fig.~\ref{subfig:concentration:rho} where obtained lowering the sparsity enforcing parameter to $\beta = 0.1$, as our standard choice $\beta = 1$ made the lowest concentration (the bottom right dot in Fig.~\ref{subfig:concentration:gt}) almost invisible. In particular we can see that the method is able to reconstruct different levels of concentration (Fig.~\ref{subfig:concentration:rho}). The shape phantom in Fig.~\ref{subfig:shape:gt} is rather well reconstructed (see Fig.~\ref{subfig:shape:rho}), but the point is a little blunted. This effect can be seen in the slight rounding off of the squares in Fig.~\ref{subfig:concentration:gt}, which are more rounded in the reconstruction Fig.~\ref{subfig:concentration:rho}.

\subsection{Experiment 5: Stability with Respect to Perturbation}

\def\imratio{0.32}
\begin{figure}[t]
	\centering
	\begin{subfigure}[t]{\imratio\linewidth}
		\includegraphics[width=\linewidth]{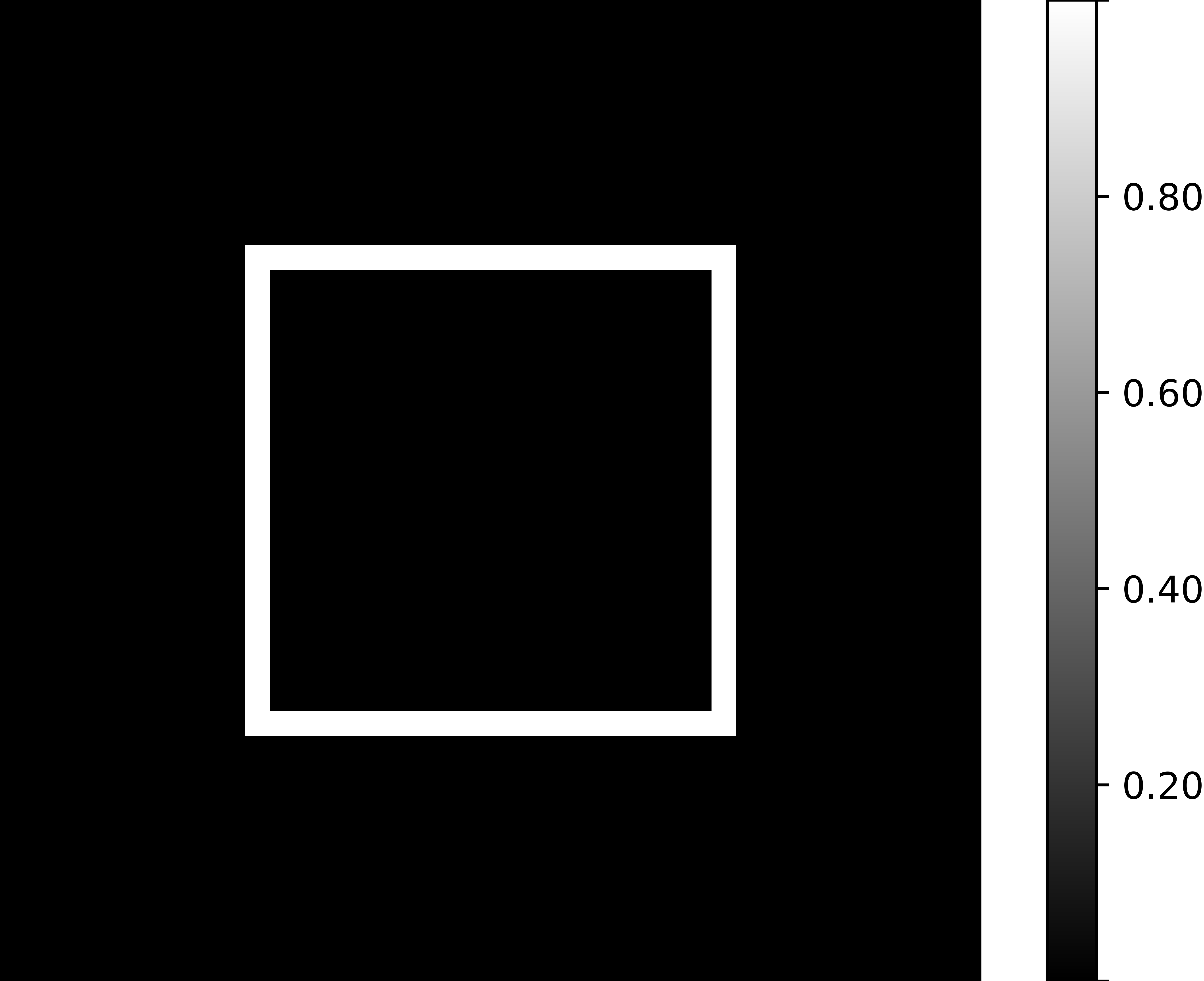}
		\caption{\centering\scriptsize Frame phantom $\rho_{\text{GT}}$.}
		\label{subfig:frame:gt}
	\end{subfigure}
	\hfil
	\begin{subfigure}[t]{\imratio\linewidth}
		\includegraphics[width=\linewidth]{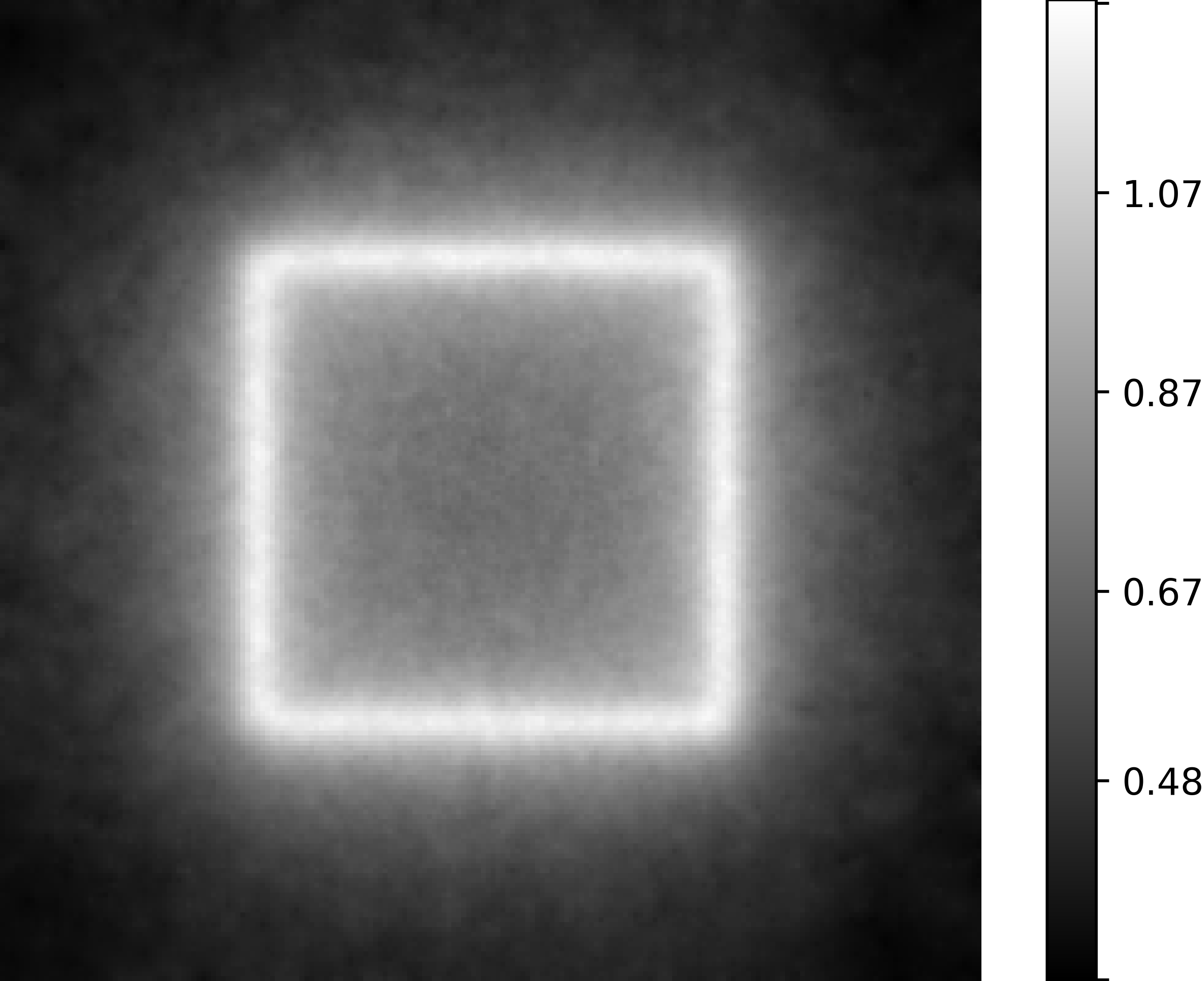}
		\caption{\centering\scriptsize  $u$ from unperturbed data, PSNR 33.88, SSIM 0.9262.}
		\label{subfig:frame:u}
	\end{subfigure}
	\hfil
	\begin{subfigure}[t]{\imratio\linewidth}
		\includegraphics[width=\linewidth]{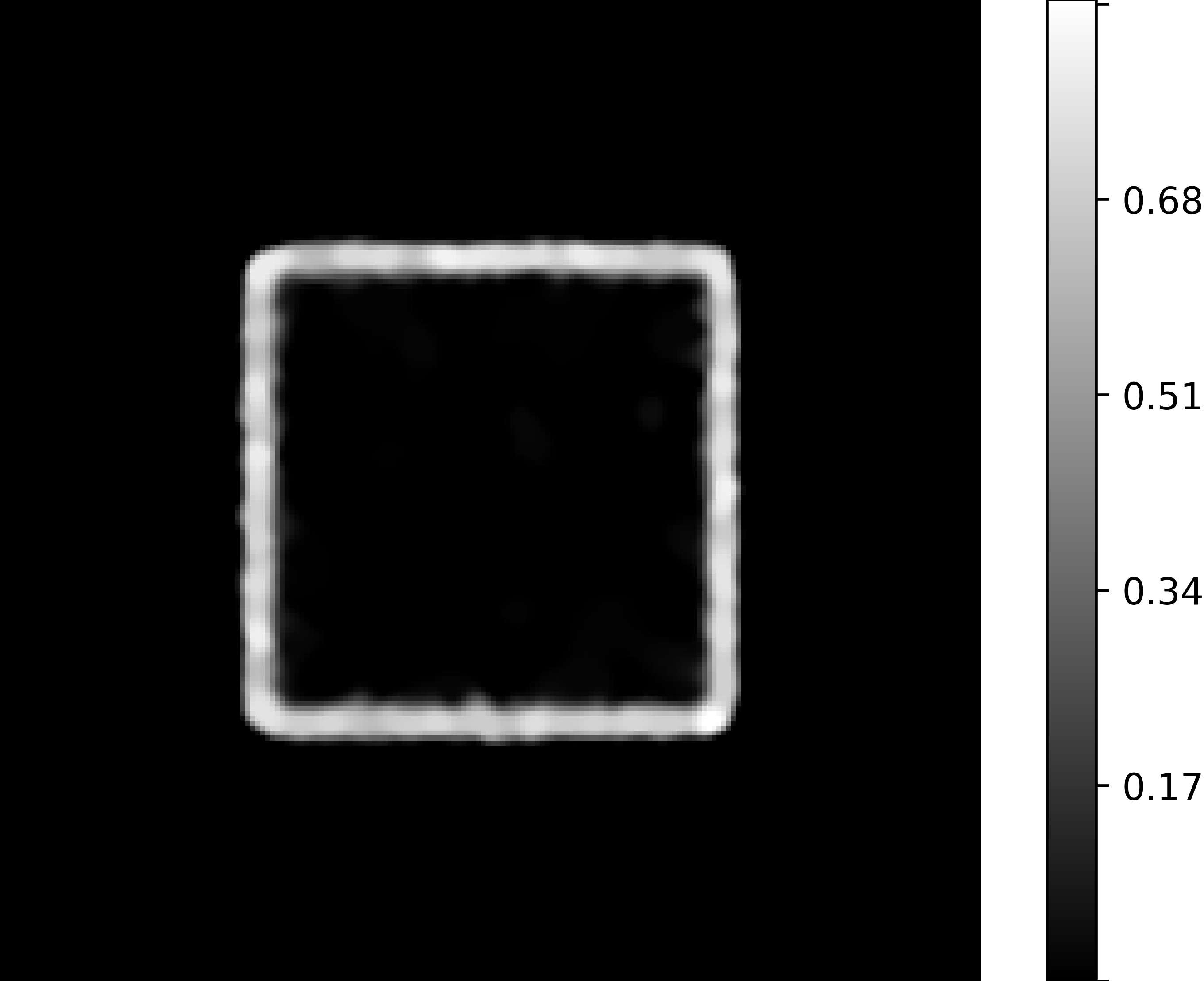}
		\caption{\centering\scriptsize $\rho$ from unperturbed data, PSNR 20.29, SSIM 0.9058.}
		\label{subfig:frame:rho}
	\end{subfigure}
	\par\medskip
	\begin{subfigure}[t]{\imratio\linewidth}
		\includegraphics[width=\linewidth]{Images/blank_dist.png}
	\end{subfigure}
	\hfil
	\begin{subfigure}[t]{\imratio\linewidth}
		\includegraphics[width=\linewidth]{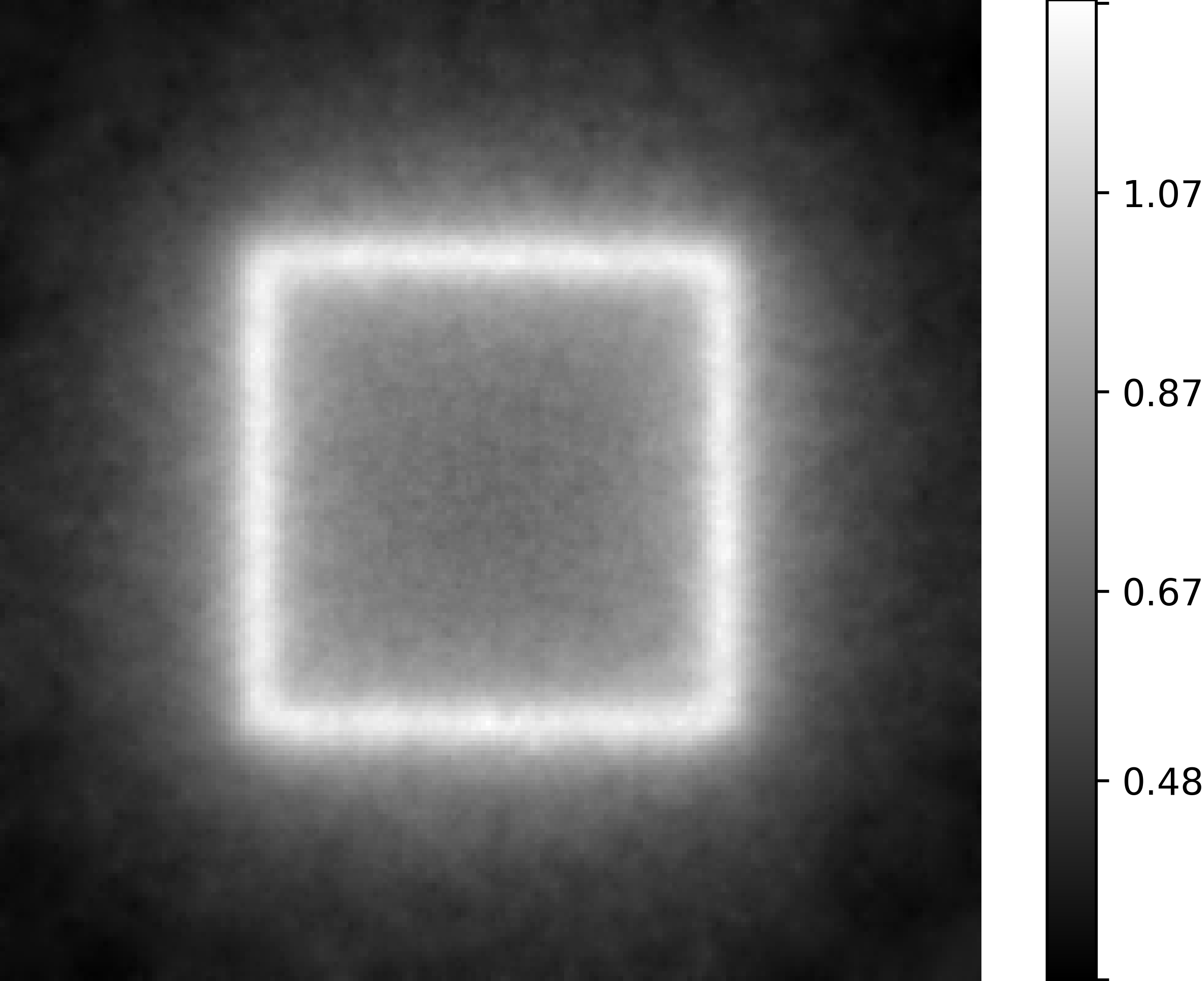}
		\caption{\centering\scriptsize $u$ from perturbed data, PSNR 33.10, SSIM 0.9181.}
		\label{subfig:frame:u:p}
	\end{subfigure}
	\hfil
	\begin{subfigure}[t]{\imratio\linewidth}
		\includegraphics[width={\linewidth}]{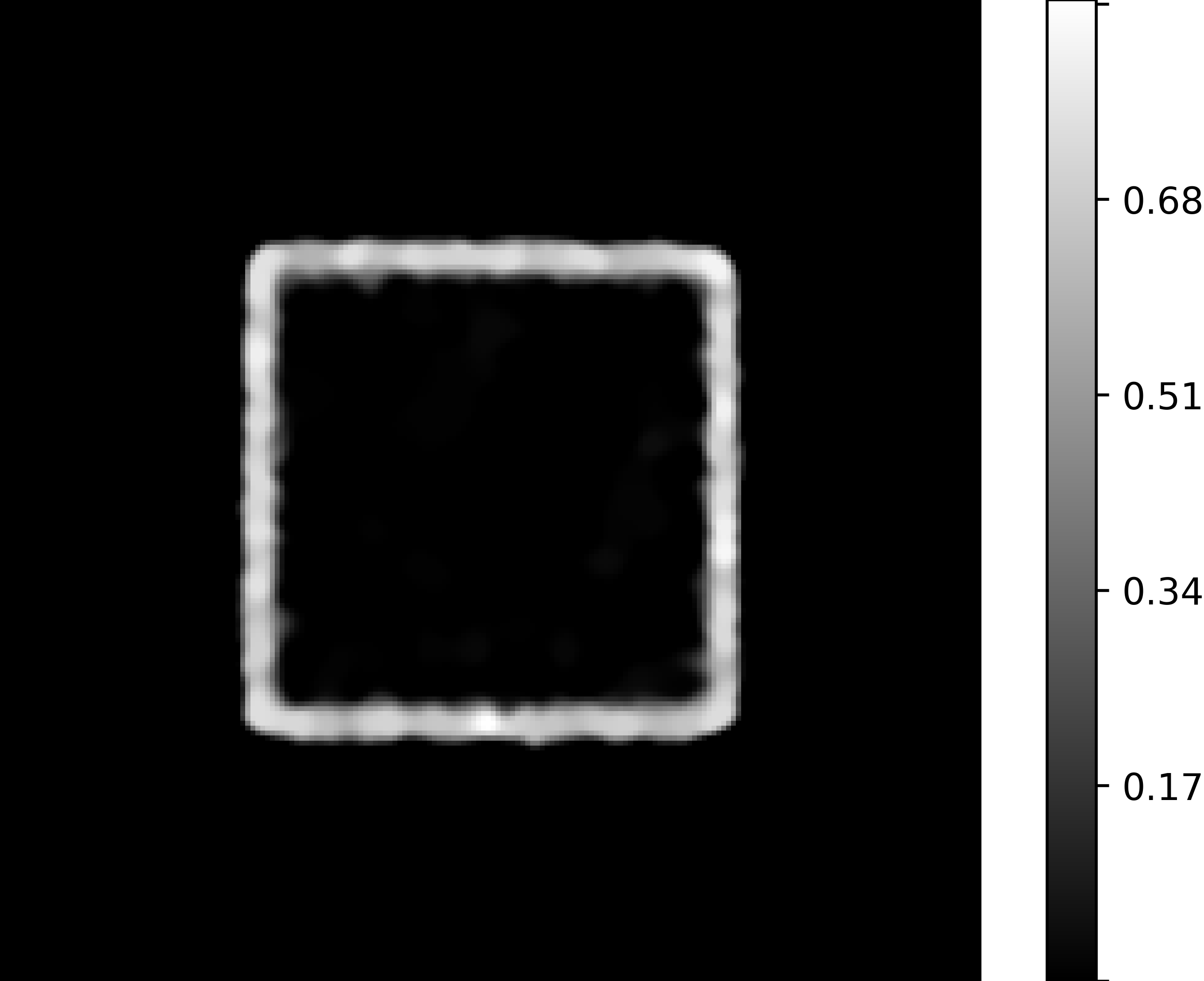}
		\caption{\centering\scriptsize $\rho$ from perturbed data, PSNR 19.97, SSIM 0.9002.}
		\label{subfig:frame:rho:p}
	\end{subfigure}
	\par\medskip
	\begin{subfigure}[t]{\imratio\linewidth}
		\includegraphics[width=\linewidth]{Images/blank_dist.png}
	\end{subfigure}
	\hfil
	\begin{subfigure}[t]{\imratio\linewidth}
		\includegraphics[width=\linewidth]{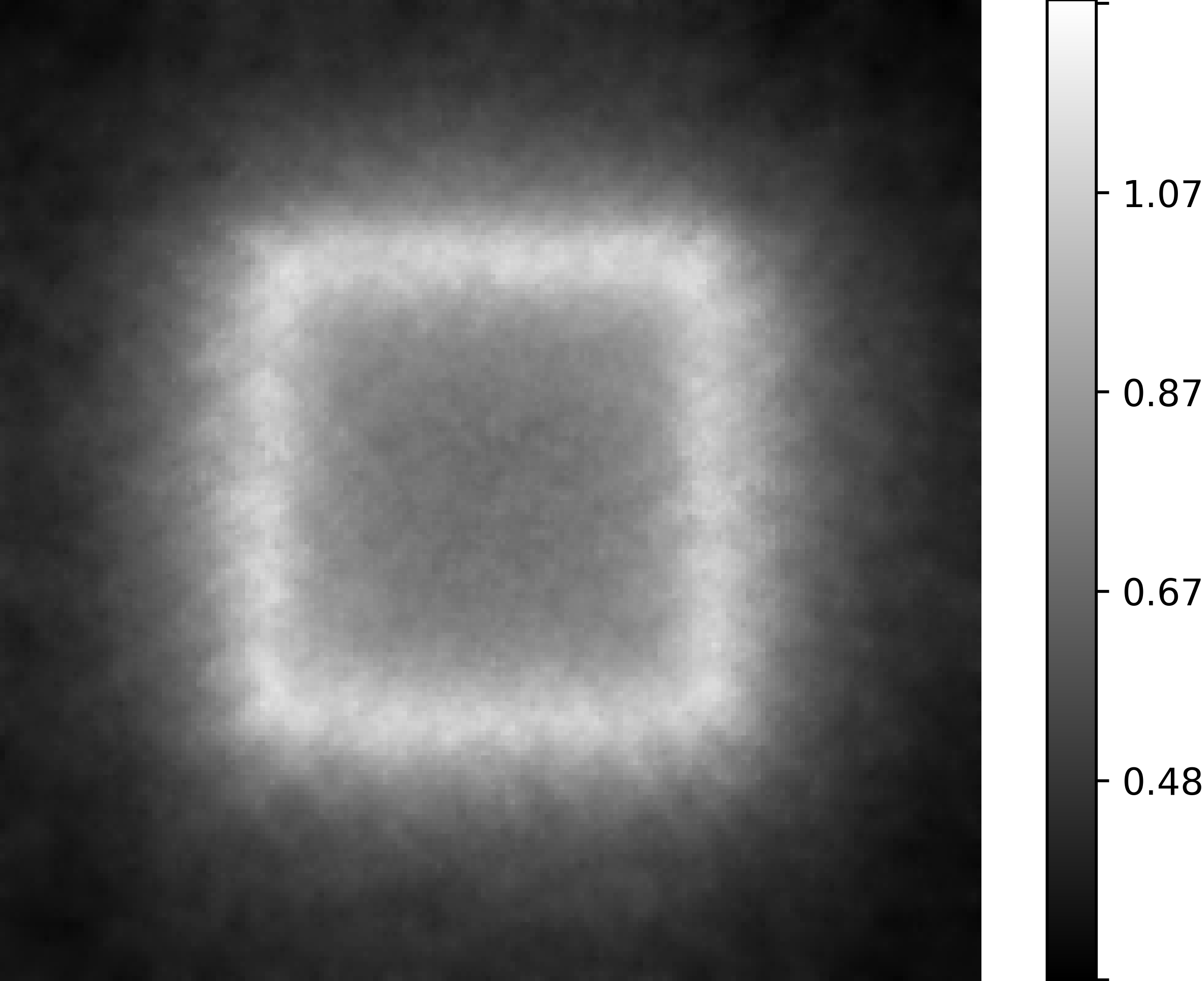}
		\caption{\centering\scriptsize $u$ from perturbed data, PSNR 26.17, SSIM 0.8534.}
		\label{subfig:frame:u:p2}
	\end{subfigure}
	\hfil
	\begin{subfigure}[t]{\imratio\linewidth}
		\includegraphics[width={\linewidth}]{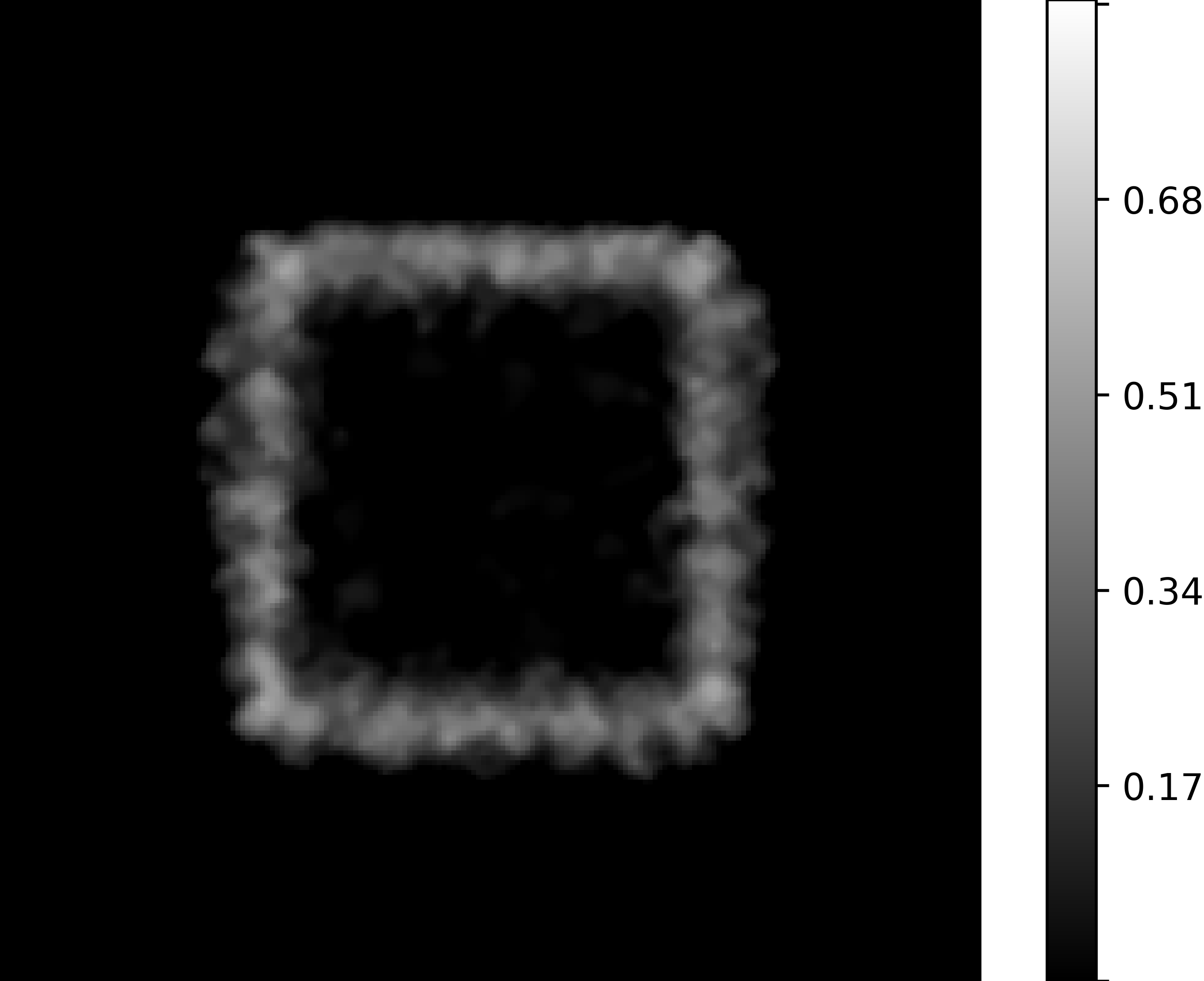}
		\caption{\centering\scriptsize $\rho$ from perturbed data, PSNR 15.63, SSIM 0.7351.}
		\label{subfig:frame:rho:p2}
	\end{subfigure}
	\caption{Reconstructions of the frame shaped phantom~\ref{subfig:frame:gt} from perturbed data (Fig.~\ref{subfig:frame:u:p},~\ref{subfig:frame:rho:p},~\ref{subfig:frame:u:p2},~\ref{subfig:frame:rho:p2}) and from unperturbed data, using $10\times 10$ standard multi-patch (Fig.~\ref{subfig:scans:1010}) scans on a $200\times 200$ grid (details of the perturbations can be found in the description to Experiment 5). The slightly perturbed reconstruction~\ref{subfig:frame:rho:p} of the frame phantom~\ref{subfig:frame:gt} has more frayed edges. An increased perturbation in the data leads to a less exact reconstruction, cf.~\ref{subfig:frame:rho:p2}; however, basic structures can still be observed.}
	\label{fig:pert:recs}
\end{figure}

\def\imratio{0.32}
\begin{figure}[t]
	\centering
	\begin{subfigure}[t]{\imratio\linewidth}
		\includegraphics[width=\linewidth]{Images/vessel_gt.png}
		\caption{\centering\scriptsize Vessel $\rho$.}
		\label{subfig:vessel:p:gt}
	\end{subfigure}
	\hfil
	\begin{subfigure}[t]{\imratio\linewidth}
		\includegraphics[width=\linewidth]{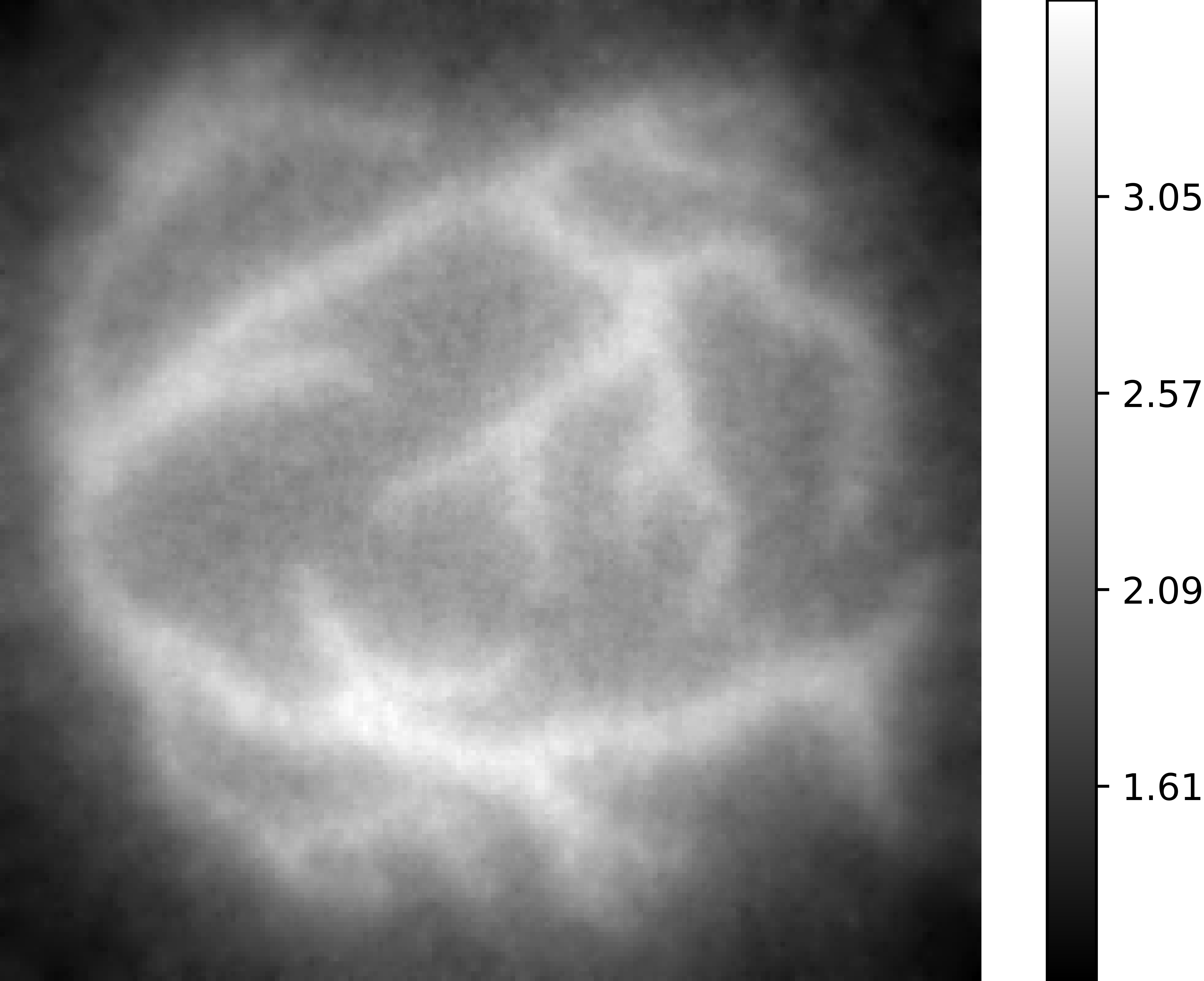}
		\caption{\centering\scriptsize $u$ from perturbed data, PSNR 32.35, SSIM 0.8864.}
		\label{subfig:vessel:u:p}
	\end{subfigure}
	\hfil
	\begin{subfigure}[t]{\imratio\linewidth}
		\includegraphics[width={\linewidth}]{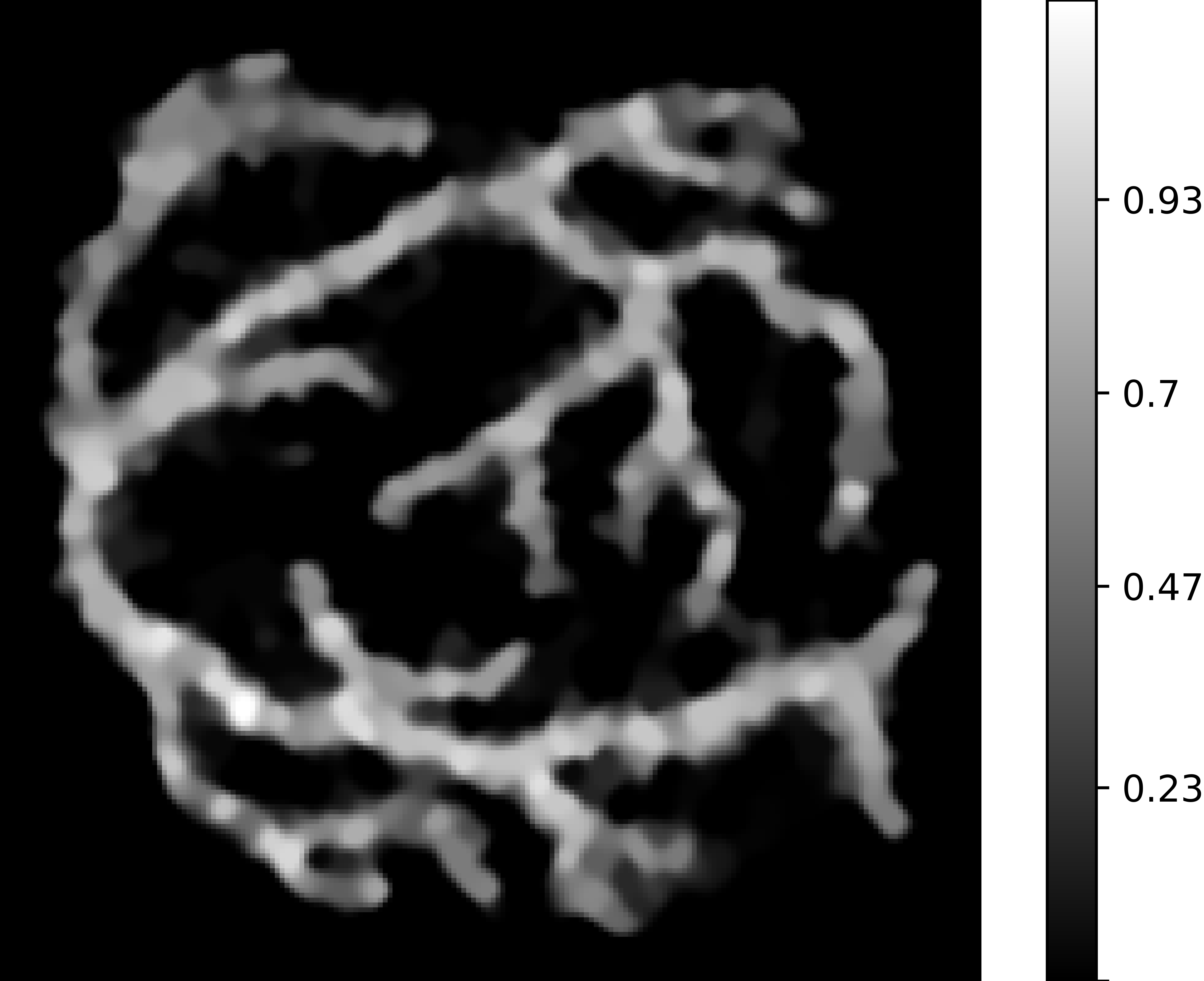}
		\caption{\centering\scriptsize $\rho$ from perturbed data, PSNR 13.23, SSIM 0.5963.}
		\label{subfig:vessel:rho:p}
	\end{subfigure}
	\caption{Reconstruction on a $200\times 200$ grid of the vessel shaped phantom with perturbations on $10\times 10$ scanning data in Experiment 5. Compared with the reconstruction of the same phantom from unperturbed data in Fig~\ref{subfig:exp:vess:multi:10:rho}, slight perturbations do affect the final result of the reconstruction, but the main features of the phantom are still recognizable, showing the stability of the reconstruction algorithm to small perturbations in the scanning trajectories even for complicated phantoms.}
	\label{fig:pert:vessel:recs}
\end{figure}

In Experiment 5 we study the stability of the reconstructions with respect to perturbations of the scanning trajectories. We considered again the vessel phantom (Fig.~\ref{subfig:vessel:p:gt}) and a frame-shaped phantom (Fig.~\ref{subfig:frame:gt}), whose straight lines help see the aberrations due to the perturbations. We considered the region $\Omega = [-2,2]^2$ with FoV with amplitudes $A_x = A_y = 1$ and a standard $10\times 10$ multi-patch scan (as in Fig.~\ref{subfig:scans:1010}) for both phantoms. Therefore we have $\Xi = 100$ scans on $\Omega$ with offset vectors $b_{\xi}$ as in Eq.~\eqref{eq:offset:std:multi} and angles $\alpha_{\xi}=0$ for every $\xi \in\lbrace 0,\dots ,\Xi -1\rbrace$, to which we added random perturbations, i.e., we used as offset vectors and angles the following quantities:
\begin{equation}
	\hat{b}_{\xi} = b_{\xi} + X_{\xi}e_1 + Y_{\xi}e_2 , \quad \hat{\alpha}_{\xi} = 0 + Z_{\xi} 
\end{equation}
where $X_{\xi}\sim \text{U}([-\frac{A_x}{100},\frac{A_x}{100}])$, $Y_{\xi}\sim \text{U}([-\frac{A_y}{100},\frac{A_y}{100}])$ and $Z_{\xi}\sim\text{U}([-\frac{2\pi}{360},\frac{2\pi}{360}])$ for the experiment in Fig.\ref{subfig:frame:u:p}, \ref{subfig:frame:rho:p}, \ref{subfig:vessel:u:p}, \ref{subfig:vessel:rho:p} and $X_{\xi}\sim \text{U}([-\frac{A_x}{10},\frac{A_x}{10}])$, $Y_{\xi}\sim \text{U}([-\frac{A_y}{10},\frac{A_y}{10}])$ and $Z_{\xi}\sim\text{U}([-2\frac{2\pi}{360},2\frac{2\pi}{360}])$ in experiment in Fig.~\ref{subfig:frame:u:p2},~\ref{subfig:frame:rho:p2}. Here $X_{\xi}$, $Y_\xi$ and $Z_\xi$ are real valued random variables and $e_1 ,e_2$ are canonical basis vectors in $\mathbb{R}^2$. This means that the offset vectors are perturbed by at most one 100-th or a 10-th of the amplitude in each direction and a random angle of at most one or two degrees. Comparing the unperturbed and the perturbed reconstructions of the frame-shape phantom in Fig.~\ref{subfig:frame:rho},~\ref{subfig:frame:rho:p},~\ref{subfig:frame:rho:p2} those of the vessel phantom in Fig.~\ref{subfig:exp:vess:multi:10:rho},~\ref{subfig:vessel:rho:p} and their PSNR values in Tab.~\ref{tab:exp}, it is clear that comparable results are obtainable even taking into account small enough perturbations. For increased perturbations (compare Fig~\ref{subfig:frame:rho:p},~\ref{subfig:frame:rho:p2}) the reconstruction show a higher level of aberrations and distortions, as expected.

\subsection{Experiment 6: Scanning While Moving}\label{sec:scan:move}

\def\imratio{0.32}
\begin{figure}[t]
	\centering
	\begin{subfigure}[t]{\imratio\linewidth}
		\includegraphics[width=\linewidth]{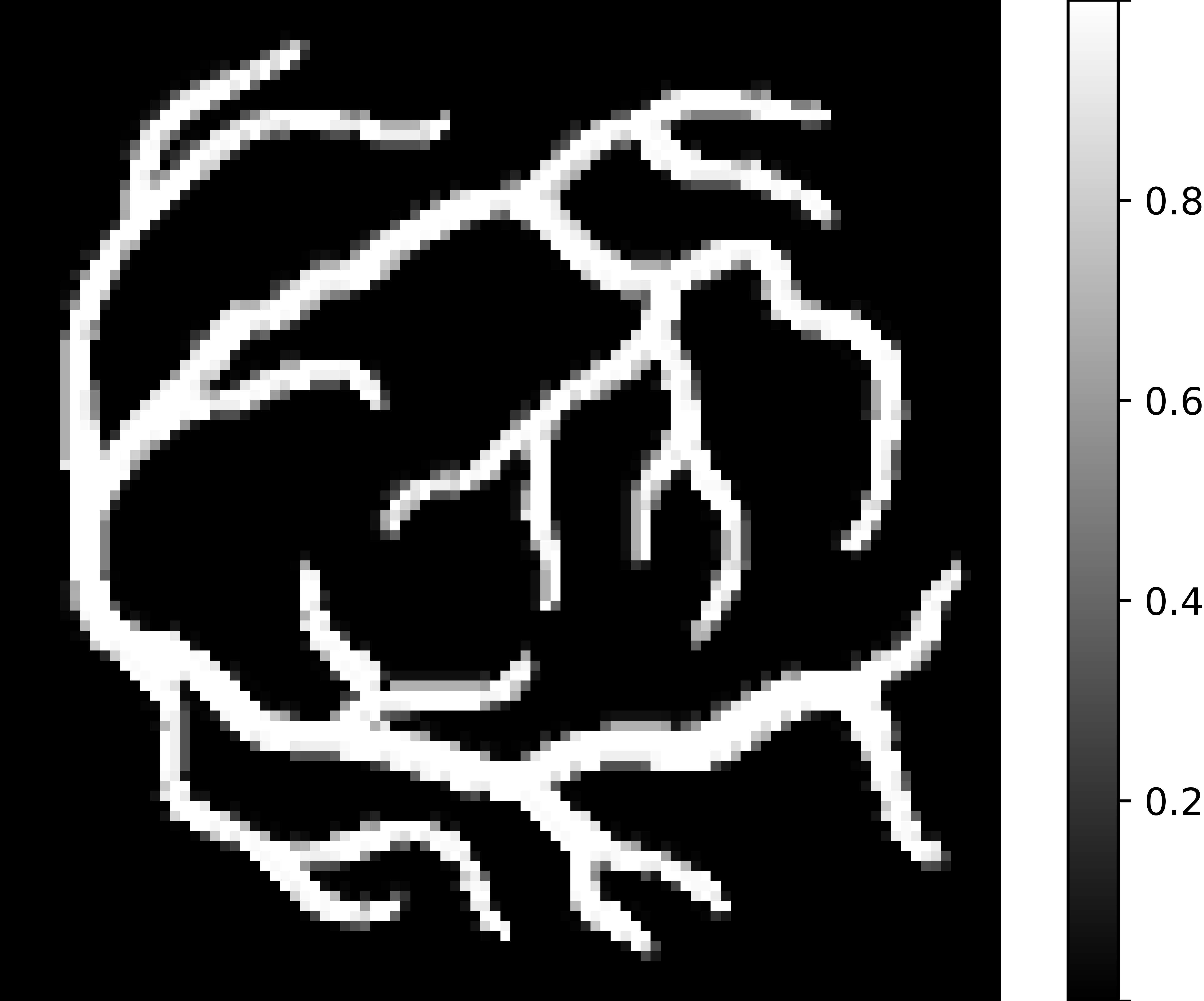}
		\caption{\centering\scriptsize Vessel $100\times 100$ $\rho_{\text{GT}}$.}
		\label{subfig:1000:gt}
	\end{subfigure}
	\hfil
	\begin{subfigure}[t]{\imratio\linewidth}
		\includegraphics[width=\linewidth]{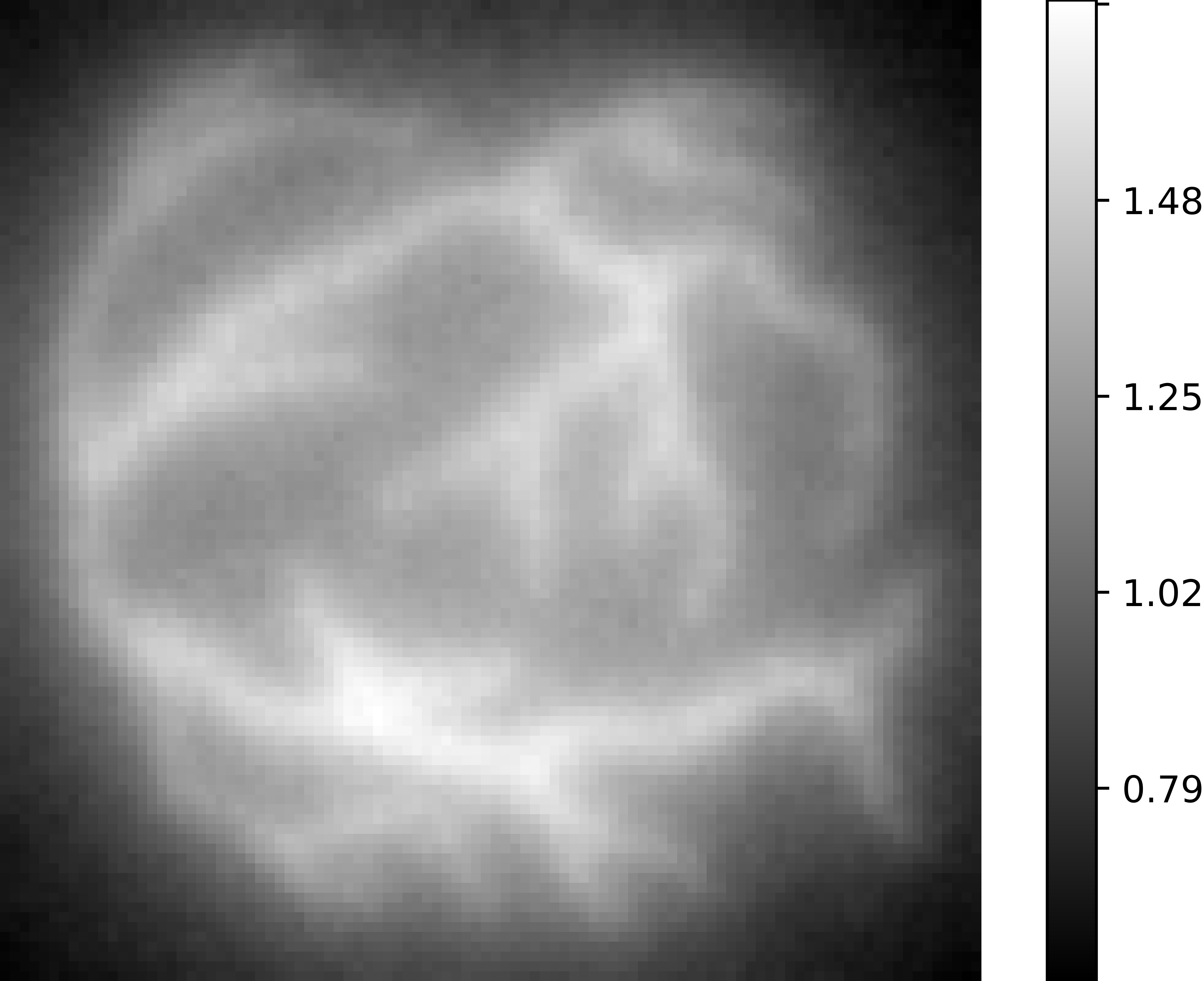}
		\caption{\centering\scriptsize $u$, PSNR 37.68, SSIM 0.9700.}
		\label{subfig:1000:u}
	\end{subfigure}
	\hfil
	\begin{subfigure}[t]{\imratio\linewidth}
		\includegraphics[width=\linewidth]{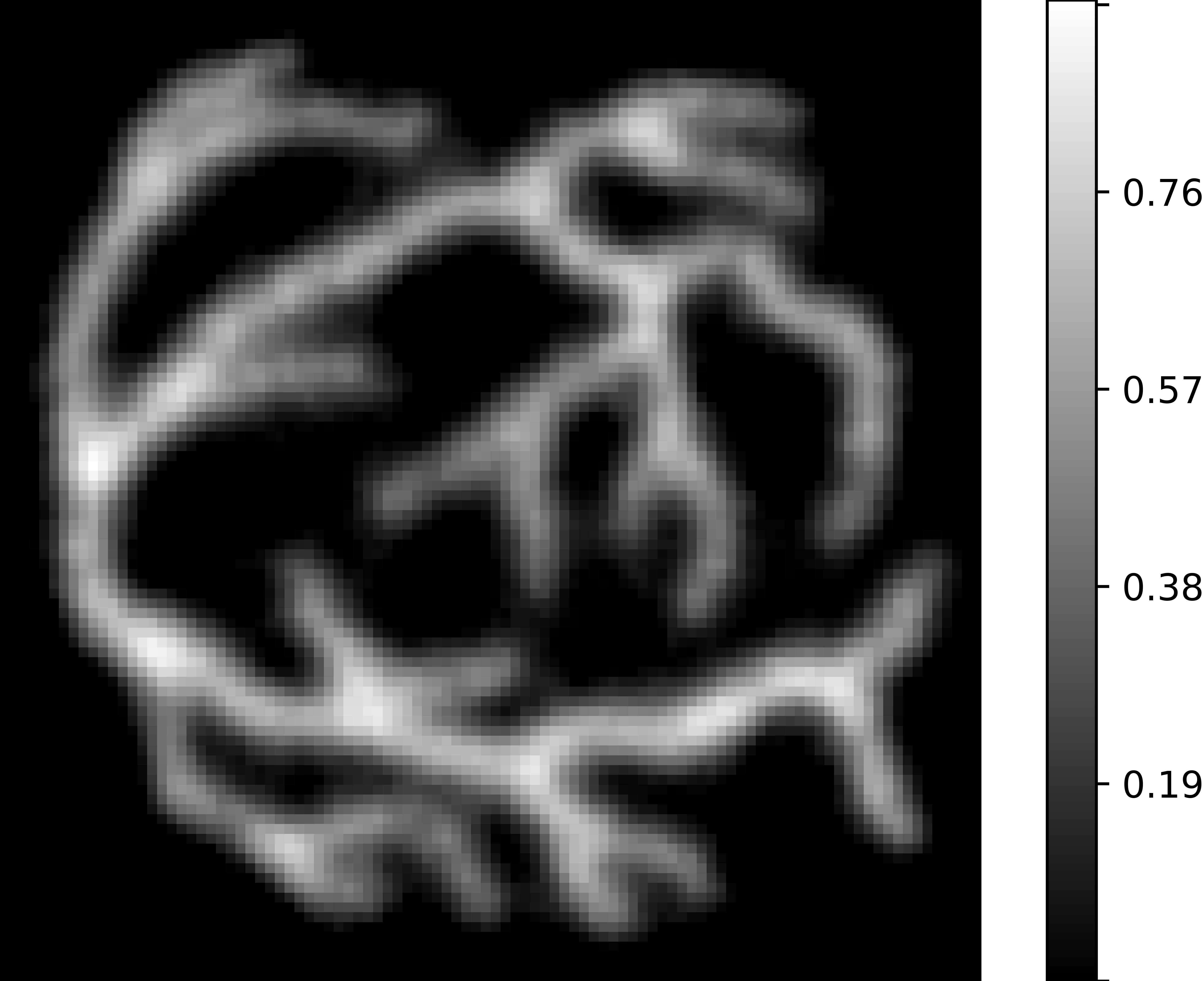}
		\caption{\centering\scriptsize $\rho$, PSNR 12.81, SSIM 0.5265.}
		\label{subfig:1000:rho}
	\end{subfigure}
	\par\medskip
	\begin{subfigure}[t]{\imratio\linewidth}
		\includegraphics[width=\linewidth]{Images/blank_dist.png}
	\end{subfigure}
	\hspace{-3mm}
	\begin{subfigure}[t]{\imratio\linewidth}
		\includegraphics[width=\linewidth]{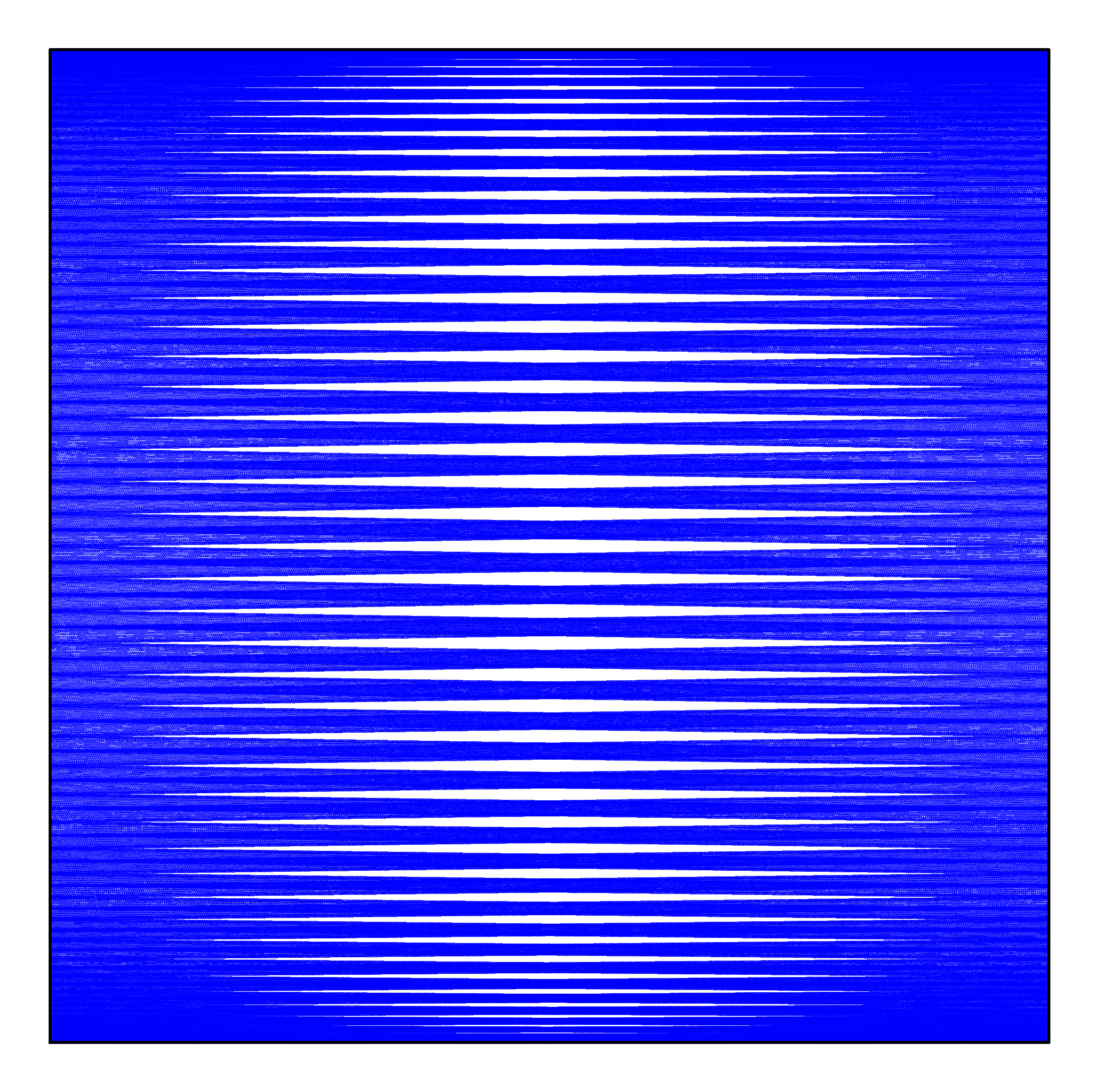}
		\caption{\centering\scriptsize 1000 Lissajous trajectories while moving.}
		\label{fig:1000:scans}
	\end{subfigure}
	\hfil
	\begin{subfigure}[t]{\imratio\linewidth}
		\includegraphics[width=\linewidth]{Images/blank_dist.png}
	\end{subfigure}
	
	\caption{Reconstruction of a $100\times 100$ vessel phantom~\ref{subfig:1000:gt} using 1000 scans performed while moving the FoV~\ref{fig:1000:scans} (Experiment 6). In this experiment the FoV and the region $\Omega$ have exactly the same size $[-1,1]^2$ and the ground truth~\ref{subfig:1000:gt} has varying levels of concentration, mimicking what would happen in a real phantom. We observe that with a high sampling, most of the features of the phantom are already clearly recognizable in the trace $u$ (\ref{subfig:1000:u}) before the deconvolution step and are rather well reconstructed after the second stage of the algorithm (\ref{subfig:1000:rho}). }
	\label{fig:1000:recs}
\end{figure}

In the standard multi-patching and in the experiments performed so far, each scan (1632 sampling) is performed with fixed offset and angle vector, then it is moved, then another scan is performed and so on. This means that between each scan there is a time window in which the offset vector is moved and in which no data are collected. We show in this experiment that it is possible to keep performing scans continuously while moving. This allows to cover big regions $\Omega$ while also optimizing the time for the whole scan. We performed the reconstruction in Fig.~\ref{fig:1000:recs} on a $100\times 100$ grid, where the $100\times 100$ pixel version of the vessel phantom in Fig.~\ref{subfig:1000:gt} is obtained by bicubic interpolation and projection onto $[0,1]$. This produced a more natural particle distribution, as it is not a characteristic function anymore, but it has values of concentration in the whole set $[0,1]$. Moreover, the region $\Omega = [a,b]\times [c,d] = [-1,1]^2$ in this experiment coincides in dimension with the FoV used of amplitudes $A_x = A_y = 1$. With this choice, the FoV is not a zoom-in of a portion of the phantom, as it was the case of the experiments in Fig.~\ref{fig:exp:vess:multi:recs}, resulting in a more challenging phantom to reconstruct, as the features of the vessel phantom have been shrunk without modifying the scanning setup of 1632 points and the Lissajous trajectories. We simulated 1000 scans while moving the FoV from left to right and both starting and finishing just outside the region $\Omega$, i.e., $b(t)$ is a uniform rectilinear motion from $b(0)=(a-A_x ,0)$ to $b(T)=(b+A_x ,0)$. We point out that, like in Experiment 3, with this choice of off-set vectors, the region $\Omega$ is over-scanned. This is done to prevent undersampling of any region in $\Omega$. More precisely, we produce data using equation Eq.~\eqref{eq:gen:liss:data} scanning along the generalized Lissajous trajectory $\Lambda (t)$ in Eq.~\eqref{eq:gen:liss} where the offset trajectory $b(t)$ and angle function $\alpha (t)$ are of the form:
\begin{equation}
	b(t) = \biggl ( 1-\frac{t}{T}\biggr ) (a-A_x)+\frac{t}{T}(b+A_x) ,\quad \alpha (t) = 0\quad \text{for all }t\in [0,T]
\end{equation}
where $T>0$ is the time needed to perform the 1000 scans. This results in 816,007 total sample points in $\Omega$, which correspond ca. 500 scans (see Fig.~\ref{fig:1000:scans}). Because the grid is smaller than in the other experiments we set the sparsity enforcing parameter to $\beta = 0.1$, and because the decreased grid size mean faster computation we let Alg.~\ref{alg:lasso} run for 100,000 iterations. The reconstruction obtained can be seen in Fig.~\ref{fig:1000:recs}.

\subsection{Experiment 7: Comparison with a Baseline Approach I}\label{sec:exp:sm:22}

\def\imratio{0.24}
\begin{figure}[t]
	\centering
	\begin{subfigure}[t]{\imratio\linewidth}
		\includegraphics[width=\linewidth]{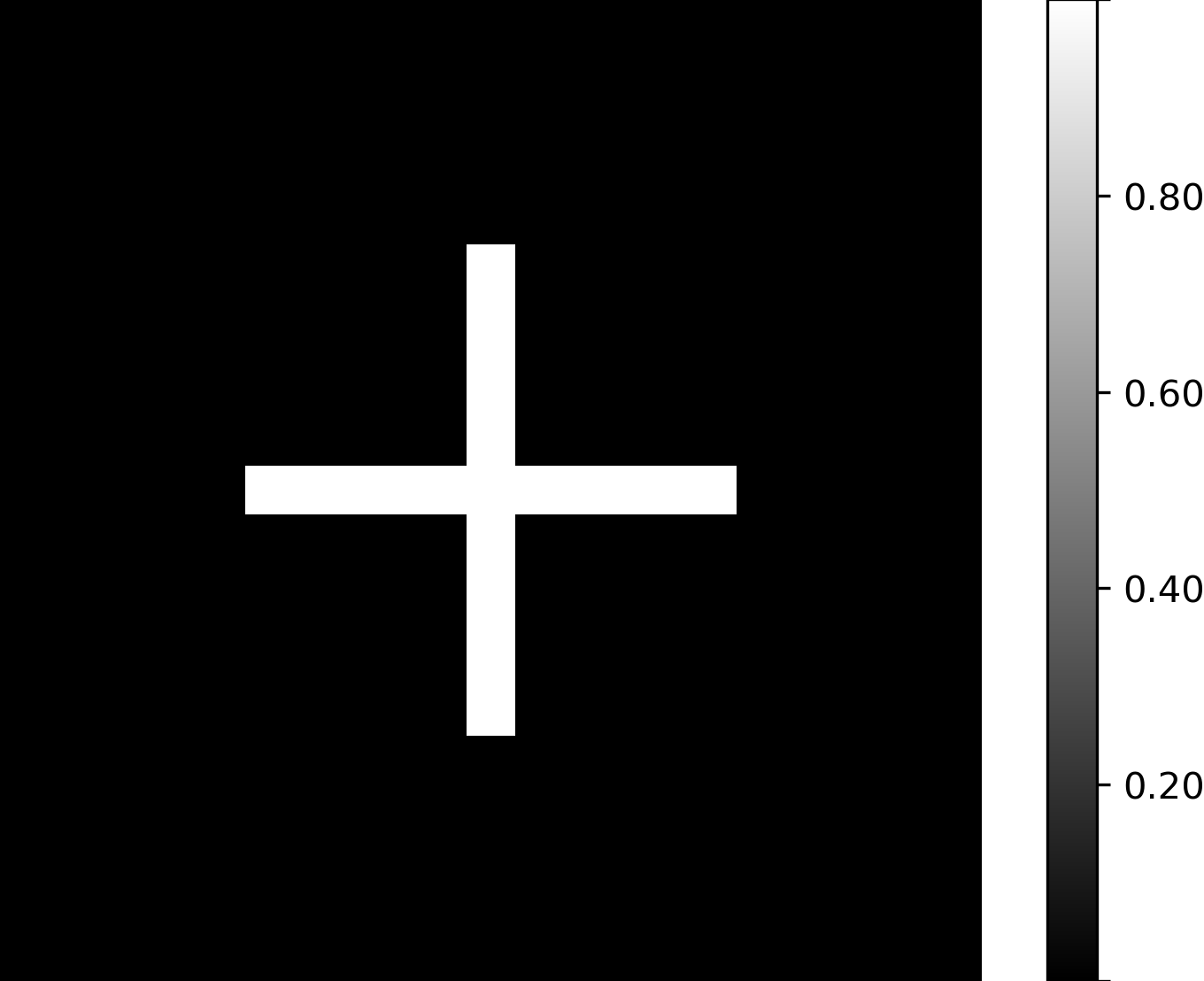}
		\caption{\centering\scriptsize Plus-shaped phantom $40\times 40$ $\rho_{\text{GT}}$.}
		\label{subfig:plus22:gt}
	\end{subfigure}
	\hfil
	\begin{subfigure}[t]{\imratio\linewidth}
		\includegraphics[width=\linewidth]{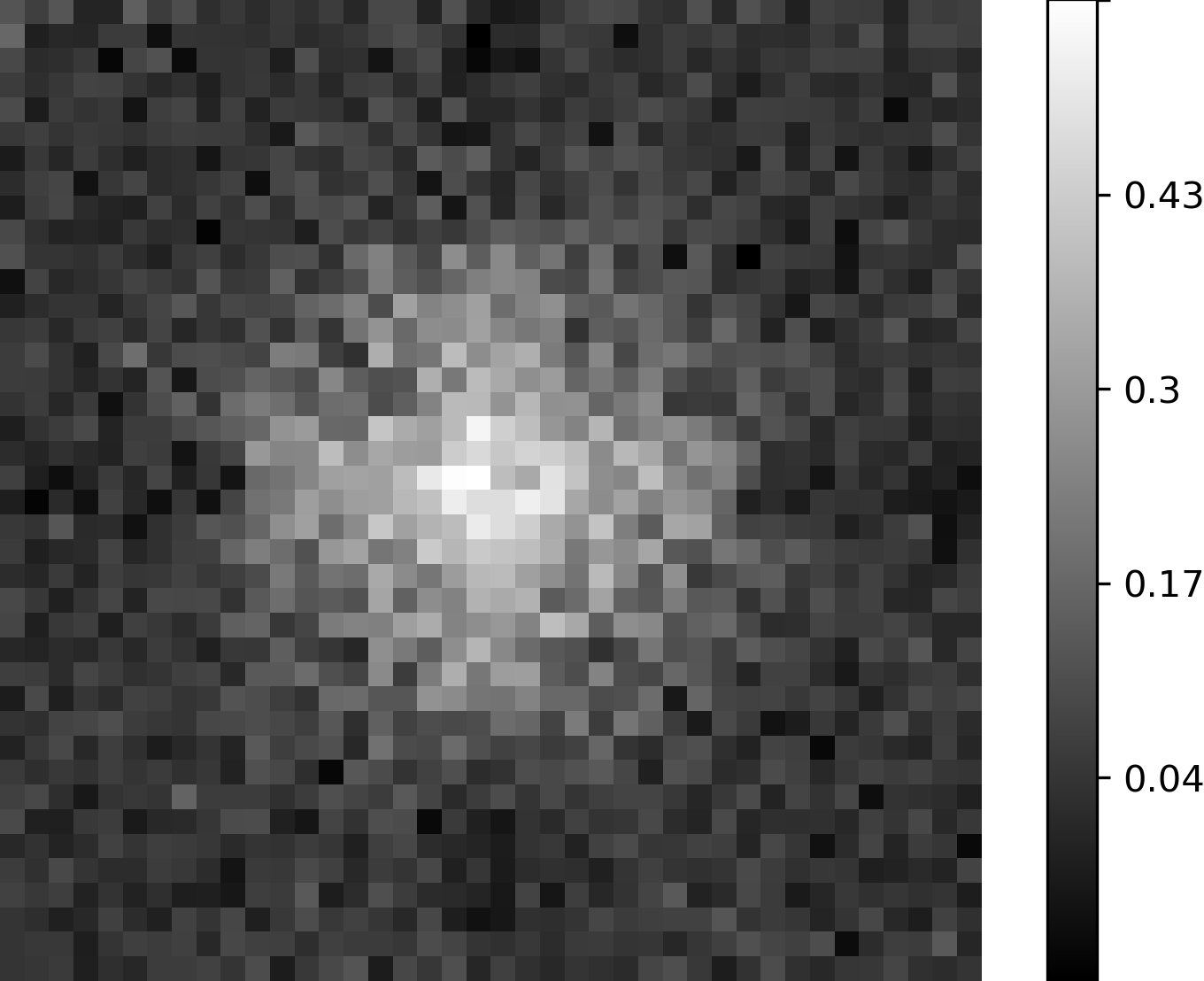}
		\caption{\centering\scriptsize Patched SM reco, PSNR 14.64, SSIM 0.0572.}
		\label{subfig:plus22:SM}
	\end{subfigure}
	\hfil
	\begin{subfigure}[t]{\imratio\linewidth}
		\includegraphics[width=\linewidth]{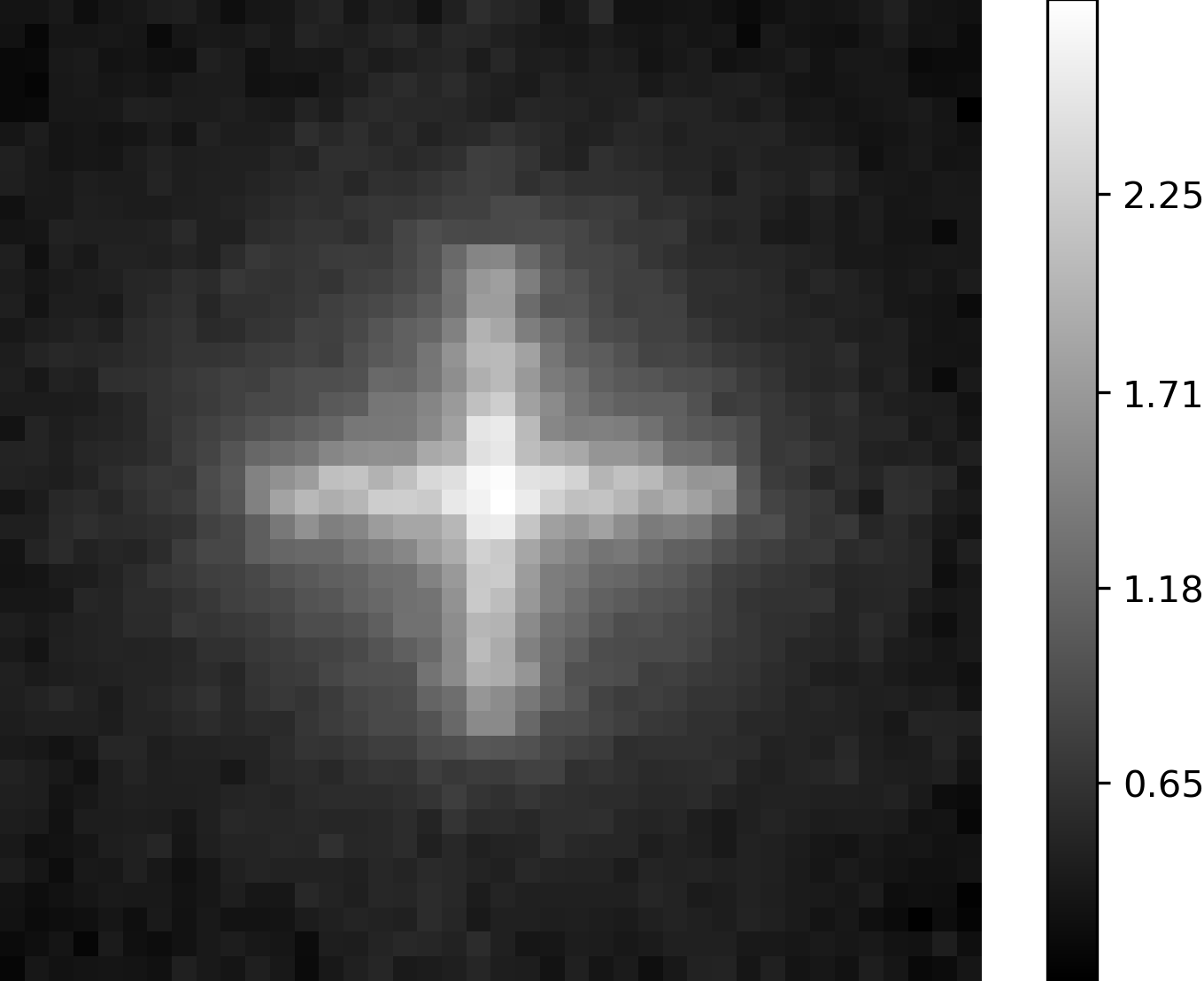}
		\caption{\centering\scriptsize First stage $u$, PSNR 30.01, SSIM 0.9092.}
		\label{subfig:plus22:u}
	\end{subfigure}
	\hfil
	\begin{subfigure}[t]{\imratio\linewidth}
		\includegraphics[width=\linewidth]{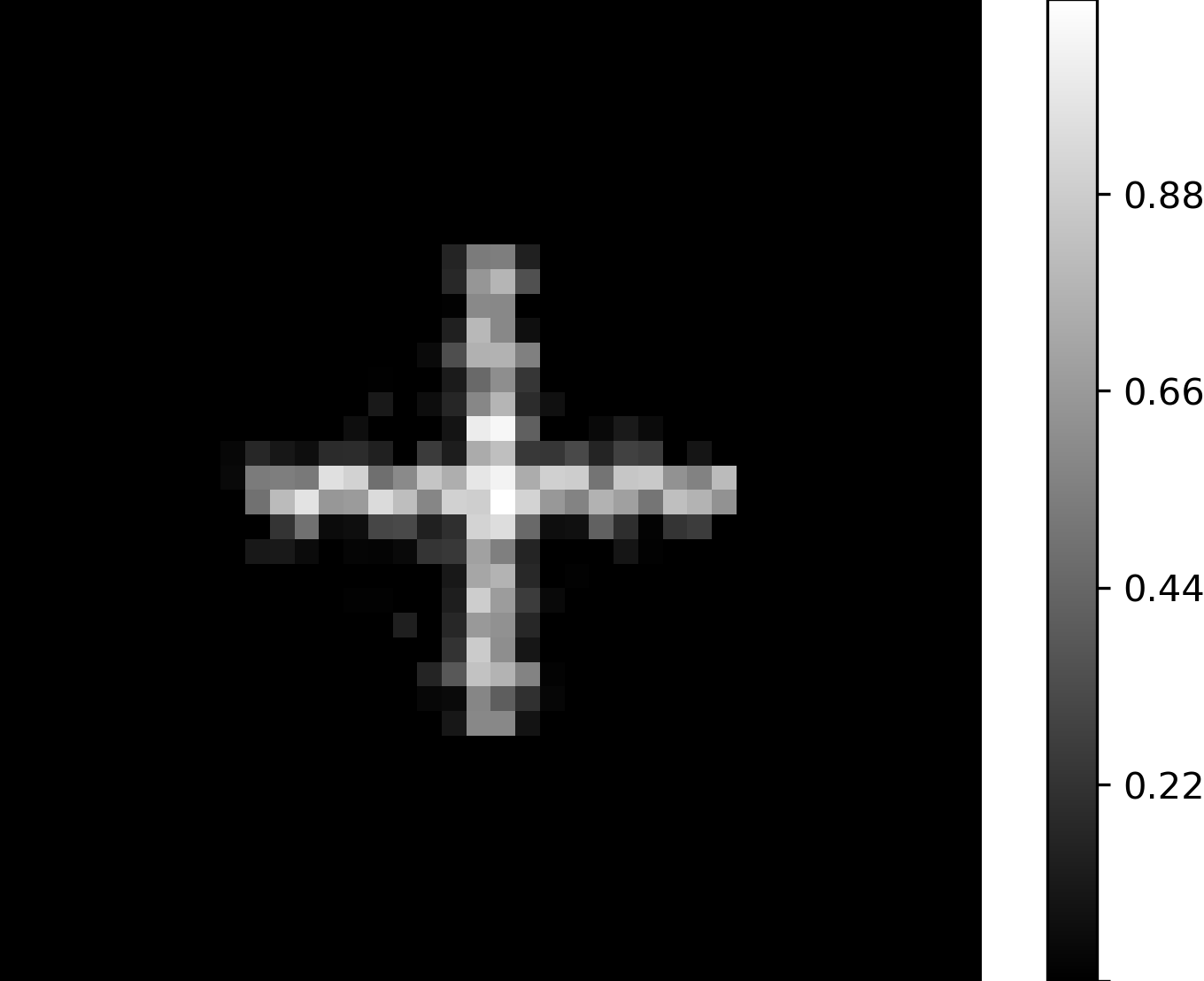}
		\caption{\centering\scriptsize  Second stage $\rho_{\text{rec}}$, PSNR 21.49, SSIM 0.8953.}
		\label{subfig:plus22:rec}
	\end{subfigure}
	
	\caption{Reconstruction of the $40\times 40$ plus-shaped phantom~\ref{subfig:plus:gt} using a $2\times 2$ multi-patch scan (Fig.~\ref{subfig:scans:22}) of the region $\Omega = [-2,2]^2$ using both the baseline system-matrix-based method and our two stage algorithm. In Fig.~\ref{subfig:plus22:SM} the reconstruction is performed using a simulated system matrix and Tikhonov regularization on each of the four disjoint regions scanned in Fig.~\ref{subfig:scans:22} and the four resulting images are stitched together. Fig.~\ref{subfig:plus22:rec} draws the final reconstruction of our 2 stage approach with the same scan data but performed in a joint fashion, taking all data points in the scanning region $\Omega$ into account. The trace field obtained in the first stage is shown in Fig.~\ref{subfig:plus22:u}. }
	\label{fig:plus22:recs}
\end{figure}

In this experiment we perform reconstructions with both our two-stage algorithm and a baseline method. The baseline method is inspired by the system-matrix-based approaches used in MPI~\cite{Ahlborg2016-ru}: we simulate a system matrix and use it to perform regularized reconstruction of the plus-shaped phantom in Fig.~\ref{subfig:plus22:gt}. We consider the $2\times 2$ multi-patch scan in Fig.~\ref{subfig:scans:22} of the region $\Omega = [-2,2]^2$ with a FoV of amplitudes $A_x = A_y =1$. In particular, let $b_1 = (-1,1)$, $b_2 = (1,1)$, $b_3 = (-1,-1)$ and $b_4 = (1,-1)$ be the four offset vectors and $\Omega_\xi\coloneq b_\xi + [-1,1]^2$ for $\xi\in\lbrace 1,2,3,4\rbrace$ the regions in which each scan is performed. The plus-shaped phantom chosen for this experiment is designed to be exactly at the center of the region $\Omega$ and to lay precisely on the contact edges of a $2\times 2$ multi-patch scans with disjoint FoVs (cf. Fig.~\ref{subfig:scans:22}). In this experiment we aim at showing the benefits of employing our method to perform reconstructions with multi-patch data when compared with methods that perform reconstructions on each region $\Omega_\xi$ and stitching of the patches after the reconstruction. More precisely, we consider a $40\times 40$ grid discretizing $\Omega$ and consequently, a $20\times 20$ grid discretization of each $\Omega_\xi$. Each subregion $\Omega_\xi$ is treated as an independent MPI reconstruction problem. Let $L$ be the number of time samples in a single scan and denote with $s^\xi=(s_k^\xi )_k\in\mathbb{R}^{2\times L}$, where $\xi$ is the quadrant index $\xi\in\lbrace 1,2,3,4\rbrace$ and $k\in\lbrace 1,\dots ,L\rbrace$ the time index, the four scans covering  the plus-shape distribution $\rho$ in $\Omega_\xi$. We want to reconstruct the $\xi$-th $20\times 20$ patch using the data $s_{k}^\xi$. To this aim, we need to simulate a system matrix for each patch. Let $\delta_{i,j}\in\mathbb{R}^{20\times 20}$ for $i,j\in\lbrace 1,\dots ,20\rbrace$ denote the the distribution $\delta_{i,j} = 1$ on the $(i,j)$-th pixel and zero otherwise. We scan each distribution $\delta_{i,j}$ and obtain the signal $s_{i,j}=(s_{i,j,k})_k\in\mathbb{R}^{2\times L}$. Reading the indexes $i,j$ in lexicographic order, we can form the matrix $S_x\in\mathbb{R}^{L\times (20\cdot 20)}$ whose  $ij$-th column is the $x$-component of the multi-patch scan $s_{i,j,k}$ of the delta distribution $\delta_{i,j}$. Analogously, a matrix $S_y\in\mathbb{R}^{L\times (20\cdot 20)}$ can be formed using the $y$-components of the the scans $s_{i,j,k}$.
Such matrices $S_x$ and $S_y$ should be in principle collected for each region $\Omega_\xi$, i.e., we obtain the matrices $S_x^\xi$ and $S_y^\xi$, and the relation between the target distribution $\rho$ and the signal $s^\xi$ is encoded in the following linear systems:
\begin{equation}
	S^\xi \rho^\xi
	= \hat{s}^\xi\quad \text{where}\quad S^\xi = \left [ (S_x^\xi )^T , (S_x^\xi )^T\right ]^T ,
\end{equation}
and $S^\xi\in\mathbb{R}^{(2\cdot L)\times (20\cdot 20)}$ is the system matrix for the patch $\Omega_\xi$, $\hat{s}^\xi\in\mathbb{R}^{2\cdot L}$ is the vector build from the scan data  $s^\xi$ by re-ordering the the vector components coherently to the component-dependent block structure of the matrix $S^\xi$ and $\rho^\xi$ is the $20\cdot 20$ vector which gives the reconstruction of the distribution $\rho$ on $\Omega_\xi$. Because we simulate the scans using ideal fields, one system matrix can be reused for each problem~\cite{szwargulski2018efficient}. This simplification leads to four problems of the form 
\begin{equation}\label{eq:sm:disjoint}
	S\rho^\xi = \hat{s}^\xi ,\quad\text{where}\quad S\coloneq S^\xi\quad\text{for all }\xi .
\end{equation}
The systems in Eq.~\eqref{eq:sm:disjoint} are then solved in a regularized fashion using Tikhonov regularization as it is typically done in MPI (see~\cite{Zdun2022multi-patch,boberggeneralized2020,szwargulski2018efficient}), i.e., we consider the following optimization problems:
\begin{equation}\label{eq:sm:minimizer}
	\rho^\xi  = \arg\min_{\hat{\rho}}\left\lVert S\hat{\rho} - \hat{s}^\xi\right\rVert^2 + \mu^\xi \lVert\hat{\rho}\rVert^2 ,
\end{equation}
and we find the minimizer by applying the CG method to the Euler-Lagrange equations associated to the minimization problems Eq.~\eqref{eq:sm:minimizer}. In this particular experiment, each problem in Eq.~\eqref{eq:sm:minimizer} has been solved for the parameter $\mu^\xi$ in the range from 35,000 and 55,000 with step 100, and the chosen $\rho^\xi$ is the solution with the highest PSNR when compared to $\rho_{\text{GT}}$ restricted to $\Omega_\xi$. For the CG method, we have considered a maximum of 10,000 iterations if the tolerance of $10^{-12}$ is not reached. Finally, the four reconstructions $\rho_\xi$ are stitched together to obtain the final reconstruction $\rho_{\text{rec}}$ in Fig.~\ref{subfig:plus22:SM} and the PSNR is computed between $\rho_{\text{rec}}$ and $\rho_{\text{GT}}$ and displayed both in Fig.~\ref{subfig:plus22:SM} and Tab.~\ref{tab:sm}.

For the two-stage algorithm, we have merged the scans $s^\xi$ in exactly the same way as in Experiment 1. Moreover, all the parameter choices for this $2\times 2$ reconstruction coincide exactly with the choices for Experiment 1. In Fig.~\ref{fig:plus22:recs} the reconstructions with both the baseline method and the two stage algorithm of this paper are shown. We point out that the same noisy data have been used for both reconstruction, but the reconstruction with the two-stage algorithm in Fig.~\ref{subfig:plus22:rec} has a higher level of noise compression and  the phantom is fairly well reconstructed using the two-stage algorithm. In particular, employing all data points at once in a joint fashion and thank to the spatial coherence of the data enforced in stage one, leads to a reconstruction which avoids the risks coming from stitching together local reconstructions, especially in the case where the particles are distributed in such a way that the region of interest can lay on the seams of the patches. Some of these effects are however due to the choice of performing separately the reconstruction of each patch. Joint reconstruction methods for system-matrix-based methods have also been developed in recent years~\cite{szwargulski2018efficient,boberggeneralized2020,Zdun2022multi-patch} as well as inversion of the systems matrices using non-negative fused lasso regularization and TV regularizer~\cite{storath2016edge}. For a fairer comparison, we describe in Experiment 8 a joint patch-wise reconstruction using the baseline method and using the algorithm developed in this paper for Stage 2 to perform a regularized inversion of the system matrix with TV and sparsity enforcing priors.

\subsection{Experiment 8: Comparison with the Baseline Approach II}\label{sec:exp:sm:lasso}

\def\imratio{0.32}
\begin{figure}[t]
	\centering
	\begin{subfigure}[t]{\imratio\linewidth}
		\includegraphics[width=\linewidth]{Images/plus_gt.png}
		\caption{\centering\scriptsize Plus-shaped phantom $40\times 40$ $\rho_{\text{GT}}$.}
		\label{subfig:plus:gt}
	\end{subfigure}
	\hfil
	\begin{subfigure}[t]{\imratio\linewidth}
		\includegraphics[width=\linewidth]{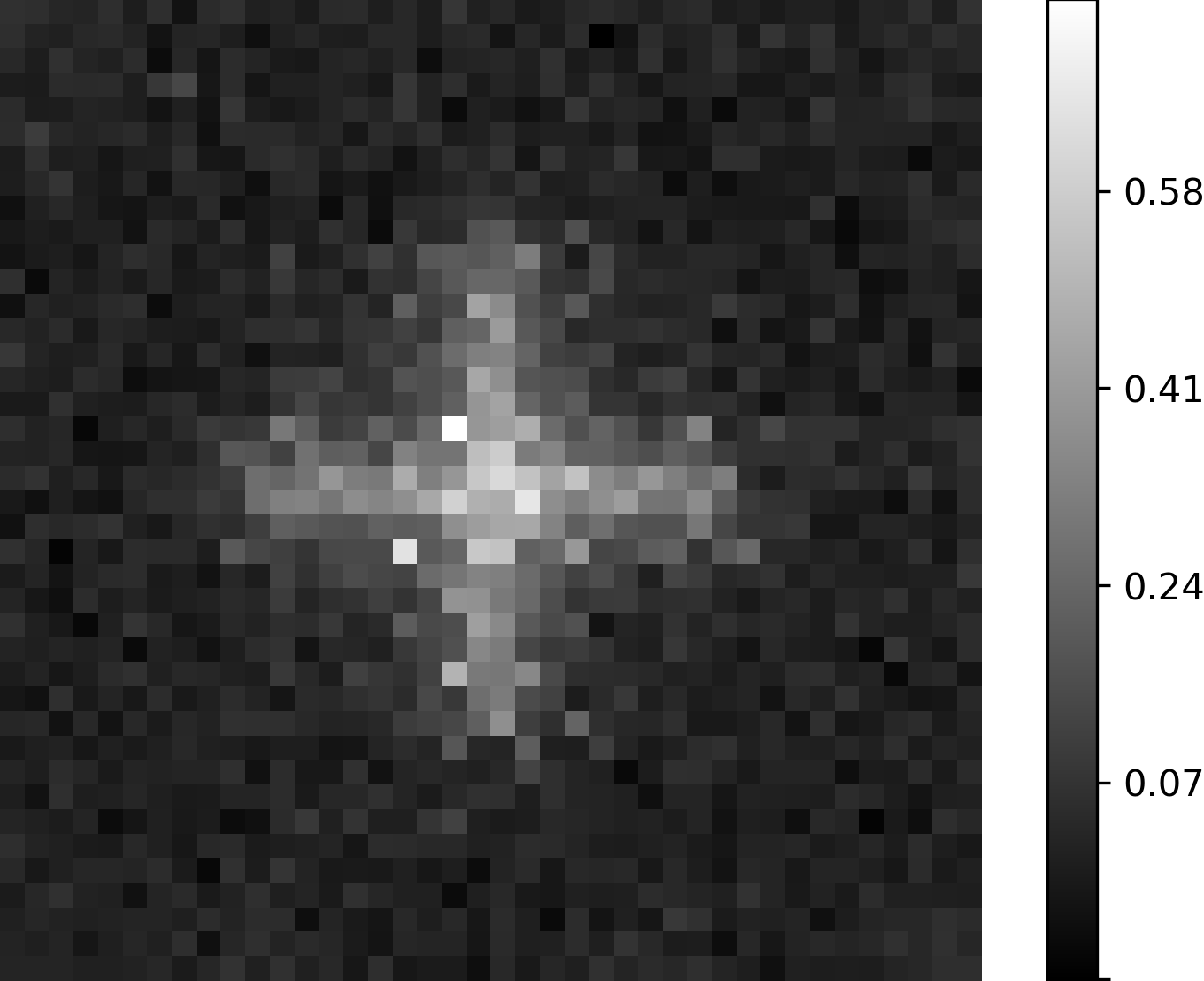}
		\caption{\centering\scriptsize $\rho_{\text{rec}}$ with SM and Tikhonov reg., PSNR 16.00, SSIM 0.2054.}
		\label{subfig:plus:SM}
	\end{subfigure}
	\hfil
	\begin{subfigure}[t]{\imratio\linewidth}
		\includegraphics[width=\linewidth]{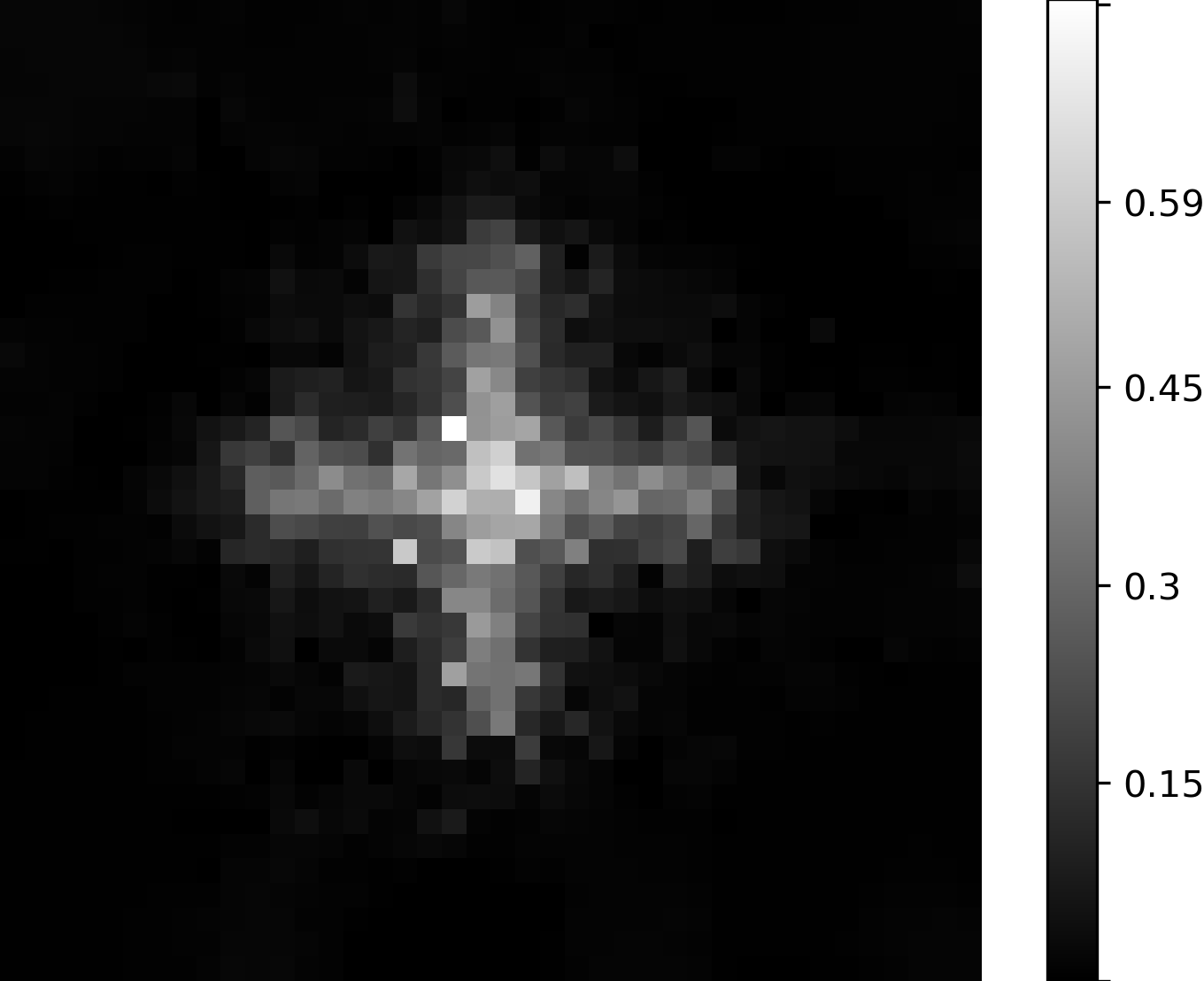}
		\caption{\centering\scriptsize $\rho_{\text{rec}}$ with SM and non-negative fused lasso, 16.38, SSIM 0.3281.}			
		\label{subfig:plus:SM:tv}
	\end{subfigure}
	\par\medskip
	\begin{subfigure}[t]{\imratio\linewidth}
		\includegraphics[width=\linewidth]{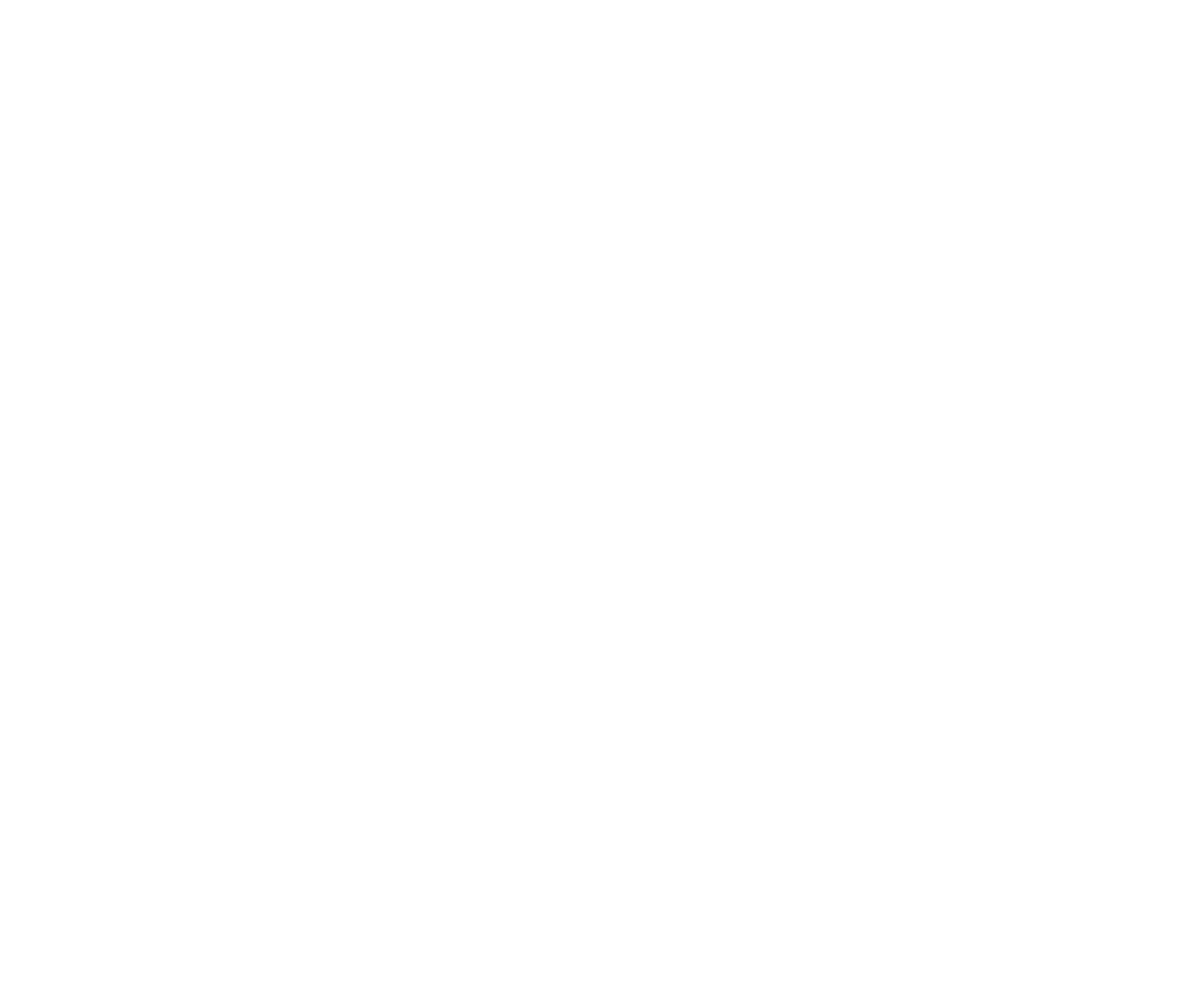}
	\end{subfigure}
	\hfil
	\begin{subfigure}[t]{\imratio\linewidth}
		\includegraphics[width=\linewidth]{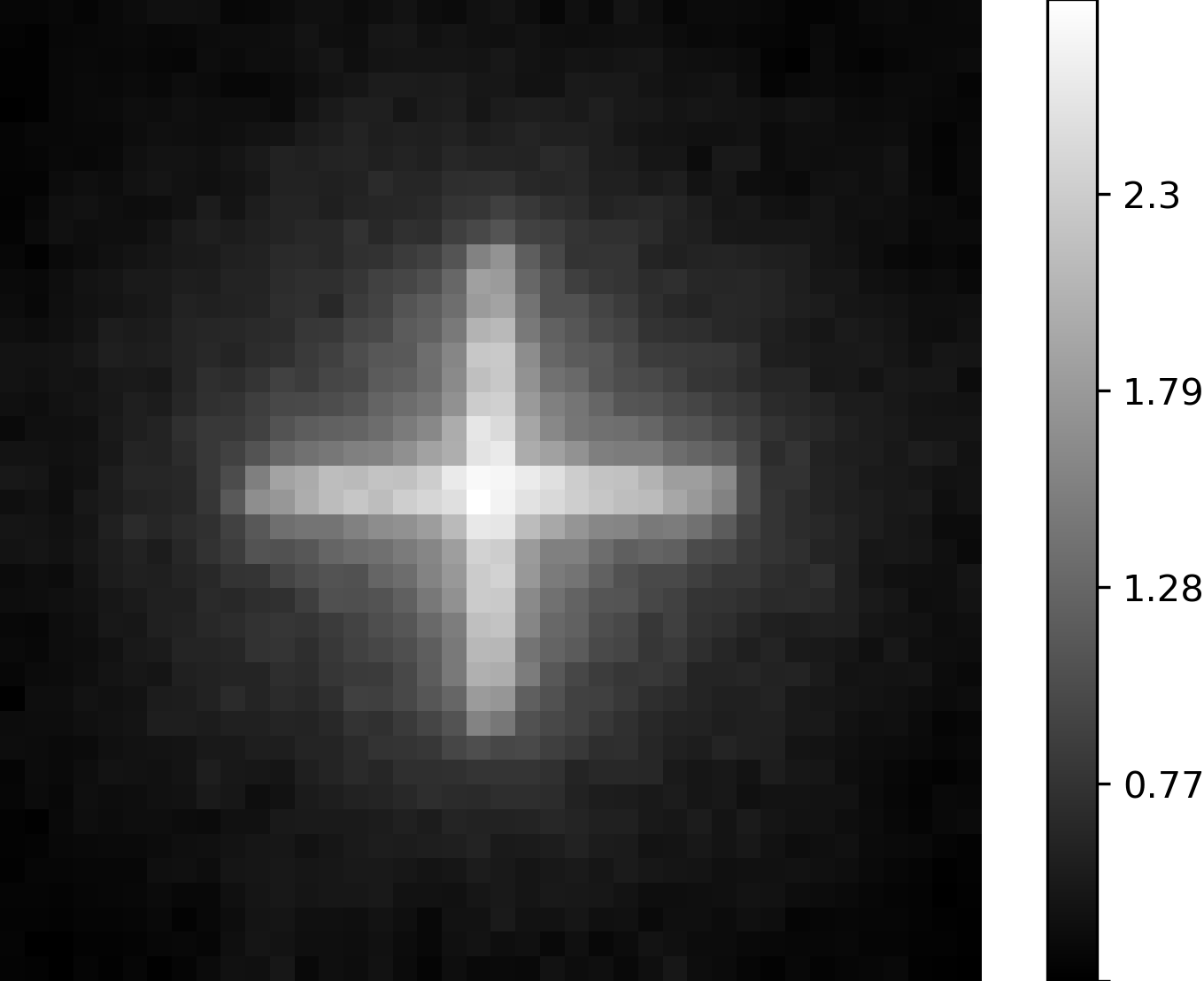}
		\caption{\centering\scriptsize First stage $u$, PSNR 33.45, SSIM 0.9467.}
		\label{subfig:plus:u}
	\end{subfigure}
	\hfil
	\begin{subfigure}[t]{\imratio\linewidth}
		\includegraphics[width=\linewidth]{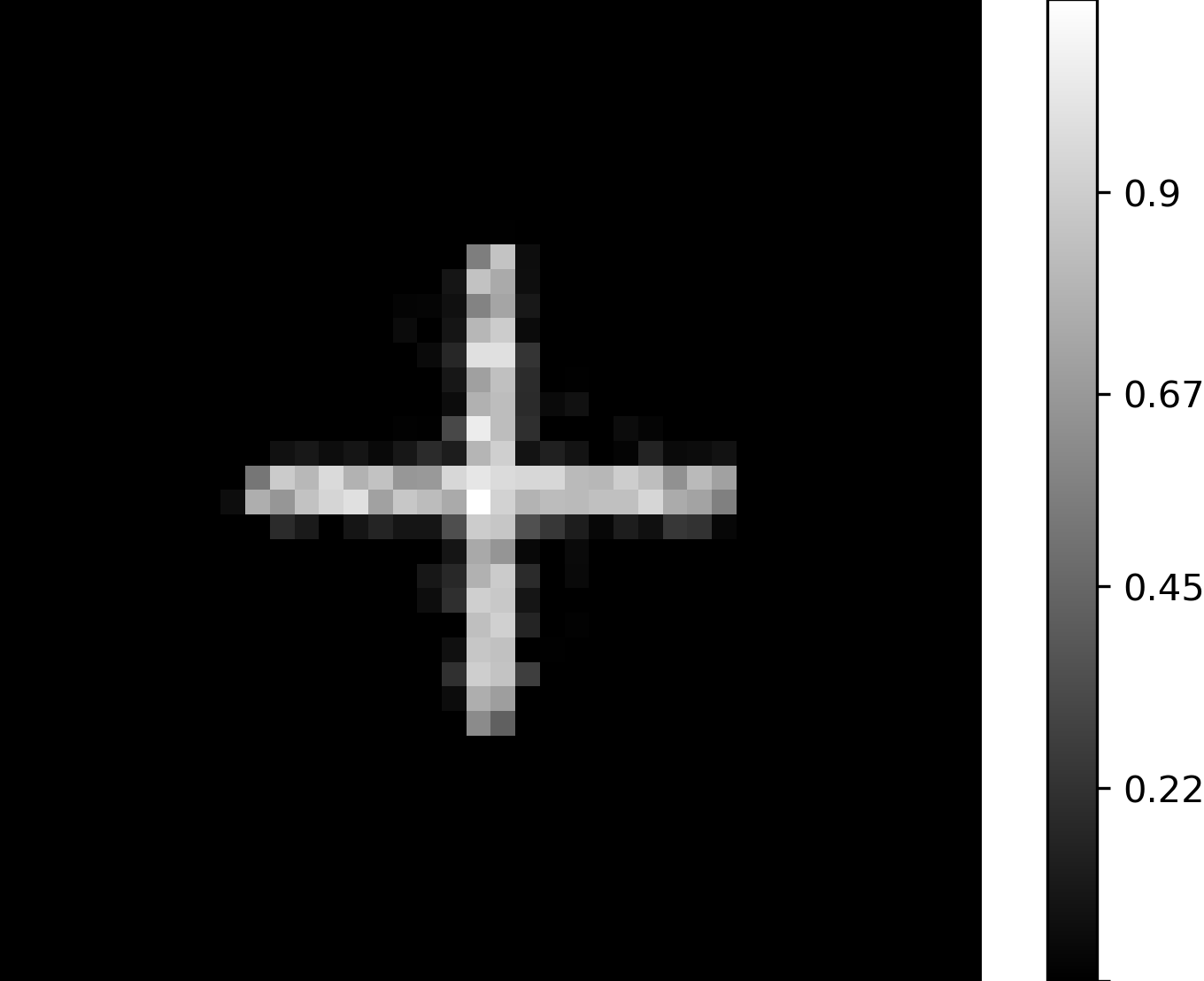}
		\caption{\centering\scriptsize Second stage $\rho_{\text{rec}}$ PSNR 24.89, SSIM 0.9375.}
		\label{subfig:plus:rec}
	\end{subfigure}
	
	\caption{Reconstruction of the $40\times 40$ plus-shaped phantom~\ref{subfig:plus:gt} using a $4\times 4$ multi-patch scan (Fig.~\ref{subfig:scans:44}) using both the baseline system-matrix-based method and our two stage algorithm. In Fig.~\ref{subfig:plus:SM} the reconstruction is performed using a simulated system matrix and Tikhonov regularization, while in Fig.~\ref{subfig:plus:SM:tv} the same system matrix is used to perform a reconstruction with the TV-Smooth regularizer and a sparsity enforcing prior. In the second row the results of the two stage algorithm are shown, in particular the trace field $u$ in Fig.~\ref{subfig:plus:u} outputted by the first stage is used as the input for the second stage, which yields the final reconstruction in Fig.~\ref{subfig:plus:rec}. }
	\label{fig:plus:recs}
\end{figure}

In this experiment we use again the baseline method described in Experiment 7, but to have a fairer comparison with the results of the two-stage algorithm we perform a joint-reconstruction for both methods inspired by~\cite{Knopp_2015,szwargulski2018efficient,boberggeneralized2020,Zdun2022multi-patch} and in particular employ the non-negative fused lasso regularization of Section~\ref{sec:rec:algorithm} for the inversion of the system matrix. 

We consider the same phantom and the same $40\times 40$ discretization of the region $\Omega = [-2,2]^2$ as in Experiment 7, but we perform a $I\times J = 4\times 4$ multi-patch scan as in Fig.~\ref{subfig:scans:44} with $L=1632$ time samples per scan. With this multi-patch scan the patches overlap and methods to perform joint reconstruction of the patches have been developed~\cite{szwargulski2018efficient,boberggeneralized2020,Zdun2022multi-patch}. The idea behind it is to collect the system matrices associated to each patch (or re-use the same system matrix) and construct a bigger block matrix, whose blocks are the system matrices of each patch composed with an index operator that refers the correct equation of the system to the corresponding pixel of the target unknown variable $\rho$.  Because our focus in this experiment is to analyze the quality of the reconstruction and not to optimize the time-efficiency of the reconstruction, we simulate a system-matrix for the multi-patch scan. More explicitly, given the $4\times 4$ multi-patch scan approach in Fig.~\ref{subfig:scans:44} and described in Experiment 1, we perform  a $4\times 4$ multi-patch scan of the phantom and of each of the delta distributions $\delta_{i,j}$ for $i,j\in\lbrace 1,\dots ,40\rbrace $ where $\delta_{i,j} =1$ on the $(i,j)$-th pixel and zero otherwise. The multi-patch scans of the $\delta_{i,j}$ distributions are collected in the columns of a system matrix $S\in\mathbb{R}^{(2\cdot L\cdot I\cdot J)\times (40^2) }$ and the scan of the plus-shaped phantom $\rho_{\text{GT}}$ in Fig.~\ref{subfig:plus:gt} is collected in the vector $s\in\mathbb{R}^{2\cdot L\cdot I\cdot J}$. The baseline method consists therefore in solving in a regularized fashion the linear system of equations $S\hat{\rho} = s$, where $\hat{\rho}\in\mathbb{R}^{40^2}$ is the unknown variable in vector form. We first perform regularized inversion of the system using Tikhonov regularization, i.e., we solve the following minimization problem:
\begin{equation}\label{eq:sm2:mininimizer}
	\rho  = \arg\min_{\hat{\rho}}\left\lVert S\hat{\rho} - s\right\rVert^2 + \mu \lVert\hat{\rho}\rVert^2 ,
\end{equation}
which is solved applying the CG method to the Euler-Lagrange equations associated to the problem in Eq.~\eqref{eq:sm2:mininimizer}. Similarly to Experiment 7, we perform the reconstruction for $\mu$ ranging between 10,000 and 80,000 with step 1, and with a maximum number of CG iteration of 10,000 if a tolerance of $10^{-12}$ is not reached. The result of the Tikhonov regularization are displayed in Fig.~\ref{subfig:plus:SM} and we observe an improvement when compared to the reconstruction in Fig.~\ref{subfig:plus22:SM}. This result is not surprising in view of the fact that the plus-shaped phantom is positioned in the region of $\Omega$ where the $4\times 4$ scans (cf. Fig.~\ref{subfig:scans:44}) overlap, producing a much denser sampling of the phantom. Moreover, the system matrix collected here contains the system response to delta impulses of the multi-patch scans, meaning that the information in the regions where the patches overlap is naturally contained in $S$ and the loss of information due to considering each patch separately is minimized.

Given the linear system $S\rho = s$, we would also perform regularized inversion with the TV-Smooth regularizer and a sparsity enforcing prior (this is inspired by similar results in~\cite{storath2016edge,Zdun2022multi-patch}). This regularized inversion can be then formulated as a minimization problem of the following form:
\begin{align}\label{eq:sm:lasso:func}
	\rho = \arg\min_{\hat{\rho}} \biggl\lbrace \lVert S \hat{\rho}-s\rVert_2^2 + \mu R_\delta [\hat{\rho} ] + \beta\lVert\hat{\rho}\rVert_1 +  \iota_+ (\hat{\rho}) \biggr\rbrace ,
\end{align}
where $R_\delta [\rho ]$ is the TV-Smooth approximation of TV as in Eq.~\eqref{eq:tv:smooth:discrete}. This minimization problem has been solved using the algorithm in Alg.~\ref{alg:lasso} where we have replaced the convolution matrix $K_h$ with the system matrix $S$ and $u$ with $s$, whose convergence is guaranteed by Theorem~\ref{thm:convergence} where a Lipschitz constant is of the form (cf. Eq.~\eqref{eq:lipschitz:constant})
\begin{equation}
	\tilde{L}\coloneq \mu L + \lVert S^T S \rVert_{2} .
\end{equation}
For the solution of the problem in Eq.~\eqref{eq:sm:lasso:func} we use Alg.~\ref{alg:lasso} but with different data fidelity with $\delta = 10^{-16}$, $\alpha = 1$, $\gamma = 10^{-9}$ and $20,000$ maximum iterations if a tolerance of $10^{-5}$ is not reached. The parameter $\mu$ ranges between 10 and 600 with step 10. The solution given by the minimization problem in Eq.~\eqref{eq:sm:lasso:func} is shown in Fig.~\ref{subfig:plus:SM:tv}, where the effect of the sparsity enforcing prior as well as of the positivity constraint are clear.

Finally, in Fig.~\ref{subfig:plus:u} and \ref{subfig:plus:rec} the trace field $u$ after the first stage and the final reconstruction after the second stage $\rho_{\text{rec}}$ are displayed.

\subsection{Discussion on the Regularization Parameters}\label{sec:disc:reg:params}
The results of all the experiments with the two stage algorithm are collected in Tab.~\ref{tab:exp} and provide us with the following observations on the regularization parameters. Concerning the parameter $\lambda$ in the first stage of the algorithm, it regulates the effect of the regularizer in Eq.~\eqref{eq:first:discrete:reg}. This regularizer is in charge of imposing smoothness of the solution and of the inpainting, i.e., of providing with reasonable values those regions where there is less (or lack of) data. Consequently, it is reasonable to suppose that $\lambda$ strongly depends on the amount of data points in each pixel, given a discretization grid: indeed, if the grid is finer than the average distance between the data points, numerous sets of pixel do not contain data and in turn, the inpainting effect must be stronger to assign suitable values to dataless pixel. One way to test this, is to fix the discretization grid and observe the behavior of $\lambda$ by increasing the amount and density of the data points. Experiment 1 provides such a scenario and in Tab.~\ref{tab:exp} can be observed that, given the same phantom, $\lambda$ decreases if the number of data points is increased considering the merged scanning data coming from $2\times 2,\cdots, 10\times 10$ patches. Moreover, letting the phantom vary and keeping both the discretization grid and the number of data points ($10\times 10$ multi-patch scans) fixed, the parameter $\lambda$ seems to be stable most of the times. The outlier $\lambda =21$ for the shape phantom could suggest a certain degree of dependency on the phantom, suggesting that methods for the optimal selection of the parameter in absence of the ground truth (as in the real scans) could involve optimization on large datasets of scans, containing various multi-patch modalities on a plethora of different phantoms. As a further confirmation of our hypotheses, we observe that $\lambda$ is the smallest in Experiment 6, where the reconstruction grid is coarser ($100\times 100$) and the amount of the data is much larger than in the previous experiments, leading to a much denser data coverage of the region $\Omega$.

Interestingly, the parameter $\mu$ in charge of the regularization in the deconvolution in the second stage, seems to be stable and does not vary much among Experiments 1 to 5 (cf. Tab.~\ref{tab:exp}). These experiments are all performed on a $200\times 200$ grid and the input of the second stage, independently of the phantom and of the specific output of the first stage, is always a $200\times 200$ array. Indeed, it is the first stage that employs the data and heavily depends on it, while the deconvolution in Experiments 1 to 5 are performed with inputs of the same size. This suggests that the grid size is one of the main variables on which the parameter $\mu$ depends. Indeed, the only outlier is in Experiment 6, where $\mu = 10^{-13}$, and the reconstruction is performed on a $100\times 100$ grid. Again, if no ground truth is available, one way of fine-tuning $\mu$ could involve the optimization on a data set of scans for different choices of the discretization grid.

Finally, the parameter $\beta$ has been fixed by convenience to $\beta =1$ for all experiments but the concentration phantom and Experiment 6. The parameter $\beta$ is in charge of the degree of minimization of the $L_1$-norm of the distribution (see Eq.~\ref{eq:second:pos:sparsity:discrete}) and in the discrete scenario, clearly depends on the discretization grid used to compute the norm. This motivates why $\beta$ has to be lowered to $0.1$ in Experiment 6, which performs the reconstruction on the coarser $100\times 100$ grid. Concerning the concentration phantom in Experiment 4, the parameter $\beta$ has again been lowered because the promotion of too much sparsity in the solution, erased the area with very low concentration of particles. This means that the parameter $\beta$ depends on the concentration levels of the distribution and can be potentially fine tuned using a data set containing scans performed on phantoms with different concentration levels.

\section{Discussion and Conclusion}\label{sec:conclusions}

\begin{table}[h]
	\caption{Summary of the results of the experiments.}\label{tab:exp}
	\newcolumntype{C}{>{\centering\arraybackslash}X}
	\begin{tabularx}{\linewidth}{CCCCCCCC}
		\toprule
		\textbf{Phantom} & \textbf{Scan type} & \multicolumn{3}{c}{\textbf{First stage}} & \multicolumn{3}{c}{\textbf{Second stage}} \\
		&   & \textbf{Reg. Param. $\lambda$} & \textbf{PSNR} & \textbf{SSIM} & \textbf{Reg. Param. $\mu$} & \textbf{PSNR} & \textbf{SSIM} \\
		\midrule
		Vessel Fig.~\ref{subfig:exp:vess:multi:gt} & $2\times 2$  & 16 & 26.21 & 0.7853 & $7.5\cdot 10^{-5}$ & 10.14 & 0.3316\\
		\cmidrule{2-8}
		&  $4\times 4$  & 12 & 28.99 & 0.8307 & $3\cdot 10^{-5}$ & 11.36 & 0.4337\\
		\cmidrule{2-8}
		& $6\times 6$  & 8 & 30.10 & 0.8521 & $10^{-4}$ & 12.29 & 0.5247 \\
		\cmidrule{2-8}
		& $8\times 8$  & 6 & 31.51 & 0.8703 &$10^{-4}$ & 12.92  & 0.5967 \\
		\cmidrule{2-8}
		&  $10\times 10$  & 5 & 32.00 &0.8857 & $10^{-4}$ & 13.41  & 0.6038\\
		\cmidrule{2-8}
		& 143 Random  & 8 & 34.53 &0.9238 & $10^{-4}$ & 13.33 & 0.5903 \\
		\cmidrule{2-8}
		& $10\times 10$ Perturbed  & 5 & 32.35 &0.8864 & $10^{-4}$ & 13.23  & 0.5963\\
		\cmidrule{2-8}
		& $10\times 10$ Fig.~\ref{subfig:abl:u}, Fig.~\ref{subfig:abl:rho:land} & 5 & 32.00 &0.8857 & $1.75\cdot 10^{-4}$ & 12.73 & 0.2811 \\
		\midrule
		Shape Fig.~\ref{subfig:shape:gt} & $10\times 10$ & 21 & 42.97 &0.9863 &  $2.5\cdot 10^{-4}$ & 26.86 & 0.9860 \\
		\midrule
		Conc. Fig.~\ref{subfig:concentration:gt} & $10\times 10$ & 6 & 39.41 & 0.9588 & $2.5\cdot 10^{-5}$ & 29.75 & 0.9743 \\
		\midrule
		Frame Fig.~\ref{subfig:frame:gt} & $10\times 10$ & 7 & 33.88 &0.9262 & $5\cdot 10^{-5}$ & 20.29 & 0.9058 \\
		\cmidrule{2-8}
		& $10\times 10$ Fig. \ref{subfig:frame:u:p},Fig. \ref{subfig:frame:rho:p} & 6 & 33.10 &0.9181 & $5\cdot 10^{-5}$ & 19.97  & 0.9002\\
		\cmidrule{2-8}
		& $10\times 10$ Fig. \ref{subfig:frame:u:p2},Fig. \ref{subfig:frame:rho:p2} & 8 & 26.17 & 0.8534 & $5\cdot 10^{-5}$ & 15.63  & 0.7351\\
		\midrule
		Vessel Fig.~\ref{subfig:1000:gt}& Generalized Fig. \ref{fig:1000:scans} & 1 & 37.68 &0.9700 & $10^{-13}$ & 12.81  & 0.5265\\
		\midrule
		Plus Fig.~\ref{subfig:plus22:gt} & $2\times 2$  & 7 & 30.01 & 0.9092 & $10^{-13}$ & 21.49 & 0.8953 \\
		\cmidrule{2-8}
		& $4\times 4$ & 14 & 33.45 & 0.9467 & $10^{-13}$ & 24.89 & 0.9375\\
		
		\bottomrule
	\end{tabularx}
\end{table}

\begin{table}[h]
	\caption{Summary of the results of experiments with a simulated system matrix (SM).}\label{tab:sm}
	\newcolumntype{C}{>{\centering\arraybackslash}X}
	\begin{tabularx}{\linewidth}{CCCCCC}
		\toprule
		\textbf{Phantom} & \textbf{Scan type} & \textbf{Method} & \textbf{Reg. parameter } & \textbf{PSNR} & \textbf{SSIM}  \\
		Plus-shaped phantom Fig.~\ref{subfig:plus22:gt} & $2\times 2$ & SM + Tikhonov & $\mu^1 =39,700$,  $\mu^2 =52,600$,  $\mu^3 =49,400$,  $\mu^4 =45,700$. & 14.64 & 0.0572 \\
		\cmidrule{2-6}
		& $4\times 4$ & SM + Tikhonov & $72,770$ & 16.00 & 0.2054 \\
		\cmidrule{3-6}
		& & SM + non-neg. Lasso ($\alpha =1$) & 380 & 16.38 & 0.3281 \\

		\bottomrule
	\end{tabularx}
\end{table}

\subsection{Discussion}

The application of magnetic fields on humans exposes them to the risk of tissue heating, unless the strength of the magnetic fields employed are limited~\cite{Saritas2013-jn}. The limitation due to safety reasons on the strength of the magnetic fields applied translates into a limitation of the size of the FoV. Multi-patching is therefore necessary to scan regions that exceed in size the covering capacity of the FoV. In the system-matrix-based approaches, methods exists to perform such multi-patch reconstructions. The first idea was to perform patch-wise reconstructions and combine the patches in a post-processing step~\cite{Ahlborg2016-ru}. Then it has been shown that it is advantageous to combine the measured data into one linear system of equations, thus formulating a joint system matrix whose inversion automatically takes care of the overlapping boundaries or prevents the artifacts coming from the stitching~\cite{Knopp_2015}. Interestingly, it has also been shown that it is possible to perform multi-patch reconstructions re-using the same system matrix~\cite{szwargulski2018efficient}, which allows to reduce the calibration time. The idea of reusing the same system matrix is based on the fact that for ideal magnetic fields the system function in MPI is shift-invariant. This holds approximately true in the central region of the scanner bore, but not when the FoV is further away from the center~\cite{szwargulski2018efficient}. For real data effort has been put in the reduction of the number of system matrices needed~\cite{boberggeneralized2020}. These steps have been simulated in Experiment 7 and 8 of this paper and considered as a baseline reference method. In this paper we propose a mathematical framework that can potentially deal with multi-patch scans without the need of a system matrix, and could potentially give more freedom to the ways the patches are collected, as they do not require any a priori grid-coherence of the patches: it is enough to collect the data samples and their positions (and velocities), the algorithm takes care of reconstructing the distribution on any chosen grid. The choice of the grid can be made after the collection of the data samples. 

Moreover, we use the principle of an MPI Core Operator to get rid of the dependency of the trajectory employed and treat different scan as set of data points that can be merged.
This is particularly useful because it allows to formulate the general multi-patching approach in Section~\ref{sec:mathematical:modeling} and for example consider modified trajectories that take into account the inherent distortion of the real scanning paths in MPI scanners:
if a particular distortion is known, it is enough to substitute the position points and velocities with the distortion-adjusted positions
and velocities without changing anything of the reconstruction algorithm but the input positions.
Such distortions can depend, for example, on the employment of additional focus fields and the position of the FoV inside the scanner
(in particular on the proximity to the boundary of the scanning region~\cite{szwargulski2018efficient}).

We next discuss the mathematical model employed.
Mathematical models of the physics behind the scanner are important because they allow for simulation
and thus for prototyping of novel MPI systems and for testing of reconstruction techniques.
Further, they may form the basis for developing analytic reconstruction techniques (as discussed above.)
Presently, models based on the Langevin theory of paramagnetism are well established and widely used in the
context of MPI~\cite{Knopp2017proofprinciple,Gdaniec2020supressionmotionartifacts,Shang_2022,soydan20213Dslicewise}.
The Langevin model is a mathematical model of the magnetization response
of superparamagnetic nanoparticles. It is based on ideal behavior of the nanoparticles.
Also the models and reconstruction schemes in this paper are based on the established Langevin theory,
and as such they inherit the limitations coming from the employment of the Langevin model.
Improving the modeling, in particular to narrow the gap between simulations based on the Langevin model
and the real scans is an ongoing topic of research in MPI; e.g.~\cite{bringout2020new,Kluth_2019_FokkerPlanck}.
For instance, to describe the magnetic response of the nanoparticles, one can solve the Fokker-Planck equation
for a probability distribution representing the particle distribution to improve the  approximation fo the mean
magnetic moment of the particles~\cite{Kluth_2019_FokkerPlanck}.
With regard to the modeling of magnetic fields, more sophisticated approaches exist and show that the field-free region in the FFP topology
is not a point but is a low field volume~\cite{bringout2020new}.

We discuss now the use of simulated data.
The generation of simulated data is a very important topic in MPI since presently
the only initiative for freely accessible MPI data is the Open MPI Dataset~\cite{knopp2020openmpidata}.
It contains measurements on three different phantoms with different FFP scanning sequences (1D, 2D and 3D) and dedicated system matrices for each scanning sequence performed using a real MPI scanner. If the data from the Open MPI Dataset does not fit with the topic considered in their paper, authors frequently refer back to their own simulated data~\cite{bringout2020new,Yagiz2020nonideal,Zdun2022multi-patch}. In this respect, we note that,
up to date, there is no common platform for simulation of MPI data which leads to the fact that most works employ their own simulation~\cite{Yagiz2020nonideal,Zdun2022multi-patch}. The topic of the present paper -- a mathematical framework for multi-patching -- does not fit with the Open MPI Dataset, as it does not contain data of multi-patch scans. Therefore, in this paper we have employed data simulated as described in Section~\ref{sec:experiments}, which is similar to~\cite{Brandt2021-rw}. An important topic of future research is the validation of the approach on real data. 

A final point of discussion is the choice of the regularization parameters in the algorithm if the ground truth is not know, for example when working with real data.
The selection of optimal parameters constitutes an interesting direction of future research.
The insights in the discussion on the regularization parameters at the end of Section~\ref{sec:experiments} show that the parameter $\lambda$
in the first stage depends on the number of data points acquired and that the parameter $\mu$ in the second stage seems stable if the reconstruction grid is kept fixed.
This suggests that optimal parameter selection strategies could be developed if a larger quantity of scan data is available.

\subsection{Conclusion}

In this paper we have introduced a mathematical framework that allows to describe general scanning sequences for multi-patching in MPI,
in particular our framework allows to describe any multi-patching modality that involves rigid movement of the MPI scanner,
the FoV and the specimen inside the scanner. With this framework we can describe the currently employed multi-patching modalities
which either use a focus field to move the FoV or physically move the specimen.
We have proved the flexibility of the mathematical framework showing how it can not only deal with possible multi-patching approaches that rely on moving the specimen
and the FoV at the same time, but also on possible future multi-patching approaches in which the scanner is allowed to move as well.
We have moreover provided formulas to transform the data acquired by the scanner into input data for our two-stage reconstruction algorithm~\cite{math10183278}.
As a second contribution of the paper, we have shown how multi-patching can be incorporated into the two-stage algorithm and added positivity constraints
and sparsity promoting priors for better quality reconstruction, proving that under reasonable hypotheses the provided algorithm converges to a minimizer
of the related energy functional. The last contribution of the paper consists of numerical experiments that simulate the scan and reconstruction
with the proposed method with a variety of different phantoms and multi-patching modalities,
showing the potential of the proposed methods applied to multi-patching with a scanning setup which resembles the one employed in the Open MPI project~\cite{knopp2020openmpidata}.

Lines of future research include improving the proposed algorithmic scheme
by employing different regularizers such as machine learning-based regularizers and Potts type regularizers~\cite{kiefer2021iterative,weinmann2013l1potts,storath2015joint} in Stage 2,
and improvements of the inpainting scheme in Stage 1.
Further, the investigation on automatic parameter choice strategies
as well as addressing the gap between simulated models and the real data are important topics of future research.

\paragraph{\textbf{Declarations of interest}} None.

\paragraph{\textbf{Acknowledgements}} This research was supported by the Hessian Ministry of Higher Education, Research, Science and the Arts within the Framework of the "Programm zum Aufbau eines akademischen Mittelbaus an hessischen Hochschulen". We would like to thank the DEAL and Elsevier Agreement for supporting the open access publication of this paper.

\bibliographystyle{elsarticle-num}

\bibliographystyle{siam}
\bibliography{References}

\end{document}